\documentclass[11pt,reqno]{amsart}

\usepackage[T1]{fontenc}
\usepackage{lmodern}
\usepackage{microtype}
\usepackage{amsmath,amssymb,amsthm,mathtools}
\usepackage{enumitem}
\usepackage{xcolor}
\usepackage{aliascnt}
\usepackage{float}
\usepackage{needspace}
\usepackage{placeins}
\usepackage{flafter}
\usepackage{booktabs,tabularx}
\usepackage{tikz}
\usepackage{mathrsfs}
\usetikzlibrary{arrows.meta,positioning,calc,fit}
\usepackage[colorlinks=true,linkcolor=blue!45!black,
 citecolor=blue!45!black,urlcolor=blue!55!black]{hyperref}
\usepackage[nameinlink,capitalise]{cleveref}

\allowdisplaybreaks
\setlist{itemsep=2pt,topsep=4pt}

\newtheorem{theorem}{Theorem}
\newtheorem{proposition}{Proposition}[section]
\newaliascnt{sectiontheorem}{proposition}
\newtheorem{sectiontheorem}[sectiontheorem]{Theorem}
\aliascntresetthe{sectiontheorem}
\newaliascnt{lemma}{proposition}
\newtheorem{lemma}[lemma]{Lemma}
\aliascntresetthe{lemma}
\newaliascnt{corollary}{proposition}
\newtheorem{corollary}[corollary]{Corollary}
\aliascntresetthe{corollary}
\theoremstyle{definition}
\newaliascnt{definition}{proposition}
\newtheorem{definition}[definition]{Definition}
\newtheorem*{openquestion}{Question}
\aliascntresetthe{definition}
\theoremstyle{remark}
\newaliascnt{remark}{proposition}
\newtheorem{remark}[remark]{Remark}
\aliascntresetthe{remark}
\newcommand{\D}{\mathbb D}
\newcommand{\T}{\mathbb T}
\newcommand{\R}{\mathbb R}
\newcommand{\Z}{\mathbb Z}
\newcommand{\Q}{\mathbb Q}
\newcommand{\N}{\mathbb N}
\newcommand{\id}{\operatorname{id}}
\newcommand{\Diff}{\operatorname{Diff}}
\newcommand{\dist}{\operatorname{dist}}
\newcommand{\Lip}{\operatorname{Lip}}
\DeclareMathOperator{\cond}{cond}
\newcommand{\rot}{\operatorname{rot}}

\newcommand{\cC}{\mathcal C}
\newcommand{\tcC}{\widetilde{\mathcal C}}
\newcommand{\normT}[1]{\lVert #1\rVert_{\T}}

\newcommand{\Per}{\operatorname{Per}}
\newcommand{\supp}{\operatorname{supp}}
\newcommand{\Bclass}{\mathcal B}
\newcommand{\Bclassbeta}[1]{\mathcal B_{#1}}

\crefname{part}{Part}{Parts}
\crefname{sectiontheorem}{Theorem}{Theorems}
\crefname{lemma}{Lemma}{Lemmas}
\crefname{proposition}{Proposition}{Propositions}
\crefname{corollary}{Corollary}{Corollaries}
\crefname{definition}{Definition}{Definitions}
\crefname{remark}{Remark}{Remarks}
\definecolor{deepblue}{RGB}{31,65,105}
\definecolor{softgray}{RGB}{245,247,249}
\definecolor{softblue}{RGB}{238,244,251}
\definecolor{coregreen}{RGB}{219,239,222}
\definecolor{coregreenline}{RGB}{42,122,69}
\definecolor{frontorange}{RGB}{255,229,194}
\definecolor{frontline}{RGB}{191,98,24}
\definecolor{collarviolet}{RGB}{235,228,247}

\tikzset{
  nipl/.style={line width=.6pt,draw=black!70},
  nipa/.style={-{Latex[length=2mm,width=1.3mm]},line width=.6pt,draw=black!70},
  nipbox/.style={draw=black!65,fill=gray!7,rounded corners=1.5pt,
    inner sep=4pt,align=center}
}

\title[Rigidity and non-rigidity below $C^2$]
{Rigidity and Non-Rigidity of Disk Pseudo-Rotations below $C^2$}
\author{JIAN WANG}
\address{School of Mathematical Sciences and LPMC, Nankai University,
Tianjin 300071, PR China}
\email{wangjian@nankai.edu.cn}
\thanks{Jian Wang was partially supported by NSFC (No.~12361141812)
and the Fundamental Research Funds for the Central University (No.~63213122).}
\author{PEIZHENG YU}
\address{School of Mathematics, Shandong University,
Jinan 250100, PR China}
\email{yupzh@sdu.edu.cn}
\date{}
\hypersetup{
pdftitle={Rigidity and Non-Rigidity of Disk Pseudo-Rotations below C2},
pdfauthor={Jian Wang and Peizheng Yu},
pdfsubject={Sequel to arXiv:2610.09986v1: non-Brjuno rigidity, C1 counterexamples and bi-Lipschitz endpoint examples}}

\begin{document}
\begin{abstract}
For every $0<\beta\le1$, we prove $C^0$-rigidity of area-preserving
$C^{1,\beta}$ disk pseudo-rotations with non-Brjuno rotation number.
At the \(C^1\) endpoint, every non-Brjuno rotation number of zero
exponential type admits non-rigid examples. We also obtain a different rigidity result when the rotation number has positive exponential type: the boundary restriction of every $C^1$ pseudo-rotation is rigid and topologically conjugate to a rotation. We show that this boundary conclusion fails in the bi-Lipschitz class by constructing non-rigid examples for every rotation number of finite positive
exponential type. The proof of rigidity is based on a finite-information closing theorem under Hölder regularity. For the construction of counterexamples, we use disjoint supports instead of a summable perturbation scheme.
\end{abstract}
\maketitle
\enlargethispage{2pt}
\raggedbottom
\setcounter{tocdepth}{2}
\tableofcontents

\section{Introduction}\label{sec:introduction}

A disk pseudo-rotation has the simplest possible periodic-point set (only one periodic point), yet its iterates need not return uniformly to the identity. The rigidity problem is to understand how the regularity of the map and the arithmetic properties of its rotation number affect this behavior.

This paper continues Liu--Wang--Yu \cite{LWY}, who constructed
non-rigid $C^1$ disk pseudo-rotations for every Brjuno rotation number
and, for $0<\beta<1$, $C^{1,\beta}$ examples in a weighted Brjuno subclass.  In
\S1.3 of that work, they asked whether every area-preserving $C^1$
disk pseudo-rotation with a non-Brjuno rotation number of finite
exponential type must be $C^0$-rigid.

The positive results build on Bramham's rigidity theorem under a
strong Liouville condition \cite{Bramham} and the non-Brjuno
$C^{k-1}$-rigidity theorem for $C^k$ maps, $k\ge2$, of
Avila--Fayad--Le Calvez--Xu--Zhang \cite[Theorem~1]{AFLXZ}.
The regularity question of whether every $C^k$
pseudo-rotation is $C^0$-rigid was recorded in
\cite[Question~3]{AFLXZ} and \cite[Question~13]{FayadKrikorian},
following Bramham's smooth formulation \cite[Question~2]{Bramham}.
Our results identify a sharp $C^1$ endpoint for non-Brjuno rigidity. We obtain new results at the \(C^1\) and \(C^{1,\beta}\) levels. We also discuss two related questions that remain open.

\subsection{Conventions and arithmetic}\label{sec:conventions}
Let $\D=\{z\in\R^2:|z|\le1\}$ be a disk with standard area, let
$\T=\R/\Z$, and write
\[
 R_\alpha(x)=x+\alpha\quad(x\in\T),\qquad
 \mathcal R_\alpha(z)=e^{2\pi i\alpha}z\quad(z\in\D).
\]
An \emph{irrational disk pseudo-rotation} is an orientation- and
area-preserving homeomorphism $F$ with
\[
 \Per(F):=\bigcup_{m\ge1}\{z:F^m(z)=z\}=\{0\}.
\]
It is \emph{$C^0$-rigid} if $F^{m_j}\to\id$ uniformly along a
sequence of positive integers $m_j\to\infty$.
For an orientation-preserving circle homeomorphism
$g$ with degree-one lift $\widehat g$, put
\[
 \rot(g)=\left[\lim_{m\to\infty}
 \frac{\widehat g^{\,m}(x)-x}{m}\right]\in\T,
 \qquad \rot(F)=\rot(F|_{\partial\D}),
 \qquad \normT{t}=\min_{k\in\Z}|t-k|.
\]
A \emph{Denjoy diffeomorphism} is an orientation-preserving circle diffeomorphism with
irrational rotation number and a nonempty open interval whose
distinct iterates are pairwise disjoint.

Differentiability is understood up to the boundary, and $C^{1,\beta}$
means that the map and its inverse have $\beta$-H\"older first
derivatives.  We do not identify $C^{1,1}$ with $C^2$.
Distances from rotations are measured in fixed Euclidean coordinates,
with the sum norm
\[
 \|u\|_{C^1}=\sup|u|+\sup\|Du\|.
\]
The full H\"older and coordinate conventions are in
\cref{H:sec:analytic-conventions}; the lifts and supports for the
construction are fixed in \cref{B:sec:circle-construction}.
All logarithms are natural.

If $\alpha$ is irrational with continued-fraction denominators
$q_n=q_n(\alpha)$, set
\[
 B(\alpha)=\sum_{n\ge0}\frac{\log q_{n+1}}{q_n},\qquad
 \mathcal E(\alpha)=\limsup_{n\to\infty}\frac{\log q_{n+1}}{q_n}.
\]
The Brjuno class is $\Bclass=\{\alpha:B(\alpha)<\infty\}$; we call
$\alpha$ \emph{non-Brjuno} when $B(\alpha)=\infty$, and we say that it
has \emph{zero exponential type} when $\mathcal E(\alpha)=0$.  As the
terms are nonnegative, the latter is the same as
$(\log q_{n+1})/q_n\to0$.
In particular, every Brjuno number has zero exponential type.
We distinguish finite positive type, $0<\mathcal E(\alpha)<\infty$,
from infinite type, $\mathcal E(\alpha)=\infty$.

\subsection{Main results}\label{sec:main-results}
\Needspace{8\baselineskip}
Our first result concerns the rigidity of \(C^{1,\beta}\) pseudo-rotations.

\begin{theorem}
\label{thm:rigidity}
Let $0<\beta\le1$.  Every area-preserving $C^{1,\beta}$ irrational
pseudo-rotation of $\D$ with non-Brjuno rotation number is $C^0$-rigid.
\end{theorem}

\Needspace{15\baselineskip}
The situation is different at the \(C^1\) level. In the following theorem, we construct non-rigid examples for every prescribed non-Brjuno rotation number of zero exponential type. This gives a negative answer to the question raised in \cite[\S1.3]{LWY}.  The remaining finite positive-type case is stated in \cref{sec:remaining-questions}.

\begin{theorem}\label{thm:nonrigidity}
Let $\alpha\in\T\setminus\Q$ be non-Brjuno with $\mathcal E(\alpha)=0$.
There exists an area-preserving $C^1$ irrational pseudo-rotation
$F:\D\to\D$ of rotation number $\alpha$ with the following properties:
\begin{enumerate}[label=\textnormal{(\roman*)}]
\item $F$ is not $C^0$-rigid;
\item $F|_{\partial\D}$ is a Denjoy diffeomorphism;
\item $F=\mathcal R_\alpha$ near $0$, and on
$\operatorname{int}\D\setminus\{0\}$ it is smoothly conjugate to
$\mathcal R_\alpha$;
\item neither $F|_{\partial\D}$ nor $F$ is $C^{1,\beta}$ for any $\beta>0$.
\end{enumerate}
\end{theorem}

\begin{remark}\label{rem:near-rotation}
For every $\varepsilon>0$, the map in \cref{thm:nonrigidity} may be
chosen so that
\[
 \|F-\mathcal R_\alpha\|_{C^1(\D)}<\varepsilon,
\]
with $\alpha$ fixed and all four properties retained.  As in
\cite[Theorem~2\textnormal{(i)}]{LWY}, this follows by starting the
same construction at a sufficiently late scale.  Here, however,
pairwise disjoint physical supports give a supremum tail estimate
in place of the summable Brjuno tail used there.  The quantitative
argument is given at the end of \S\ref{B:sec:disk-completion}.
\end{remark}

\Needspace{12\baselineskip}
Together with the Brjuno examples of \cite{LWY},
\cref{thm:nonrigidity} extends non-rigidity to the whole
zero-exponential-type class. We next consider the bi-Lipschitz case. We show that non-rigid examples also exist when the rotation number has finite positive exponential type.

\begin{theorem}
\label{thm:bilip-positive-type}\label{B:prop:bilip-positive-type}
For every irrational $\alpha$ with
$0<\mathcal E(\alpha)<\infty$, there is an orientation- and
standard-area-preserving bi-Lipschitz homeomorphism $F$ of $\D$ such that
\[
 \Per(F)=\{0\},\qquad \rot(F)=\alpha,
 \qquad \inf_{n\ge1}\|F^n-\id\|_{C^0}>0.
\]
It is the rigid rotation near $0$ and is smoothly conjugate by an area-preserving map to
$\mathcal R_\alpha$ on $\operatorname{int}\D\setminus\{0\}$.  Its boundary
restriction has a wandering interval.  The constructed map is smooth
in the open disk, but it is not $C^1$ on the closed disk.
\end{theorem}

\Needspace{5\baselineskip}
The examples constructed in \cref{thm:bilip-positive-type} have a Denjoy-type boundary behavior. However, they are not \(C^1\) on the closed disk. We next show that when the rotation number has positive exponential type, the \(C^1\) regularity of the map implies $C^0$-rigidity on the boundary.

\begin{theorem}
\label{thm:boundary-rigidity}\label{cor:boundary-rigidity}\label{B:prop:boundary-rigidity}
Let $F$ be an area-preserving $C^1$ irrational pseudo-rotation of
$\D$, with rotation number $\alpha$.  If $\mathcal E(\alpha)>0$,
then $F|_{\partial\D}$ is $C^0$-rigid, and hence is topologically
conjugate to $R_\alpha$ and has no wandering interval.
\end{theorem}

The proof is given in \cref{B:sec:scope}. Therefore, the construction using a Denjoy map on the boundary cannot produce a \(C^1\) non-rigid pseudo-rotation when the rotation number has positive exponential type. In the case of infinite exponential type, we can prove an even stronger result: every bi-Lipschitz pseudo-rotation is rigid;
see \cref{B:rem:bilip-scope}.

For comparison, let $\mathcal I_\alpha$ be the class of area-preserving
$C^1$ disk pseudo-rotations that equal $\mathcal R_\alpha$ near $0$
and are smoothly conjugate to it on the punctured interior.
This class contains a non-rigid map exactly when
$\mathcal E(\alpha)=0$ (\cref{cor:criterion}).  The criterion concerns
this specified interior structure; it does not classify all $C^1$
pseudo-rotations.

The picture is summarized in \cref{tab:scope}; the weighted Brjuno class
of \cite{LWY} that appears there is
\[
 \Bclassbeta{\beta}=
 \left\{\alpha:
 \sum_n\frac{q_{n+1}^{\beta}(1+\log q_{n+1})^{2\beta}}{q_n}
 <\infty\right\}\subset\Bclass,\qquad 0<\beta<1.
\]
For $0<\beta<1$, the non-$C^0$-rigid examples in this class show that
the non-Brjuno hypothesis in \cref{thm:rigidity} cannot be dropped.
They do not settle the $C^{1,1}$ or $C^2$ Brjuno case.

\begin{table}[htbp]
\caption{Rigidity and non-rigidity in the arithmetic ranges considered here.}
\label{tab:scope}
\centering\small
\renewcommand{\arraystretch}{1.13}
\begin{tabularx}{\textwidth}{@{}l >{\raggedright\arraybackslash}X >{\raggedright\arraybackslash}X@{}}
\toprule
Regularity & Rotation number & Conclusion\\
\midrule
$C^{1,\beta}$, $0<\beta\le1$ & Non-Brjuno
 & Every pseudo-rotation is rigid (\cref{thm:rigidity})\\
$C^1$ & Non-Brjuno, $\mathcal E=0$
 & Non-rigid examples (\cref{thm:nonrigidity})\\
$C^1$ & Brjuno
 & Non-rigid examples \cite{LWY}\\
$C^1$ & $0<\mathcal E<\infty$
 & Boundary rigid (\cref{thm:boundary-rigidity}); full-disk rigidity open\\
$C^k$, $k\ge2$ & Brjuno
 & Full-disk rigidity open\\
Bi-Lipschitz & $0<\mathcal E<\infty$
 & Non-rigid examples, not $C^1$ (\cref{thm:bilip-positive-type})\\
Bi-Lipschitz & $\mathcal E=\infty$
 & Every pseudo-rotation is rigid (\cref{B:rem:bilip-scope})\\
$C^{1,\beta}$, $0<\beta<1$ & Weighted Brjuno $\Bclassbeta{\beta}$
 & Non-rigid examples \cite{LWY}\\
\bottomrule
\end{tabularx}
\end{table}
\FloatBarrier

\Needspace{5\baselineskip}
\subsection{Two remaining rigidity questions}\label{sec:remaining-questions}
The first question concerns the remaining range of the problem in
\cite[\S1.3]{LWY}.
Positive exponential type implies non-Brjuno type, so
\cref{thm:rigidity} gives rigidity with any positive H\"older exponent.
Infinite exponential type is already rigid even in the bi-Lipschitz
class (\cref{B:rem:bilip-scope}).

\Needspace{8\baselineskip}
\begin{openquestion}
\leavevmode
\begin{enumerate}[label=\textnormal{(\roman*)}]
\item Let $\alpha$ be irrational with
$0<\mathcal E(\alpha)<\infty$.  Is every area-preserving $C^1$ disk
pseudo-rotation of rotation number $\alpha$ $C^0$-rigid?
\item Let $\alpha$ be Brjuno.  Is every area-preserving $C^k$ ($k\geq 2$) disk
pseudo-rotation of rotation number $\alpha$ $C^0$-rigid?
\end{enumerate}
\end{openquestion}

We plan to investigate these questions in subsequent work.

Neither question includes a derivative-growth hypothesis.
The second isolates the arithmetic range not covered by the
$C^k$ ($k\geq 2$) non-Brjuno rigidity theorem of \cite{AFLXZ}.

The counterexamples in \cite{LWY} and in \cref{thm:nonrigidity} obtain non-rigidity from a Denjoy map on the boundary of the disk. In \cref{B:prop:internal-denjoy-circle}, we show that the same mechanism can also be realized on a smooth invariant circle in the interior of the disk. The key point is to match the cotangent boundary jets in the construction. This allows us to preserve the standard area form and the existence of a unique periodic point.

The above construction requires the rotation number to have zero exponential type and the map is not $C^2$, therefore it does not answer either of the two questions above. Moreover, when the rotation number has positive exponential type, the boundary rigidity argument also implies that Denjoy dynamics cannot happen on any \(C^1\) embedded invariant circle. In the \(C^2\) case, Denjoy's theorem gives a similar obstruction: a \(C^2\) map cannot have a Denjoy restriction on a \(C^2\) embedded invariant circle (see \cref{B:rem:internal-circle-obstructions}).

These results suggest that we may need other methods to construct non-rigid pseudo-rotations in the remaining cases. These results do not imply that every non-rigid pseudo-rotation must contain a Denjoy invariant circle. The examples in \cref{thm:bilip-positive-type} do not answer either of the two questions above, since they are not \(C^1\).

\subsection{Proofs and organization}\label{sec:proof-strategy}\label{sec:organization}
The proof of the rigidity theorem is based on the finite-information closing argument developed in \cite[\S\S2--3]{AFLXZ}. The technical part is \cref{thm:finite-closing} which gives some conditions to show the existence of a hyperbolic periodic point. And the pseudo rotation condition contradicts the existence of a hyperbolic periodic point implies estimation of derivative-growth. Together with the arithmetic properties of the rotation number,we can prove \cref{thm:rigidity}. To prove \cref{thm:finite-closing}, we first select a suitable finite orbit segment with large derivative growth. And then we construct local coordinates near this orbit segment and study the images of certain cones and strips. The main difficulty is to carry out this argument with only Hölder continuity of the first derivative. To overcome this difficulty, we use a polynomial extension near the boundary and choose suitable scales for the local coordinate boxes.

The construction of non-rigid examples is based on the method developed in \cite[Lemmas~2.1--2.2 and~5.1]{LWY}. We follow their original local squeeze construction and its Hamiltonian realization. The main new ingredient is a different way of choosing the perturbation supports. And we use the two properties
\[
 \|\Delta_n\|_{C^1}\to0,\qquad
 \supp\Delta_n\cap\supp\Delta_m=\varnothing\quad(n\ne m),
\]
which give a $C^1$ Cauchy sequence by taking the largest remaining
increment (\cref{B:lem:localized-Cauchy}). Finally, the desired map is obtained by taking limit.

The main idea of the proof of \cref{thm:bilip-positive-type} is to group a fixed number of consecutive micro-squeezes so that the phase allocation remains possible. We then use disjoint supports to obtain a bi-Lipschitz limit. And the proof of \cref{thm:boundary-rigidity} relies on estimating the derivatives of iterates near the boundary.

\Cref{part:rigidity} proves the closing theorem and
\cref{thm:rigidity}.  \Cref{part:construction} constructs the
counterexamples and proves \cref{thm:nonrigidity}.
\Cref{B:sec:scope} proves \cref{thm:boundary-rigidity} and discusses the
integrable-interior and the internal circle criterion. \Cref{B:sec:bilip-positive-type} gives the bi-Lipschitz comparison and proves \cref{thm:bilip-positive-type}.

\Needspace{10\baselineskip}
\part{Rigidity at positive H\"older regularity}\label{part:rigidity}
\section{Finite-information closing and good orbit segments}\label{H:sec:finite-orbits}
Throughout this part, we assume that $(S,\omega)$ is a compact, connected, oriented smooth
surface, possibly with boundary.  All of the geometric data below are
fixed before any map is chosen.

\subsection{Analytic and directional conventions}\label{H:sec:analytic-conventions}
For a scalar, vector, or matrix field $u$ on a fixed coordinate domain
$E$, we write
\begin{gather*}
 \|u\|_{0;E}=\sup_{x\in E}|u(x)|,\qquad
 [u]_{\beta;E}=\sup_{\substack{x,y\in E\\x\ne y}}
           \frac{|u(x)-u(y)|}{|x-y|^\beta},\\
 \|u\|_{C^{k,\beta}(E)}=\sum_{j=0}^{k}\|D^ju\|_{0;E}+[D^ku]_{\beta;E}.
\end{gather*}
For a map $v$ defined on a coordinate domain $E$, we write
\[
\operatorname{Lip}(v;E)
:=
\sup_{\substack{x,y\in E\\x\ne y}}
\frac{|v(x)-v(y)|}{|x-y|}.
\]
We omit $E$ when the domain is clear.
For matrix fields and higher differentials, $|\cdot|$ is the induced operator norm.  We use $D^0u=u$.  When the domain is omitted, it is understood to be the coordinate set under discussion, and we also write $\|u\|_\infty$ for the same supremum norm.

Given nonnegative families $X=(X_i)$ and $Y=(Y_i)$ and a list
$\mathcal P$ of permitted fixed parameters, we write
\[
 X\lesssim_{\mathcal P}Y
 \iff \exists C=C(\mathcal P)>0\ \forall i,\ X_i\le CY_i,
 \qquad
 X\asymp_{\mathcal P}Y
 \iff X\lesssim_{\mathcal P}Y\text{ and }Y\lesssim_{\mathcal P}X.
\]
For sequences it is enough that these comparisons hold eventually.
We also write $X_n=O_{\mathcal P}(Y_n)$ when $|X_n|\lesssim_{\mathcal P}Y_n$,
and $X_n=o(Y_n)$ when $X_n/Y_n\to0$, with $Y_n>0$ eventually.
In this part the fixed parameters are $\beta$ and the surface geometry
chosen in \S\ref{H:sec:ambient}; unless another dependence is stated,
constants are uniform in the map, in the norm bounds and in the orbit
length.

We work with unoriented directions, that is, in the real projective
line
\[
 \mathbb{RP}^1=\{\R\xi:\xi\in\R^2\setminus\{0\}\},\qquad
 \mathbb P(\xi)=\R\xi,\qquad
 \mathbb P(A)=\{\R\xi:\xi\in A\setminus\{0\}\}.
\]
An invertible matrix $B$ acts on it by $(\mathbb P B)(\R\xi)=\R(B\xi)$.
If $u,u'$ are unit representatives of $\ell,\ell'\in\mathbb{RP}^1$, we put
\[
 d_{\mathbb{RP}^1}(\ell,\ell')=\arccos|\langle u,u'\rangle|,
 \qquad \cond B=\|B\|\|B^{-1}\|.
\]
For a cone $C$ (see Section \ref{H:sec:moving-boxes}), we write $BC=\{B\xi:\xi\in C\}$, and we adopt the
convention
\begin{equation}\label{eq:projective-cone-convention}
 BC\Subset C\quad\Longleftrightarrow\quad
 \operatorname{cl}_{\mathbb{RP}^1}\mathbb P(BC)
      \subset\operatorname{int}_{\mathbb{RP}^1}\mathbb P(C).
\end{equation}
Here $\operatorname{cl}$ denotes closure. For ordinary
coordinate domains, $E\Subset U$ keeps its usual meaning that
$\operatorname{cl}E$ is a compact subset of $U$.

\subsection{Ambient geometry and coordinate estimates}\label{H:sec:ambient}
\Needspace{5\baselineskip}
We first extend the surface $S$ by thickening the boundary and choose area coordinates on both sides of each boundary component.

\begin{lemma}
\label{H:lem:ambient-surface}
There are a smooth surface $\widehat S$ without boundary and a smooth
area form $\widehat\omega$ such that $S$ embeds as a compact domain in
$\widehat S$ and $\widehat\omega|_S=\omega$.
If $\partial S=\bigsqcup_{j=1}^{b}C_j$, there are disjoint collars
\[
 \chi_j:\T\times(-2\tau_*,2\tau_*)\longrightarrow\widehat S,
 \quad \chi_j^*\widehat\omega=ds\wedge dt,
 \quad \chi_j(\T\times\{0\})=C_j,
\]
with $t>0$ on the side of $S$.
\end{lemma}
\begin{proof}
If $\partial S=\varnothing$ we take $(\widehat S,\widehat\omega)=(S,\omega)$.
Otherwise every component $C_j$ is a circle, and we choose disjoint
smooth inward collars $\chi_j^+:\T\times[0,u_*)\to S$, with the
angular variable oriented so that
\[
 (\chi_j^+)^*\omega=a_j(s,u)\,ds\wedge du,\qquad a_j(s,u)>0.
\]
Put
\[
 t=t_j(s,u):=\int_0^u a_j(s,v)\,dv.
\]
Then $\partial_ut_j=a_j>0$ and
$ds\wedge dt_j=a_j\,ds\wedge du$.
Since there are only finitely many boundary circles and each is
compact, there is a common $\tau_*>0$ such that the reparameterized
inward collars are defined on $\T\times[0,3\tau_*)$.
To each of them, we attach the open exterior cylinder
$\T\times(-3\tau_*,0)$, keeping the same $(s,t)$ coordinates across
$t=0$.  The original interior charts and these two-sided collar
charts define a boundaryless smooth surface $\widehat S$ containing
$S$; no boundary appears at $t=-3\tau_*$, since that endpoint is not
included.  We let
$\widehat\omega=\omega$ on $S$ and
$\widehat\omega=ds\wedge dt$ on each added collar; the two definitions
agree on the overlaps, so $\widehat\omega$ is smooth and nonvanishing.
Restricting the collars to $|t|<2\tau_*$ gives the lemma.
\end{proof}

We fix one such ambient surface and a relatively compact neighborhood
$\mathcal N$ of $S$ in it, taking the uniform chart and exponential-map
bounds on a slightly larger compact neighborhood.  There is a
Riemannian metric compatible with the area form: if $g_0$ is any
smooth metric and $\widehat\omega=f\,\mathrm{vol}_{g_0}$ with $f>0$,
then in dimension two $g_1=f g_0$ has volume form $\widehat\omega$.
We fix such a metric $\widehat{\mathfrak g}$ and give $S$ the induced
metric $\mathfrak g=\iota^*\widehat{\mathfrak g}$, where
$\iota:S\hookrightarrow\widehat S$ is the inclusion, so that
$\mathrm{vol}_{\mathfrak g}=\omega$.
We write $d_{\widehat{\mathfrak g}}$ for the ambient
Riemannian distance and $d_{\mathfrak g}^S$ for the
intrinsic length distance on $S$.

Since $S$ is compact, connected, and smoothly embedded,
there exists a constant $C_{\mathrm{inc}}\ge 1$,
depending only on the fixed geometry, such that
for all $x,y\in S$,
\begin{equation}\label{H:eq:intrinsic-ambient-comparison}
	d_{\widehat{\mathfrak g}}(\iota x,\iota y)
	\le d_{\mathfrak g}^S(x,y)
	\le C_{\mathrm{inc}}\,
	d_{\widehat{\mathfrak g}}(\iota x,\iota y).
\end{equation}
The first inequality is immediate, while the second
follows from uniform local distance comparisons in
smooth interior and boundary charts, together with
compactness.

From now on, all tangent norms, operator norms, and angles are computed
with respect to $\widehat{\mathfrak g}$ on the ambient surface
and its restriction $\mathfrak g$ on $S$.

For $0<a\le\tau_*$ we write
\[
 \mathcal C_j^-(a):=\chi_j(\T\times(-a,0)),\qquad
 \mathcal N_a:=S\cup\bigcup_{j=1}^b\mathcal C_j^-(a).
\]
The width $a$ of a collar is thus always measured in the fixed normal
coordinate $t$ and never in a map-dependent reparameterization.
A diffeomorphism $g:S\to S$ induces a permutation $\pi$ of the boundary
circles, with $g(C_j)=C_{\pi(j)}$, and the corresponding calculations
are carried out in the source--target pair $(\chi_j,\chi_{\pi(j)})$.
The constants attached to these finitely many pairs are absorbed into
the fixed geometric constants.

\subsubsection*{Coordinate H\"older seminorms}\label{sec:coordinate-calculus}
We fix a finite atlas on the chosen neighborhood, including local
lifts of the collars and smaller chart domains compactly contained in
larger ones.  For a compatible local representative
$f_{ji}=\chi_j^{-1}\circ f\circ\chi_i$ the derivative seminorm is
\[
 [Df_{ji}]_{\beta;E_{ji}}
 =\sup_{\substack{x,y\in E_{ji}\\x\ne y}}
       \frac{\|Df_{ji}(x)-Df_{ji}(y)\|}{|x-y|^\beta},
\]
where
\[
E_{ji}
=
\left\{
x\in U_i^{\mathrm{sm}}:
\chi_i(x)\in S,\;
f(\chi_i(x))\in
\chi_j(U_j^{\mathrm{sm}})
\right\},
\]
and $U_i^{\mathrm{sm}}$ and $U_j^{\mathrm{sm}}$
denote the fixed smaller source and target coordinate
domains, respectively.
We write $[Df]_\beta$ for the maximum over these finitely many
compatible representatives, while $\|Df\|_0$, with no chart subscript,
stands for the Riemannian operator norm.

\Needspace{5\baselineskip}
The following calculus estimates keep track of the losses under the fixed coordinate changes and under inversion.

\begin{lemma}\label{H:lem:chart-change}
Let $E,F$ be coordinate domains. The product estimate
for bounded $\beta$-H\"older maps $u,v$ on $E$,
and the composition estimate for Lipschitz
$v:E\to F$ and $\beta$-H\"older $u$ on $F$,
are respectively
\begin{align}
 [uv]_\beta&\le[u]_\beta\|v\|_0+\|u\|_0[v]_\beta,
 & [u\circ v]_\beta&\le[u]_\beta\Lip(v)^\beta.
 \label{eq:holder-calculus}
\end{align}
For fixed smooth coordinate changes $a,b$, put
$\widehat f=b\circ f\circ a$ and $M_f=\|Df\|_0$.
On common smaller domains where $\Lip(f)\le C M_f$,
\begin{equation}\label{eq:atlas-change}
 \|D\widehat f\|_0\le C M_f,\qquad
 [D\widehat f]_\beta\le
 C\bigl([Df]_\beta+M_f+M_f^{1+\beta}\bigr).
\end{equation}
If a $C^{1,\beta}$ diffeomorphism $\Psi$ has inverse derivative norm
at most $M_-$ and inverse Lipschitz constant at most $L_-$, then
\begin{equation}\label{eq:inverse-holder-calculus}
 [D\Psi^{-1}]_\beta\le M_-^2[D\Psi]_\beta L_-^\beta.
\end{equation}
\end{lemma}
\begin{proof}
The inequalities \eqref{eq:holder-calculus} follow from the difference
identities for products and compositions.  By definition of $\widehat f$,
$D\widehat f=(Db\circ f\circ a)(Df\circ a)Da$.
The three factors are bounded in $C^0$ by $C,CM_f,C$, respectively,
and their H\"older seminorms satisfy
\[
 [Db\circ f\circ a]_\beta\le C M_f^\beta,\quad
 [Df\circ a]_\beta\le C[Df]_\beta,\quad [Da]_\beta\le C.
\]
For the first of these, we use that $Db$ is bounded and Lipschitz, and
hence uniformly $\beta$-H\"older on the fixed domain.  The product
inequality now gives \eqref{eq:atlas-change}.
Finally, for $x=\Psi^{-1}(u)$ and $y=\Psi^{-1}(v)$ the inverse-matrix
identity gives
\[
 \begin{split}
 \|D\Psi^{-1}(u)-D\Psi^{-1}(v)\|
 &\le M_-^2\|D\Psi(x)-D\Psi(y)\|\\
 &\le M_-^2[D\Psi]_\beta L_-^\beta|u-v|^\beta,
 \end{split}
\]
which is \eqref{eq:inverse-holder-calculus}.
\end{proof}
\noindent\emph{Uniformity on moving coordinate domains.}
The sets $E_{ji}$ above may vary with the map, but the larger and
smaller source--target charts, the smooth coordinate changes between
them and their positive buffers are all fixed in advance.  We take the
coordinate cores to be convex disks, half-disks or lifted collar
rectangles, as appropriate.
Let $f:S\to S$ be a map with
$M=\max\{1,\|Df\|_0\}$.
For $x,y\in E_{ji}$, put $h=|x-y|$.
If $h\le c_0/M$, where $c_0>0$ depends only on
the fixed atlas and its geometric buffers, the
coordinate segment joining $x$ and $y$ lies
in a larger source chart and represents a path
in $S$ of length at most $Ch$.
Its image under $f$ has length at most $CMh$.
Since $f(\chi_i(x))$ belongs to the smaller target
chart, choosing $c_0$ sufficiently small ensures
that the entire image path stays in the larger
target chart. Writing $a=\chi_i$ and $b=\chi_j^{-1}$ for the
fixed source and target coordinate maps, respectively,
and computing derivatives in fixed local coordinates,
the uniform smooth bounds for $b$ and the image-path
length estimate give
\[
 \|Db(f(a(x)))-Db(f(a(y)))\|\le CMh
       \qquad (h\le c_0/M).
\]
For $h>c_0/M$ we use instead the boundedness of $Db$, which gives
\[
 \frac{\|Db(f(a(x)))-Db(f(a(y)))\|}{h^\beta}
 \le C M^\beta.
\]
Near and far pairs together give
$[Db\circ f\circ a]_\beta\le CM^\beta$.  Using the product rule for
near pairs and the bound
$\|D\widehat f\|_0\le CM$ for far pairs, we obtain, with no convexity
assumption on $E_{ji}$,
\[
 [D\widehat f]_{\beta;E_{ji}}
 \le C\bigl([Df]_\beta+M+M^{1+\beta}\bigr).
\]
This covers the original surface map on its compatible moving
chart domains.  For maps that are only defined on nonconvex collar
images, such as the inverse coordinate maps of the extension
construction, we check the required Lipschitz controls directly in
\cref{H:lem:symplectic-extension} instead of deducing them from
derivative bounds.  The full $M^{1+\beta}$ loss in
\eqref{eq:atlas-change} is carried along in the polynomial bounds
below.

We write $\Diff^{1,\beta}_{\mathrm{vol}}(S)=\Diff^{1,\beta}_{\omega}(S)$
for the $C^{1,\beta}$ diffeomorphisms preserving the fixed area form,
and we use the same convention for class $C^1$. The constants
$C_S,C_{S,\beta}$ may increase as the argument proceeds, whereas named thresholds are fixed after all requirements on them have been collected. 

\Needspace{6\baselineskip}
For area-preserving $C^{1,\beta}$ maps, the theorem below
uses first-derivative and H\"older bounds to show that we can obtain a hyperbolic periodic point through a sufficiently long finite orbit segment with a sufficiently rapid derivative growth.

\begin{sectiontheorem}\label{thm:finite-closing}
Fix $0<\beta\le1$.  There are constants
\[
 A_*=A_*(S,\beta)>1,\qquad
 \vartheta_* =\vartheta_*(\beta)\in(0,1),\qquad
 H_*=H_*(S,\beta)>1
\]
with the following property.  Suppose that
$g\in\Diff^{1,\beta}_{\mathrm{vol}}(S)$ and that
\begin{equation}\label{H:eq:finite-data}
 \|Dg\|_0\le A,\qquad [Dg]_{\beta}\le D,\qquad
 D\ge A\ge A_*,\qquad q\in\N,\qquad q\ge D^{H_*}.
\end{equation}
If
\begin{equation}\label{H:eq:large-growth}
                       \|Dg^q\|_0\ge A^{\vartheta_*q},
\end{equation}
then $g$ has a hyperbolic periodic point.  If $\partial S\ne\varnothing$
and $g$ has no periodic point on $\partial S$, this periodic point lies in
$\operatorname{int}S$.
\end{sectiontheorem}

This is a H\"older version of the moving-box closing argument of
\cite[\S3]{AFLXZ}; compare also the frameworks of \cite{Katok,CP}.
The scales are changed according to
\[
 p=1+\frac2\beta,\qquad
 \psi_p(t)=\begin{cases}t,&t\le0,\\t/p,&t>0,\end{cases}
 \qquad
 c_{n+1}=\min\{e^{\psi_p(\lambda_n^e-\delta)}c_n,100\},
 \qquad r_n=\bar r c_n^p,
\]
where $\lambda_n^e$ is the effective expansion defined below.
Positive growth is slowed down so that the horizontal box crosses,
while $p\beta-1=\beta+1$ lets the nonlinear error absorb the
inverse angle.  The cone estimates then produce a returned
hyperbolic branch, and the boundary is handled by a two-sided
symplectic extension.

\FloatBarrier

\subsection{Good finite orbit segments}
We now recall the construction of good finite orbit segments from \cite{AFLXZ}. And we will use a slightly different choice of constants to fit the scales introduced above. Write
$x_i=g^i(x_0)$, and transport unit vectors $v_0^u,v_0^s\in T_{x_0}S$ along
the orbit by
\[
 v_i^{u,s}=\frac{Dg^i(v_0^{u,s})}{\|Dg^i(v_0^{u,s})\|}.
\]
Set
\[
 \lambda_i^u=\log\|Dg(v_i^u)\|,
 \qquad
 \lambda_i^s=\log\|Dg(v_i^s)\|,
 \qquad
 \bar\lambda_i^e=\min\{\lambda_i^u,-\lambda_i^s\},
\]
and
\begin{equation}\label{H:eq:refined-rate}
 \lambda_i^e=\min\{\lambda_i^u,
              \lambda_i^u-\lambda_i^s,-2\lambda_i^s\}.
\end{equation}

For the rest of the proof we fix
\begin{equation}\label{H:eq:small-constants}
 p=1+\frac2\beta,\qquad
 \eta=\frac1{100p^2},\qquad
 \rho=\frac1{100p},\qquad
 \delta=\rho a,\qquad a=\log A.
\end{equation}

\begin{definition}\label{H:def:good-triple}
A triple $(x_0,v_0^s,v_0^u)$ is forward $(L,a)$-good if
$|\lambda_i^{u,s}|\le a$ for $0\le i<L$ and
\[
 \frac1k\sum_{i=0}^{k-1}\bar\lambda_i^e>(1-\eta)a
       \qquad(1\le k\le L).
\]
It is backward $(L,a)$-good if $|\lambda_i^{u,s}|\le a$ for $0\le i<L$
and
\[
 \frac1k\sum_{i=L-k}^{L-1}\bar\lambda_i^e>(1-\eta)a
       \qquad(1\le k\le L).
\]
Given an ambient integer $q\ge L$, a \emph{returned good segment relative
to $q$} is both forward and backward $(L,a)$-good and satisfies the following two conditions, with $\sigma>0$ fixed:
\begin{equation}\label{H:eq:endpoint-return}
 d(x_0,x_L)+d_{T^1S}(v_0^u,v_L^u)
       +d_{T^1S}(v_0^s,v_L^s)<q^{-\sigma}
\end{equation}
\begin{equation}\label{H:eq:endpoint-angle}
 \cot\angle(v_0^u,v_0^s)+\cot\angle(v_L^u,v_L^s)\le2A^3.
\end{equation}
Here and below, angles are taken in $(0,\pi/2]$
after forgetting orientations.
We write
$d(x,y):=d_{\widehat{\mathfrak g}}(\iota x,\iota y)$
for $x,y\in S$, and let $d_{T^1S}$ denote the
distance of a fixed Sasaki-type metric induced by
$\mathfrak g$ on the oriented unit tangent bundle.  Any other fixed smooth
bundle metric would give equivalent return estimates, after
adjusting fixed geometric constants.  Thus
\eqref{H:eq:endpoint-return} also aligns the orientations of the endpoint
vectors.
\end{definition}

We take $\sigma=1/100$ and $0<4\sigma<1$.

\Needspace{5\baselineskip}

\begin{lemma}\label{H:lem:good-segment}
Let $g\in\Diff^1_{\mathrm{vol}}(S)$.  For every $\eta>0$ there are
$0<\vartheta_-(\eta)<\vartheta_+(\eta)<1$, $A_0(\eta,S)>1$,
and $H_0(\eta,S)>1$ such that the following holds.  If
\[
 \|Dg\|_0\le A,\qquad
 \|Dg^q\|_0\ge A^{\vartheta q},\qquad
 A\ge A_0,\quad \vartheta\ge\vartheta_+,\quad q\ge A^{H_0},
\]
then there are $x_0,v_0^s,v_0^u$ and $1\le L\le q$ forming a returned
good segment relative to $q$ in the sense of \cref{H:def:good-triple}.
If $S$ has boundary, the segment can be chosen in $\operatorname{int}S$.
\end{lemma}

\begin{proof}
It is enough to treat $0<\eta\le1/100$, since for larger $\eta$ one may
impose the stronger goodness condition with $\eta'=\min\{\eta,1/100\}$.
We choose $\vartheta_-$ close enough to one for the counting argument
below to go through, and then fix any $\vartheta_+\in(\vartheta_-,1)$.
Since $y\mapsto\|Dg^q(y)\|$ is continuous and $\operatorname{int}S$ is
dense in $S$, the strict exponent reserve lets us pick
$x\in\operatorname{int}S$ with
\[
                    \|Dg^q(x)\|>A^{\vartheta_-q}.
\]
(When $S$ is closed this just says $x\in S$.)  Let $v_0^s$ be the most
contracted direction of $Dg^q(x)$ and $v_0^u$ the direction orthogonal
to it.  The preserved area form is the Riemannian volume, so every
$Dg(y)$ has determinant one; the two singular values of $Dg^q(x)$ are
therefore reciprocal, and area preservation gives
\[
 \frac1q\sum_{i=0}^{q-1}\lambda_i^u>\vartheta_- a,
 \qquad
 \frac1q\sum_{i=0}^{q-1}\lambda_i^s<-\vartheta_- a,
 \qquad |\lambda_i^{u,s}|\le a.
\]
We include the short counting argument because it fixes the dependence
on $\eta$.  Put $\gamma_{\rm sel}=\eta/32$.  The
total deficiency of $\lambda_i^u$ from $a$, and likewise that of
$-\lambda_i^s$, is at most $(1-\vartheta_-)qa$.  The set of indices at
which either $\lambda_i^u<(1-\gamma_{\rm sel})a$ or
$-\lambda_i^s<(1-\gamma_{\rm sel})a$ therefore has at most
$2(1-\vartheta_-)q/\gamma_{\rm sel}$ elements.  As
$\bar\lambda_i^e\ge-a$ everywhere, this gives
\[
 \frac1q\sum_{i=0}^{q-1}\bar\lambda_i^e
 \ge\left(1-\gamma_{\rm sel}-\frac{4(1-\vartheta_-)}{\gamma_{\rm sel}}\right)a.
\]
We chose $\vartheta_-$ so close to one that the right-hand side exceeds
$(1-\eta/8)a$.

The form of Pliss' lemma \cite{Pliss} we need is elementary: if
$z_i\le a$ and $q^{-1}\sum z_i>\ell'$, then for every $\ell''<\ell'$
there are at least $(\ell'-\ell'')q/(a-\ell'')$ indices $n$ with
\[
 \frac1k\sum_{i=n}^{n+k-1}z_i>\ell''
 \quad\text{for every }1\le k\le q-n.
\]
We apply this to $z_i=\bar\lambda_i^e$, and again to the reversed
sequence, both times with $\ell'=(1-\eta/8)a$ and $\ell''=(1-\eta)a$.
The forward times so obtained lie in $\{0,\ldots,q-1\}$ and the backward
endpoint times in $\{1,\ldots,q\}$; we restrict both sets to
$\{1,\ldots,q-1\}$, which discards at most one endpoint from each.  Since
Pliss' lemma supplies at least $7q/8$ indices before the restriction,
the restricted sets satisfy
\[
 \min\{|\mathcal F|,|\mathcal B|\}
 \ge\frac{7q}{8}-1>\frac{3q}{4}\qquad(q>8).
\]
By construction,
\begin{align*}
 n\in\mathcal F&\Longrightarrow
 \frac1k\sum_{i=n}^{n+k-1}\bar\lambda_i^e>(1-\eta)a
       \quad(1\le k\le q-n),\\
 n\in\mathcal B&\Longrightarrow
 \frac1k\sum_{i=n-k}^{n-1}\bar\lambda_i^e>(1-\eta)a
       \quad(1\le k\le n).
\end{align*}
Hence $\mathcal F\cap\mathcal B$ contains more than $q/2$ indices.

For $n\in\mathcal B$ we also control the angle: the area-preserving
angle recurrence
\begin{equation}\label{H:eq:cot-recurrence}
 \cot\angle(v_{i+1}^u,v_{i+1}^s)
 \le e^{2\lambda_i^s}\cot\angle(v_i^u,v_i^s)+A^2
\end{equation}
combined with the backward stable averages yields
\[
 \cot\angle(v_n^u,v_n^s)
 \le A^2\sum_{j=0}^{n-1}
 \exp\left(2\sum_{k=n-j}^{n-1}\lambda_k^s\right)
 \le \frac{A^2}{1-A^{-2(1-\eta)}}\le A^3
\]
once $A_0$ is large; the initial cotangent vanishes because
$v_0^u\perp v_0^s$.  It is here that area preservation, in
dimension two, is essential.

We now produce the return.  Consider the compact fiber product
\[
 \mathcal E=\{(x,v^s,v^u):x\in S,\ v^s,v^u\in T_x^1S\}.
\]

Equip $\mathcal E$ with the sum of the ambient
base-point distance and the two unit-tangent-bundle
distances. By subdividing a fixed finite atlas,
we can cover $\mathcal E$ by at most
$C_S q^{4\sigma}$ sets of diameter smaller
than $q^{-\sigma}/3$.
  As $4\sigma<1$, once
$q\ge A^{H_0}$ with $H_0$ large, two of the more than $q/2$ indices in
$\mathcal F\cap\mathcal B$, say $i<i+L$, fall into the same covering
set.  We rebase the segment at $g^i(x)$.  Because $i$ lies in
$\mathcal F$, every forward prefix inequality holds up to time $L$, and
because $i+L$ lies in $\mathcal B$, every backward suffix inequality
holds as well.  The common covering set yields
\eqref{H:eq:endpoint-return}, and the two angle estimates yield
\eqref{H:eq:endpoint-angle}.  This is the argument of
\cite[Proposition~3.1]{AFLXZ}, with its numerical threshold replaced by
$\eta$.  Finally, a diffeomorphism of a manifold with boundary maps the
interior to itself, so every point of the rebased segment stays in
$\operatorname{int}S$.
\end{proof}

The estimates below will be applied to every prefix of a forward-good
segment and to every suffix of a backward-good segment.

\begin{lemma}\label{H:lem:effective-rates}
Let $(x_0,v_0^s,v_0^u)$ be a returned good segment relative to some $q$.
For every $1\le k\le L$, take either
\[
 J=J_k^+:=\{0,\ldots,k-1\}\quad\text{or}\quad
 J=J_k^-:=\{L-k,\ldots,L-1\}.
\]
Then
\begin{align}
 &\sum_{i\in J}\lambda_i^s
     \le-(1-\eta)|J|a,\label{H:eq:stable-average}\\
 &\sum_{i\in J}\lambda_i^e
     \ge(1-4\eta)|J|a,\label{H:eq:effective-average}\\
 &\sum_{i\in J}\min\{\lambda_i^e,0\}
     \ge-8\eta|J|a.\label{H:eq:negative-part}
\end{align}
\end{lemma}

\begin{proof}
Fix one of the admissible sets $J$ and write $m=|J|$.  Goodness, together
with the pointwise bounds, gives
\[
 \sum_{i\in J}\bar\lambda_i^e>(1-\eta)ma,
 \qquad -a\le\bar\lambda_i^e\le a,
 \qquad -2a\le\lambda_i^e\le a.
\]
As $-\lambda_i^s\ge\bar\lambda_i^e$, the first of these already gives
\eqref{H:eq:stable-average}.  Next, whenever $\lambda_i^s\le0$, every
entry of the minimum defining $\lambda_i^e$ is at least
$\bar\lambda_i^e$.  Let $E=\{i\in J:\lambda_i^s>0\}$.  On $E$ we have
$\bar\lambda_i^e\le0$, whence
\[
 (1-\eta)ma<\sum_{i\in J}\bar\lambda_i^e
       \le(m-|E|)a,\qquad |E|<\eta m.
\]
On $E$ we use instead $\lambda_i^e-\bar\lambda_i^e\ge-3a$, and altogether
\[
 \sum_{i\in J}\lambda_i^e
 \ge\sum_{i\in J}\bar\lambda_i^e-3a|E|
 >(1-4\eta)ma,
\]
which proves \eqref{H:eq:effective-average}.
Finally, $\lambda_i^e\le0$ forces $\bar\lambda_i^e\le0$, since if
$\lambda_i^u>0$ and $\lambda_i^s<0$ then all three entries defining
$\lambda_i^e$ are positive.  The same count therefore gives
\[
 \#\{i\in J:\lambda_i^e\le0\}<\eta m,
 \qquad
 \sum_{i\in J}\min\{\lambda_i^e,0\}\ge-2\eta ma.
\]
This is stronger than \eqref{H:eq:negative-part}, whose constant is
the one that enters the later scale estimates.
\end{proof}

\section{A polynomial boundary extension}\label{H:sec:extension}
In this section, we construct a local extension of the map to a neighborhood of \(S\) and make some estimates in the following lemma. 

We call each boundary circle
$C_j=\chi_j(\mathbb T\times\{0\})$ a seam,
separating $S$ from its attached exterior collar
in $\widehat S$. For a source--target collar pair
$(\chi_j,\chi_{\pi(j)})$, the corresponding seams
are represented by $t=0$ and $T=0$, respectively.

\begin{lemma}
\label{H:lem:symplectic-extension}
Let $(S,\omega)\subset(\widehat S,\widehat\omega)$ be as above and fix
$0<\beta\le1$.  There are constants
$C_{\rm ext}=C_{\rm ext}(S,\omega,\beta)\ge1$,
$D_{\rm ext}=D_{\rm ext}(S,\omega,\beta)>1$, and an exponent
$m_{\rm ext}=m_{\rm ext}(\beta)\ge1$ with the following property.
Suppose
\[
 g\in\Diff^{1,\beta}_\omega(S),\qquad
 \|Dg\|_0\le A,\qquad [Dg]_\beta\le D,\qquad D\ge A\ge D_{\rm ext}.
\]
Set $\Delta=C_{\rm ext}D^{m_{\rm ext}}$.  There are neighborhoods
$\mathcal N_{\rm src},\mathcal N_{\rm tar}$ of $S$ in $\widehat S$
and a $C^{1,\beta}$ symplectic diffeomorphism
$\widetilde g:\mathcal N_{\rm src}\to\mathcal N_{\rm tar}$ such that:
\begin{enumerate}[label=\textnormal{(\roman*)}]
\item $\widetilde g|_S=g$;
\item both neighborhoods contain the respective exterior collars of
      width $\Delta^{-1}$;
\item $\widetilde g$ sends the interior side to the interior side and the
      exterior side to the exterior side of each source--target seam;
\item in the fixed atlas, on the corresponding source and target domains,
\[
 \|D\widetilde g^{\pm1}\|_0+[D\widetilde g^{\pm1}]_\beta\le\Delta.
\]
\end{enumerate}
One admissible, non-optimal choice is $m_{\rm ext}=100+20/\beta$.
\end{lemma}

\begin{proof}
We divide the proof into several steps. We first use the symplectic identity to construct a generating function in mixed coordinates. We then extend this function across the boundary and use it to define the desired map.
We extend its scalar correction across the seam and use positivity
of the mixed derivative to recover an injective symplectic collar
map with the same boundary jet.
The final estimates control the map, its inverse, and the
reciprocal collar widths by fixed powers of the original bound $D$.

If $\partial S=\varnothing$, we take
$\widehat S=S$ and $\widetilde g=g$.
The required inverse bounds follow from
\cref{H:lem:chart-change} and area preservation, after
enlarging the fixed constants.
Henceforth assume $\partial S\ne\varnothing$.

\textit{1. Boundary coordinates and a common mixed-coordinate strip.}
Fix a source--target collar pair
$(\chi_j,\chi_{\pi(j)})$, where
$g(C_j)=C_{\pi(j)}$.
In these coordinates, write
\[
(\chi_{\pi(j)}^{-1}\circ g\circ\chi_j)(s,t)
=(\Sigma(s,t),T(s,t)),
\qquad t\ge 0.
\]
Here $\Sigma$ denotes the angular coordinate function; it should not be
confused with the surface $S$.  By \eqref{eq:atlas-change} we may fix
$C_{\rm at}\ge1$ and work with
\begin{equation}\label{eq:extension-coordinate-budget}
 D_{\rm c}:=C_{\rm at}(1+D+A^{1+\beta}),\qquad
 D\le D_{\rm c}\le3C_{\rm at}D^2,
\end{equation}
which gives, in these collar coordinates,
\[
 \|D(\Sigma,T)\|_0\le C_bA,\qquad
 [D(\Sigma,T)]_\beta\le D_{\rm c}.
\]
The constants $C_b\ge1$ and $c_b>0$ of this step are fixed.  In
area-compatible orthonormal frames a
two-dimensional matrix of determinant one has singular values
$\sigma,\sigma^{-1}$, so that $\|Dg^{-1}\|_0=\|Dg\|_0\le A$ in the
Riemannian norm.  On the seam,
\begin{gather*}
 T(s,0)=T_s(s,0)=0,\qquad
 Dg(s,0)=\begin{pmatrix}\Sigma_s&\Sigma_t\\0&T_t\end{pmatrix},\\
 \Sigma_s(s,0)T_t(s,0)=1,\qquad T_t(s,0)>0,\qquad
 c_bA^{-1}\le T_t(s,0)\le C_bA.
\end{gather*}
The strict sign comes from invertibility together with preservation of
the inner side.  Choose $c_h>0$ with $c_h^\beta\le c_b/2$ and put
\begin{equation}\label{H:eq:first-extension-width}
 h_0=c_h(AD_{\rm c})^{-1/\beta}
       \ge c_hD_{\rm c}^{-2/\beta}.
\end{equation}
The source collar is chosen so that its image under
$g$ remains in the fixed target collar. Indeed,
each radial source segment of coordinate length
$t$ has Riemannian length at most $Ct$, and its
image under $g$ has length at most $CAt$.
The fixed target collar has a uniform tubular
margin, so this containment follows by choosing
$c_h$ sufficiently small. For $0\le t\le h_0$,
\[
 |T_t(s,t)-T_t(s,0)|\le D_{\rm c}h_0^\beta
       =c_h^\beta A^{-1}\le\tfrac12c_bA^{-1}.
\]
This yields
\begin{equation}\label{H:eq:normal-monotonicity}
 \tfrac12c_bA^{-1}\le T_t(s,t)\le2C_bA.
\end{equation}
Since $T(s,0)=0$, integration gives
$T(s,h_0)\ge(c_b/2)A^{-1}h_0$, so the image of
$\Psi(s,t)=(s,T(s,t))$ contains the product strip
\[
 \mathcal C_+=\T\times[0,T_*],\qquad
 T_*=(c_b/4)A^{-1}h_0.
\]
Strict monotonicity in $t$ makes $\Psi$ one-to-one onto its image, and
on $\mathcal C_+$ we may write
\[
 \Psi^{-1}(s,T)=(s,P(s,T)),\qquad
 Q(s,T)=\Sigma(s,P(s,T)).
\]

\textit{2. A normalized periodic generating function.}
We normalize the real lift by requiring $Q(0,0)\in[0,1)$.  The boundary
circle map is a circle homeomorphism with $\Sigma_s(s,0)=1/T_t(s,0)>0$,
so it has degree one.  The degree cannot change through the collar,
hence
\begin{equation}\label{H:eq:Q-monodromy}
 P(s+1,T)=P(s,T),\qquad Q(s+1,T)=Q(s,T)+1.
\end{equation}
In the mixed variables the symplectic identity becomes
\[
 dQ\wedge dT=ds\wedge dP
 \quad\Longleftrightarrow\quad Q_s=P_T.
\]
Because $P(s,0)=0$, the explicit primitive
\begin{equation}\label{eq:explicit-generating-function}
 W(s,T)=\int_0^TQ(s,u)\,du,\qquad U(s,T)=W(s,T)-sT
\end{equation}
satisfies
\begin{align*}
 W_T&=Q,&W_s&=\int_0^TP_T(s,u)\,du=P,\\
 W(s+1,T)-W(s,T)&=T,&U(s+1,T)&=U(s,T),\qquad U(s,0)=0.
\end{align*}
In particular $W,U\in C^{2,\beta}$: the extra derivative sits in the
\emph{generating function}, while the map it generates remains
$C^{1,\beta}$.

\textit{3. Polynomial bounds before extension.}
The inverse differential is
\[
 (D\Psi)^{-1}=\begin{pmatrix}1&0\\-T_s/T_t&1/T_t\end{pmatrix},
 \qquad \|D\Psi^{-1}\|_0\le CD_{\rm c}^{2}.
\]
Since the inverse lift is defined on the rectangle
$[0,1]\times[0,T_*]$, its Lipschitz constant obeys the same bound.  By
\eqref{eq:inverse-holder-calculus},
\begin{align*}
 [DP]_\beta&\le C D_{\rm c}^{4}
          D_{\rm c}(D_{\rm c}^{2})^\beta
       =CD_{\rm c}^{5+2\beta}\le CD_{\rm c}^{7},\\
 \|DP\|_0&\le CD_{\rm c}^{2}.
\end{align*}
In addition,
\[
 P_T=(T_t^{-1})\circ\Psi^{-1},\qquad
 Q_s=P_T,\qquad Q_T=(\Sigma_t\circ\Psi^{-1})P_T.
\]
These identities give
\begin{align*}
 \|D(Q-s)\|_0&\le CD_{\rm c}^{2},\\
 [D(Q-s)]_\beta
 &\le C\bigl(D_{\rm c}^{3+2\beta}+D_{\rm c}^{8}\bigr)
 \le CD_{\rm c}^{8}.
\end{align*}
The additive constant in the lift is controlled by the normalization
$Q(0,0)\in[0,1)$, and $U_s=P-T$, $U_T=Q-s$, $U(0,0)=0$.  Integrating
over a fundamental strip, whose diameter is bounded independently of
its height, we obtain
\begin{equation}\label{H:eq:generating-polynomial-bound}
 \|P\|_{C^{1,\beta}}+\|Q-s\|_{C^{1,\beta}}
             +\|U\|_{C^{2,\beta}}\le C_{\rm gen}D_{\rm c}^{8}.
\end{equation}
These bounds are uniform across the angular
identification $s=0\sim s=1$, using compatible
real lifts. The normalization $Q(0,0)\in[0,1)$
removes the additive integer ambiguity in $Q$.

\textit{4. Extension of the scalar function and positivity.}
For $-T_*/3\le T<0$ define
\begin{equation}\label{eq:seeley-explicit}
 \widetilde U(s,T)=6U(s,-T)-8U(s,-2T)+3U(s,-3T),
 \qquad \widetilde U=U\quad(T\ge0).
\end{equation}
The coefficients $(a_1,a_2,a_3)=(6,-8,3)$ satisfy
\[
 \sum_{j=1}^3a_j=1,\qquad
 \sum_{j=1}^3(-j)a_j=1,\qquad
 \sum_{j=1}^3j^2a_j=1.
\]
Every trace $\partial_s^u\partial_T^vU$ with $u+v\le2$ therefore
matches across $T=0$.  On each side separately, composing with
$(s,T)\mapsto(s,-jT)$ costs at most a fixed factor $j^{2+\beta}$ in a
second-derivative H\"older bound.  For two points on opposite sides we
split the segment joining them at the seam; continuity of the traces
then gives
\[
 |D^2\widetilde U(z)-D^2\widetilde U(z')|
 \le C\|U\|_{C^{2,\beta}}
       (|z-z_0|^\beta+|z'-z_0|^\beta)
 \le C_\beta\|U\|_{C^{2,\beta}}|z-z'|^\beta.
\]
Hence, with $\widetilde W=sT+\widetilde U$,
\[
 \|\widetilde U\|_{C^{2,\beta}}
       +\|\widetilde W\|_{C^{2,\beta}([0,1]\times[-T_*/3,T_*])}
       \le C_eD_{\rm c}^{8}.
\]
At $T=0$,
\[
 \widetilde W_{sT}(s,0)=P_T(s,0)=T_t(s,0)^{-1}
       \ge(2C_bA)^{-1}\ge(2C_bD_{\rm c})^{-1}.
\]
Now choose $c_1>0$, independent of $A,D$, such that
\begin{equation}\label{H:eq:second-extension-width}
 h_1=c_1D_{\rm c}^{-9/\beta},\qquad
 3h_1\le T_*,\qquad C_ec_1^\beta\le(4C_b)^{-1}.
\end{equation}
The middle condition can be met because
$T_*\ge cD_{\rm c}^{-(1+2/\beta)}$ and
$9/\beta>1+2/\beta$.  For $|T|\le h_1$,
\[
 |\widetilde W_{sT}(s,T)-\widetilde W_{sT}(s,0)|
 \le C_eD_{\rm c}^{8}h_1^\beta
 =C_ec_1^\beta D_{\rm c}^{-1}.
\]
Therefore, for $|T|\le h_1$,
\begin{equation}\label{H:eq:mixed-generating-positive}
 m_b:=(4C_bD_{\rm c})^{-1}
       \le \widetilde W_{sT}(s,T)\le C_eD_{\rm c}^{8}.
\end{equation}
We keep the three constants $c_b,c_h,c_1$ distinct: they govern the
boundary lower bound, the first collar width, and the extension width,
respectively.

\textit{5. The collar map, its inverse, and the boundary jet.}
On $\mathcal C=\T\times(-h_1,h_1)$ set
\[
 \Phi_{\rm src}(s,T)=(s,\widetilde W_s(s,T)),\qquad
 \Phi_{\rm tar}(s,T)=(\widetilde W_T(s,T),T).
\]
With $a=\widetilde W_{ss}$, $b=\widetilde W_{sT}$ and
$c=\widetilde W_{TT}$ we have
\[
 D\Phi_{\rm src}=\begin{pmatrix}1&0\\a&b\end{pmatrix},\qquad
 D\Phi_{\rm tar}=\begin{pmatrix}b&c\\0&1\end{pmatrix},\qquad b\ge m_b>0.
\]
The source map is one-to-one because $T\mapsto\widetilde W_s(s,T)$
is increasing.  The target map is one-to-one and onto the whole
$T$-cylinder, since for fixed $T$ its angular lift has positive
derivative $b$ and satisfies
\[
 \widetilde W_T(s+1,T)=\widetilde W_T(s,T)+1.
\]
On the image of $\Phi_{\rm src}$ we now define
\begin{equation}\label{eq:collar-factorization}
 \widetilde g=\Phi_{\rm tar}\circ\Phi_{\rm src}^{-1},\qquad
 t=\widetilde W_s(s,T),\quad
 \widetilde\Sigma=\widetilde W_T(s,T).
\end{equation}
This makes $\widetilde g$ a single injective $C^{1,\beta}$ collar map,
with differential in source variables
\begin{equation}\label{eq:extension-jacobian}
 D\widetilde g=
 \begin{pmatrix}b-ac/b&c/b\\-a/b&1/b\end{pmatrix},
 \qquad \det D\widetilde g=1.
\end{equation}
Equivalently,
$ds\wedge dt=b\,ds\wedge dT=d\widetilde\Sigma\wedge dT$.
Furthermore,
\[
 t=\int_0^T\widetilde W_{sT}(s,u)\,du
 \quad\Longrightarrow\quad \operatorname{sign}t=\operatorname{sign}T.
\]
For $T\ge0$ the defining equations are those of $g$, and at $T=0$
\[
 a=0,\qquad b=T_t(s,0)^{-1},\qquad c=\Sigma_t(s,0)/T_t(s,0),
\]
so that \eqref{eq:extension-jacobian} coincides with $Dg(s,0)$.

Thus $\widetilde g$ and $g$ have the same first-order
boundary jet. Since $\widetilde W\in C^{2,\beta}$
across $T=0$ and $b\ge m_b>0$, the glued map
is $C^{1,\beta}$ across the seam.

\begingroup
\setlength{\abovedisplayskip}{3pt plus 1pt minus 1pt}
\setlength{\belowdisplayskip}{3pt plus 1pt minus 1pt}
\setlength{\abovedisplayshortskip}{2pt plus 1pt minus 1pt}
\setlength{\belowdisplayshortskip}{2pt plus 1pt minus 1pt}
\setlength{\jot}{2pt}
\textit{6. Global Lipschitz bounds and the final polynomial estimate.}
Since the image of $\Phi_{\rm src}$ need not be a product cylinder, we
estimate the inverses directly, without any convexity assumption on
that image.  For two inverse source points with mixed coordinates
$(s_i,T_i)$,
\[
 m_b|T_1-T_2|
 \le |t_1-t_2|+C_eD_{\rm c}^{8}|s_1-s_2|.
\]
For two inverse target points,
\[
 m_b|s_1-s_2|
 \le |\widetilde\Sigma_1-\widetilde\Sigma_2|
                +C_eD_{\rm c}^{8}|T_1-T_2|.
\]
We work in compatible real lifts; periodicity then gives
the same local estimates on the cylinder, so both inverse
maps have Lipschitz constant at most $CD_{\rm c}^{9}$.
These bounds supply the Lipschitz control needed in the
subsequent coordinate changes. From the matrices above and
\eqref{eq:inverse-holder-calculus},
\begin{samepage}
	\[
	\begin{array}{@{}l@{}}
		\displaystyle
		\|D\Phi_{\rm src}\|_0+\|D\Phi_{\rm tar}\|_0
		+[D\Phi_{\rm src}]_\beta+[D\Phi_{\rm tar}]_\beta
		\le CD_{\rm c}^{8},\\[2pt]
		\displaystyle
		\|D\Phi_{\rm src}^{-1}\|_0+\|D\Phi_{\rm tar}^{-1}\|_0
		\le CD_{\rm c}^{9},\\[2pt]
		\displaystyle
		[D\Phi_{\rm src}^{-1}]_\beta+[D\Phi_{\rm tar}^{-1}]_\beta
		\le CD_{\rm c}^{18+8+9\beta}
		\le CD_{\rm c}^{35}.
	\end{array}
	\]
\end{samepage}
Composing as in \eqref{eq:collar-factorization}, and repeating the
computation with source and target interchanged, we get
\[
\begin{array}{@{}l@{}}
	\displaystyle
	\|D\widetilde g^{\pm1}\|_0\le CD_{\rm c}^{17},\\[2pt]
	\displaystyle
	[D\widetilde g^{\pm1}]_\beta
	\le C\bigl(D_{\rm c}^{8+9\beta+9}+D_{\rm c}^{8+35}\bigr)
	\le CD_{\rm c}^{43}.
\end{array}
\]
Returning to the fixed atlas adds, by \eqref{eq:atlas-change}, at most
$C(D_{\rm c}^{43}+D_{\rm c}^{17(1+\beta)}+D_{\rm c}^{17})$ to these
bounds.
The derivative traces agree across the seam.
Put $r_*:=m_bh_1$, so that
$r_*\ge cD_{\rm c}^{-(1+9/\beta)}$.
For two points on opposite sides of the seam,
let $h$ denote their distance in compatible
fixed local coordinates.
If $h\le c_0r_*$, with $c_0>0$ depending
only on the fixed geometry, the joining
segment stays in a common two-sided collar
chart.
Splitting this segment at the corresponding seam
($t=0$ for $\widetilde g$ and $T=0$ for
$\widetilde g^{-1}$) and using the matching
derivative traces gives the required H\"older estimate.

If $h>c_0r_*$, the uniform derivative
bound instead gives a H\"older quotient
bounded by
\[
CD_{\rm c}^{17}r_*^{-\beta}
\le CD_{\rm c}^{26+\beta}
\le CD_{\rm c}^{45}.
\]
  On the original interior, the inverse-matrix identity and the
area-compatible bound on the inverse norm give
$[Dg^{-1}]_\beta\le C(D^{3+\beta}+D^2)$.  The glued map and its
inverse therefore satisfy the common estimate
\begin{equation}\label{H:eq:extension-final-polynomial}
 \|D\widetilde g^{\pm1}\|_0+[D\widetilde g^{\pm1}]_\beta
          \le CD_{\rm c}^{45}\le C D^{90}.
\end{equation}
The finitely many source collars are pairwise disjoint, as are their
target collars, and $\pi$ is a permutation.  By side preservation, the
exterior images are disjoint from $S$.  The collar maps therefore patch
injectively with $g$ on $S$, and no further global injectivity
assumption is needed.

Finally, the source image contains $|t|<m_bh_1$ and the target
contains $|T|<h_1$, so both exterior widths are at least
\[
 cD_{\rm c}^{-(1+9/\beta)}\ge c'D^{-(2+18/\beta)}.
\]
We take
\[
 m_{\rm ext}=100+20/\beta
       >\max\{90,\,2+18/\beta\}
\]
and then choose $C_{\rm ext}$ large enough to dominate all the fixed
constants and the reciprocals of the collar constants.  Equations
\eqref{eq:extension-coordinate-budget}--\eqref{H:eq:extension-final-polynomial}
now give all four assertions, uniformly over the finitely many collar
pairs.
\endgroup
\end{proof}

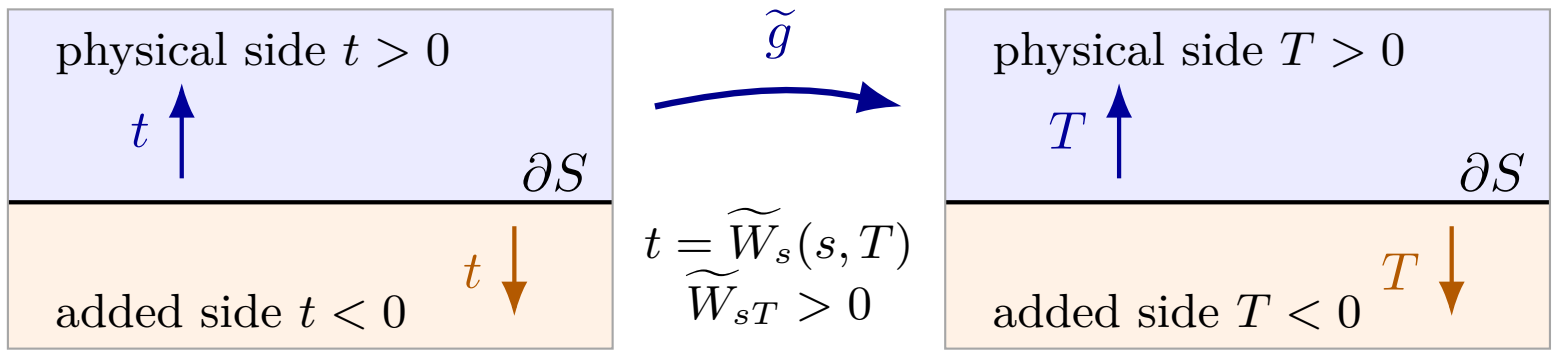
\begin{figure}[H]
\centering
\begin{tikzpicture}[x=1cm,y=.88cm,>=Latex,font=\small]
  \path[use as bounding box] (-0.1,-1.55) rectangle (10.3,1.55);
  \fill[blue!8] (0,0) rectangle (4,1.45);
  \fill[orange!10] (0,-1.1) rectangle (4,0);
  \draw[thick] (0,0)--(4,0);
  \draw[gray!70] (0,-1.1) rectangle (4,1.45);
  \node[anchor=west,font=\footnotesize] at (0.18,1.12)
       {physical side $t>0$};
  \node[anchor=west,font=\footnotesize] at (0.18,-0.82)
       {added side $t<0$};
  \node[anchor=south east,fill=blue!8,inner sep=1pt]
       at (3.88,0.04) {$\partial S$};
  \draw[->,thick,blue!60!black] (1.15,0.18)--(1.15,0.90)
       node[midway,left=2pt] {$t$};
  \draw[->,thick,orange!70!black] (3.35,-0.18)--(3.35,-0.86)
       node[midway,left=2pt] {$t$};
  \node[font=\footnotesize] at (2,-1.48) {(a) source collar};

  \draw[->,very thick,blue!55!black] (4.28,0.72) to[bend left=12]
       node[above=2pt] {$\widetilde g$} (5.92,0.72);

  \fill[blue!8] (6.2,0) rectangle (10.2,1.45);
  \fill[orange!10] (6.2,-1.1) rectangle (10.2,0);
  \draw[thick] (6.2,0)--(10.2,0);
  \draw[gray!70] (6.2,-1.1) rectangle (10.2,1.45);
  \node[anchor=west,font=\footnotesize] at (6.38,1.12)
       {physical side $T>0$};
  \node[anchor=west,font=\footnotesize] at (6.38,-0.82)
       {added side $T<0$};
  \node[anchor=south east,fill=blue!8,inner sep=1pt]
       at (10.08,0.04) {$\partial S$};
  \draw[->,thick,blue!60!black] (7.35,0.18)--(7.35,0.90)
       node[midway,left=2pt] {$T$};
  \draw[->,thick,orange!70!black] (9.55,-0.18)--(9.55,-0.86)
       node[midway,left=2pt] {$T$};
  \node[font=\footnotesize] at (8.2,-1.48) {(b) image collar};

  \node[align=center,font=\footnotesize] at (5.1,-0.48)
       {$t=\widetilde W_s(s,T)$\\[-1pt]
        $\widetilde W_{sT}>0$};
\end{tikzpicture}

\caption{\small The polynomial local extension near a
	boundary seam. Since
	$\widetilde W_s(s,0)=0$ and
	$\widetilde W_{sT}>0$, the relation
	$t=\widetilde W_s(s,T)$ preserves the two sides
	of the seam. For each fixed $T$, the map
	$s\mapsto\widetilde\Sigma=\widetilde W_T$
	is a degree-one circle diffeomorphism by
	positivity and the periodic lift relation.
	The resulting map
	$\widetilde g=\Phi_{\rm tar}\circ\Phi_{\rm src}^{-1}$
	is symplectic.}
\label{H:fig:symplectic-collar}
\end{figure}

\enlargethispage{1.0cm}
\Needspace{4\baselineskip}
The two triangular coordinate changes in
\eqref{eq:collar-factorization} are illustrated in
\cref{H:fig:symplectic-collar}.

\section{H\"older moving boxes}\label{H:sec:moving-boxes}

All estimates in this section are local.
Let $D$ be the original regularity bound in
\eqref{H:eq:finite-data}. We define the
regularity bound used in the moving boxes by
\[
\mathfrak D=
\begin{cases}
	D, & \partial S=\varnothing,\\
	\Delta=C_{\rm ext}D^{m_{\rm ext}},
	& \partial S\ne\varnothing,
\end{cases}
\]
where the extension in the boundary case is
provided by \cref{H:lem:symplectic-extension}.

We write $\widetilde g=g$ in the boundaryless
case and use the local symplectic extension
$\widetilde g$ in the boundary case.
All local scale and regularity estimates below
are expressed in terms of $\mathfrak D$.
Since $\widetilde g$ agrees with $g$ on $S$,
the derivative rates along the selected orbit
remain unchanged, and $a=\log A$ retains
its original meaning.

We assume that
$\mathfrak D\ge A\ge A_*(S,\beta)$,
with $A_*$ sufficiently large that
$\delta=\rho\log A\ge1$.

The moving coordinates are normalized as follows.  Put
\[
 i_n:\R^2\longrightarrow T_{x_n}\widehat S,\qquad
 i_n(v,w)=vv_n^u+wv_n^s,\qquad
 \Phi_n=\exp_{x_n}\circ i_n.
\]
Here $\exp_{x_n}$ is the exponential map of the
fixed ambient Riemannian metric
$\widehat{\mathfrak g}$.

For a closed surface, $\widehat S$ can be replaced by $S$.
Both frame vectors have unit length, so
\[
 \|i_n\|\le\sqrt2,\qquad
 \|i_n^{-1}\|\le C_S(1+\cot\angle(v_n^u,v_n^s)).
\]
Throughout, $(v,w)$ are coefficients with respect to this frame, so
that a coordinate box $U_n\subset\R^2$ stands for the physical set
$\Phi_n(U_n)$. 
Let $K=100$, and let
\begin{equation}\label{H:eq:psi}
 \psi_p(t)=
 \begin{cases}
  t,&t\le0,\\
  t/p,&t>0.
 \end{cases}
\end{equation}
Starting from $c_0=1$, define recursively
\begin{equation}\label{H:eq:c-recurrence}
 c_{n+1}=\min\{e^{\psi_p(\lambda_n^e-\delta)}c_n,K\}.
\end{equation}

Fix $M=1000$ and introduce the base scales
\begin{equation}\label{H:eq:base-scales}
	\bar r=\mathfrak D^{-5M/\beta},\qquad
	\bar\kappa=\mathfrak D^{-M},\qquad
	\bar b=\mathfrak D^M.
\end{equation}
In this section $D_0=D_0(S,\beta)\ge1$ denotes
a threshold for $\mathfrak D$ that may be
increased a finite number of times.
In the proof of \cref{thm:finite-closing},
we choose $A_*\ge D_0$.
Since $\mathfrak D\ge D\ge A\ge A_*$,
this ensures $\mathfrak D\ge D_0$.

The box and cone parameters are then
\begin{equation}\label{H:eq:box-parameters}
	\begin{split}
		r_n&=\bar r c_n^p,\qquad
		\tau_n=\bar r\exp\!\left(
		\sum_{i=0}^{n-1}(\lambda_i^s+\delta)
		\right),\\
		\kappa_n&=\bar\kappa c_n^{-1},\qquad
		\widetilde\kappa_n=\bar\kappa c_n,\qquad
		b_n=\bar b c_n^{-1}.
	\end{split}
\end{equation}

\paragraph{The coordinate boxes and their four sides.}\label{def:boxes-and-strips}
For $r,\tau,\kappa>0$ let
\[
 U(r,\tau,\kappa)=\{(v,w):|v|\le r,\ |w|\le\tau+\kappa|v|\}.
\]
Its four distinguished sides are
\begin{align*}
 \partial_v^\pm U&=\{(\pm r,w):|w|\le\tau+\kappa r\},\\
 \partial_h^\pm U&=\{(v,\pm(\tau+\kappa|v|)):|v|\le r\}.
\end{align*}
The upper side is thus V-shaped and the lower side is its mirror
image; the vertical half-height is smallest at $v=0$, where it
equals $\tau$.  The homeomorphism
\[
 (x,y)\longmapsto(rx,y(\tau+\kappa r|x|)),\qquad (x,y)\in[-1,1]^2,
\]
identifies these four sides with the sides of a square.
Thus $U(r,\tau,\kappa)$ is a topological
rectangle, star-shaped about the origin,
although it need not be convex.
We write $U_n=U(r_n,\tau_n,\kappa_n)$.

For a slope $h>0$ the \emph{open directional cones} are
\[
 \cC(h)=\{\xi\in\R^2:|\xi_2|<h|\xi_1|\},\qquad
 \tcC(h)=\{\xi\in\R^2:|\xi_1|<h|\xi_2|\}.
\]
We set $\cC_n=\cC(\kappa_n)$ and $\tcC_n=\tcC(\widetilde\kappa_n)$.
These are cones of tangent vectors and are not regions in the box;
their closures contain the zero vector.

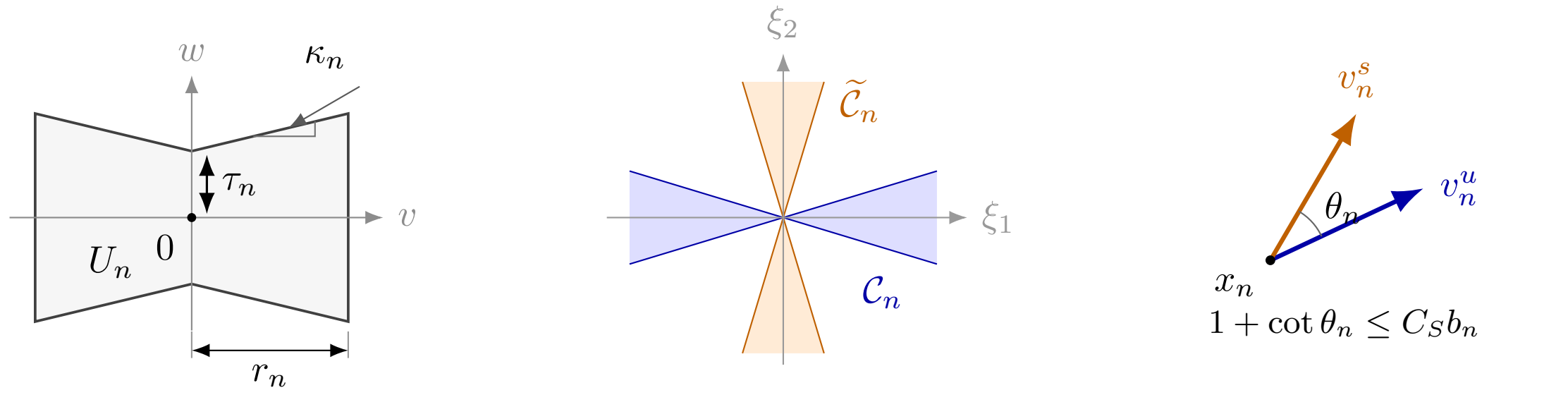
\begin{figure}[!htbp]
\centering
\begin{tikzpicture}[x=.89cm,y=.89cm,>=Latex,font=\small]
 \begin{scope}
  \filldraw[fill=gray!7,draw=black!75,line width=.7pt]
   (-1.65,-1.096)--(0,-.7)--(1.65,-1.096)--(1.65,1.096)
   --(0,.7)--(-1.65,1.096)--cycle;
  \draw[->,black!45] (-1.92,0)--(2.02,0) node[right] {$v$};
  \draw[->,black!45] (0,-1.19)--(0,1.50) node[above] {$w$};
  \fill (0,0) circle (1.2pt);
  \node[anchor=north east] at (-.03,-.04) {$0$};
  \draw[<->,line width=.6pt] (.16,.035)--(.16,.665)
        node[midway,right] {$\tau_n$};
  \draw[black!40] (0,-1.2)--(0,-1.48);
  \draw[black!40] (1.65,-1.2)--(1.65,-1.48);
  \draw[<->] (0,-1.40)--(1.65,-1.40)
        node[midway,below] {$r_n$};
  \draw[black!55] (.65,.856)--(1.30,.856)--(1.30,1.012);
  \draw[->,black!65] (1.77,1.38)--(1.02,.945);
  \node[anchor=south east] at (1.78,1.42) {$\kappa_n$};
  \node at (-.85,-.47) {$U_n$};
  \node[font=\footnotesize,align=center] at (.04,-2.43)
       {(a) Coordinate box\\$|w|\le\tau_n+\kappa_n|v|$};
 \end{scope}
 \begin{scope}[xshift=5.55cm]
  \fill[blue!13] (0,0)--(1.62,.49)--(1.62,-.49)--cycle;
  \fill[blue!13] (0,0)--(-1.62,.49)--(-1.62,-.49)--cycle;
  \fill[orange!16] (0,0)--(.43,1.43)--(-.43,1.43)--cycle;
  \fill[orange!16] (0,0)--(.43,-1.43)--(-.43,-1.43)--cycle;
  \draw[blue!65!black] (-1.62,-.49)--(1.62,.49);
  \draw[blue!65!black] (-1.62,.49)--(1.62,-.49);
  \draw[orange!75!black] (-.43,-1.43)--(.43,1.43);
  \draw[orange!75!black] (-.43,1.43)--(.43,-1.43);
  \draw[->,black!40] (-1.86,0)--(1.94,0) node[right] {$\xi_1$};
  \draw[->,black!40] (0,-1.55)--(0,1.73) node[above] {$\xi_2$};
  \node[blue!65!black] at (1.05,-.79) {$\cC_n$};
  \node[orange!75!black] at (.81,1.23) {$\tcC_n$};
  \node[font=\footnotesize,align=center] at (0,-2.43)
       {(b) Tangent-vector cones\\slopes $\kappa_n,\widetilde\kappa_n$};
 \end{scope}
 \begin{scope}[xshift=10.12cm]
  \coordinate (o) at (0,-.45);
  \draw[->,very thick,blue!65!black] (o)--(1.62,.31)
        node[anchor=west] {$v_n^u$};
  \draw[->,very thick,orange!75!black] (o)--(.91,1.10)
        node[anchor=south] {$v_n^s$};
  \draw[black!60] (.54,-.198) arc[start angle=25,end angle=59.6,radius=.60];
  \node at (.77,.10) {$\theta_n$};
  \fill (o) circle (1.3pt);
  \node[anchor=north east] at (o) {$x_n$};
  \node[font=\footnotesize,align=center] at (.78,-1.12)
       {$1+\cot\theta_n\le C_Sb_n$};
  \node[font=\footnotesize,align=center] at (.80,-2.43)
       {(c) Physical moving frame\\$b_n$ controls its degeneracy};
 \end{scope}
\end{tikzpicture}
\caption{The moving-box parameters.  Panel (a) shows the actual
V-shaped upper side and its reflected lower side: $r_n$ is the
horizontal half-width, $\tau_n$ is the half-height at $v=0$, and the
outer half-height is $\tau_n+\kappa_nr_n$.  Panel (b) is a separate
vector space: $\kappa_n$ and $\widetilde\kappa_n$ bound the horizontal
and vertical cone slopes.  Panel (c) uses unit frame vectors and the
unoriented angle $\theta_n=\angle(v_n^u,v_n^s)$; the inverse-angle
bound is proved in \cref{H:lem:scale-recurrence}.  The parameter $b_n$
is not a side length of $U_n$.  Dimensions and cone apertures are
enlarged for visibility.}
\label{H:fig:moving-box-parameters}
\end{figure}

\Cref{H:fig:moving-box-parameters} shows side by side the dimensions
of the box, the tangent-vector cones, and the physical moving frame.

An $h$-horizontal graph is a set $\{(v,f(v)):v\in J\}$ with
$\Lip(f)\le h$.  It is \emph{full in $U$} if $J=[-r,r]$ and the graph
lies in $U$, joining $\partial_v^-U$ to $\partial_v^+U$.
A full $k$-vertical graph is a set
$\{(a(w),w):w\in[w_-,w_+]\}\subset U$ with $\Lip(a)\le k$ whose
endpoints lie on $\partial_h^-U$ and $\partial_h^+U$.  The interval
$[w_-,w_+]$ is part of the data of the graph; in a nonrectangular box
it need not be $[-\tau,\tau]$. We write
$\partial_h U=\partial_h^-U\cup\partial_h^+U$
and
$\partial_v U=\partial_v^-U\cup\partial_v^+U$.

A \emph{horizontal strip} is the closed region bounded by two disjoint
full horizontal graphs, neither of which meets $\partial_hU$, together
with the corresponding subarcs of $\partial_vU$.  A \emph{vertical
strip} is the analogous region with the roles reversed: its two
internal sides are full vertical graphs disjoint from $\partial_vU$,
and its horizontal sides are subarcs of $\partial_hU$.  A homeomorphism
$F:V\to H$ is \emph{regular} when
\[
 F(\partial_hV)=\partial_hH,\qquad
 F(\partial_vV)=\partial_vH,
\]
each connected side being mapped onto a single connected side of the
same type; the two sides of a given type may be exchanged.  Regularity
is a condition on the distinguished \emph{boundary arcs} alone and does
not say that every horizontal curve is mapped horizontally.

\subsection{Scale identities}

\Needspace{6\baselineskip}
\begin{lemma}\label{H:lem:scale-recurrence}
For $0\le n<L$,
\begin{align}
 r_{n+1}&<e^{\lambda_n^u-\delta/2}r_n,
                                                \label{H:eq:radius-crossing}\\
 \kappa_{n+1}&\ge
   \max\{e^{-\lambda_n^e+\delta}\kappa_n,\bar\kappa/K\},
                                                \label{H:eq:kappa-rec}\\
 \widetilde\kappa_{n+1}&\le
      e^{\lambda_n^e-\delta}\widetilde\kappa_n,
                                                \label{H:eq:tkappa-rec}\\
 b_{n+1}&\ge
   \max\{e^{-\lambda_n^e+\delta}b_n,\bar b/K\}.
                                                \label{H:eq:b-rec}
\end{align}
Moreover, for $0\le n\le L$, $b_n$ bounds the inverse angle:
\begin{equation}\label{H:eq:angle-bound}
             1+\cot\angle(v_n^u,v_n^s)\le C_Sb_n.
\end{equation}
\end{lemma}

\begin{proof}
Write $s_n=\lambda_n^e-\delta$.  If $s_n>0$, then \eqref{H:eq:psi}
gives $p\psi_p(s_n)=s_n\le\lambda_n^u-\delta$, while if $s_n\le0$, then
$p\psi_p(s_n)=ps_n\le s_n\le\lambda_n^u-\delta$. Since $\delta>0$, this proves
\eqref{H:eq:radius-crossing}.

Since $\psi_p(s)\le s$, taking reciprocals
in the recurrence for $c_n$ gives the
bounds for $\kappa_{n+1}$ and $b_{n+1}$.
The bound for $\widetilde\kappa_{n+1}$
follows directly from
$c_{n+1}\le e^{\psi_p(s_n)}c_n$.
Because $\lambda_n^e\le-2\lambda_n^s$, we have in particular
\[
 b_{n+1}\ge e^{2\lambda_n^s+\delta}b_n.
\]

To prove \eqref{H:eq:angle-bound}, write
$\theta_n=\angle(v_n^u,v_n^s)$.
The initial bound follows from
\eqref{H:eq:endpoint-angle}, since
\[
1+\cot\theta_0\le 1+2A^3
\le C_S\bar b=C_Sb_0
\]
once $A_*$ is sufficiently large, using
$\bar b=\mathfrak D^M\ge A^M$.

For the induction, enlarge $A_*$ so that
$\delta\ge1$ and $\bar b/K\ge10A^2$.
Put $x=e^{2\lambda_n^s}b_n$. Then
\[
x+A^2\le\max\{e^\delta x,\bar b/K\}
\le b_{n+1}.
\]
Indeed, for $x\le\bar b/(2K)$ the first
inequality follows from the second term of
the maximum, while for $x>\bar b/(2K)$
it follows from $A^2<x/5$ and $\delta\ge1$.

Assuming $1+\cot\theta_n\le C_Sb_n$,
the recurrence \eqref{H:eq:cot-recurrence}
therefore gives, with $C_S\ge2$,
\[
\begin{aligned}
	1+\cot\theta_{n+1}
	&\le C_Sx+A^2+1\\
	&\le C_S(x+A^2)\\
	&\le C_Sb_{n+1}.
\end{aligned}
\]
This completes the induction.
\end{proof}

\Needspace{5\baselineskip}
These inequalities also control the aspect ratios of the moving boxes and give the endpoint scales needed for the return.

\begin{lemma}\label{H:lem:box-shape}
After increasing $A_*(S,\beta)$ and $D_0(S,\beta)$,
\begin{equation}\label{H:eq:shape-bounds}
 r_n\ge\tau_n\quad(0\le n\le L),\qquad
 c_L=K,\qquad \tau_L\le\frac1{10}\bar r.
\end{equation}
Consequently
\begin{equation}\label{H:eq:end-scales}
 r_L=K^p\bar r,\qquad
 \kappa_L=\bar\kappa/K,\qquad
 \widetilde\kappa_L=K\bar\kappa.
\end{equation}
All $U_n$ lie in coordinate balls of radius
$O(\mathfrak D^{-1})$.
\end{lemma}

\begin{proof}
For $n=0$ we have $r_0=\tau_0=\bar r$, so in the strict estimates
below we take $n\ge1$.  The ingredients are the prefix and suffix
bounds of \cref{H:lem:effective-rates} and the recurrences of
\cref{H:lem:scale-recurrence}.  From \eqref{H:eq:psi}, for every
$s\in\R$,
\begin{equation}\label{H:eq:psi-identity}
 \psi_p(s)=\frac{s}{p}+\left(1-\frac1p\right)\min\{s,0\}.
\end{equation}
Moreover, the truncation in \eqref{H:eq:c-recurrence} never decreases
$c_n$ on a positive step, so that
\[
 \log c_n\ge\sum_{i<n}\min\{\lambda_i^e-\delta,0\}.
\]
In view of \eqref{H:eq:negative-part} this implies
\[
 p\log c_n\ge-p(8\eta+\rho)na>-\frac1{10}na,
\]
whereas \eqref{H:eq:stable-average} gives
\[
 \log(\tau_n/\bar r)\le-(1-\eta-\rho)na<-\frac9{10}na.
\]
Hence $r_n=\bar r c_n^p\ge\tau_n$.  The same stable estimate, applied
with $n=L\ge1$, yields $\tau_L\le\bar r/10$ for $A_*$ large.

It remains to show that $c_L=K$.  On every terminal block
$J=\{j,\ldots,L-1\}$, \eqref{H:eq:effective-average},
\eqref{H:eq:negative-part} and the identity \eqref{H:eq:psi-identity}
give
\begin{align}
 \sum_{i=j}^{L-1}\psi_p(\lambda_i^e-\delta)
 &\ge\left[
 \frac{1-4\eta-\rho}{p}
 -\left(1-\frac1p\right)(8\eta+\rho)
 \right](L-j)a\notag\\
 &\ge\frac{L-j}{2p}a.\label{H:eq:suffix-growth}
\end{align}
Now let $j$ be the largest index in $\{0,\ldots,L-1\}$ with $c_j\ge1$.
By maximality of $j$, the truncation at $K$ cannot be active on
$\{j,\ldots,L-1\}$ before the last step, so \eqref{H:eq:suffix-growth}
gives
\[
 c_L=\min\left\{K,
 c_j\exp\left(\sum_{i=j}^{L-1}
       \psi_p(\lambda_i^e-\delta)\right)\right\}=K
\]
as soon as $A_*^{1/(2p)}\ge K$.

Finally,
\[
r_n+\tau_n+\kappa_nr_n
\le 2K^p\bar r+K^{p-1}\bar\kappa\bar r.
\]
Since $\bar r=\mathfrak D^{-5M/\beta}$,
the right-hand side is smaller than
$\mathfrak D^{-1}$ times the fixed coordinate
radius once $D_0$ has been increased.
\end{proof}

\subsection{The H\"older remainder}
With the normalized maps
$\Phi_n=\exp_{x_n}\circ i_n$ fixed above, set
\[
g_n=\Phi_{n+1}^{-1}\circ\widetilde g\circ\Phi_n
=i_{n+1}^{-1}\circ\exp_{x_{n+1}}^{-1}
\circ\widetilde g\circ\exp_{x_n}\circ i_n.
\]
Here $\widetilde g=g$ in the boundaryless case,
as agreed at the beginning of this section.
The domain check in
\cref{H:lem:coordinate-holder} shows that $g_n$
is defined on a neighborhood of $U_n$.
We do not claim that $g_n(U_n)\subset U_{n+1}$.
Since $\widetilde g$ agrees with $g$ along
the selected orbit, we have
\begin{equation}\label{H:eq:coordinate-map}
	\begin{split}
		g_n(v,w)&=
		\bigl(e^{\lambda_n^u}v+F_n(v,w),
		e^{\lambda_n^s}w+G_n(v,w)\bigr),\\
		F_n(0)&=G_n(0)=0,\qquad
		DF_n(0)=DG_n(0)=0.
	\end{split}
\end{equation}

\Needspace{5\baselineskip}

\begin{lemma}\label{H:lem:coordinate-holder}
	There is $C_{S,\beta}\ge1$ such that, for
	$0\le n<L$, on $U_n$,
	\begin{equation}\label{H:eq:coordinate-holder}
		[Dg_n]_\beta
		\le C_{S,\beta}\mathfrak D^2b_{n+1}.
	\end{equation}
	Consequently, if
	\begin{equation}\label{H:eq:epsilon-def}
		\varepsilon_n=
		C_{S,\beta}\mathfrak D^2b_{n+1}
		\bigl(2r_n(1+\kappa_n)\bigr)^\beta,
	\end{equation}
	then
	\begin{equation}\label{H:eq:remainder-derivative}
		\|DF_n\|_{C^0(U_n)}
		+\|DG_n\|_{C^0(U_n)}
		\le\varepsilon_n.
	\end{equation}
\end{lemma}

\begin{proof}
	Put $f=\widetilde g$.
	By the definition of $\mathfrak D$ and
	\cref{H:lem:symplectic-extension},
	\[
	\|Df\|_0\le\mathfrak D,\qquad
	[Df]_\beta\le\mathfrak D
	\]
	on the relevant fixed coordinate domains.

	We first verify the domains of the coordinate
	maps. In the boundary case,
	\cref{H:lem:symplectic-extension} provides
	exterior collars of width at least
	$\mathfrak D^{-1}$. By the fixed collar
	geometry, there is $c_S>0$ such that
	every $x\in S$ has an ambient ball of radius
	$c_S\mathfrak D^{-1}$ contained in
	$\mathcal N_{\rm src}$, once $\mathfrak D$
	is sufficiently large.
	In the boundaryless case no extension-domain
	restriction is needed.

	By \cref{H:lem:box-shape} and
	$\|i_n\|\le\sqrt2$, the physical boxes have
	radius at most
	\[
	C_{S,\beta}\mathfrak D^{-5M/\beta}.
	\]
	Since $5M/\beta>2$, after increasing $D_0$
	they lie, with a positive buffer, in ambient
	balls of radius $c_S\mathfrak D^{-2}$
	on which $f$ is defined.

	Let $\iota_S>0$ be a fixed uniform
	exponential-chart radius. Since
	$f(x_n)=x_{n+1}$ and $\|Df\|_0\le\mathfrak D$,
	the image of each physical box is contained
	in the ambient ball about $x_{n+1}$ of radius
	\[
	C_{S,\beta}\mathfrak D^{1-5M/\beta}
	<\iota_S/10
	\]
	once $D_0$ has been increased.
	Thus the source and target exponential
	coordinates in the definition of $g_n$
	are simultaneously valid on a neighborhood
	of $U_n$. In fact, after decreasing a fixed constant
	$c>0$, the same domain argument applies to
	the Euclidean ball $B(0,c\mathfrak D^{-2})$
	in the moving coordinates: its physical image
	lies in $\mathcal N_{\rm src}$, and its image
	under $f$ lies in the target exponential chart.
	For sufficiently large $D_0$, this convex ball
	contains $U_n$. We apply the coordinate-change
	estimates on this larger ball, without requiring
	$U_n$ itself to be convex.

	The coordinate comparisons now take place
	in fixed exponential charts with uniform
	smooth bounds, restricted to these admissible
	neighborhoods. The chain and product estimates
	of \cref{H:lem:chart-change} give
	\[
	\begin{aligned}
		&[D(\exp_{x_{n+1}}^{-1}
		\circ f\circ\exp_{x_n})]_\beta\\
		&\qquad\le C_{S,\beta}
		\bigl([Df]_\beta+
		\|Df\|_0^{1+\beta}+\|Df\|_0\bigr)\\
		&\qquad\le C_{S,\beta}\mathfrak D^2.
	\end{aligned}
	\]
	The term $\|Df\|_0^{1+\beta}$ comes from
	the variation of the derivative of the
	target exponential coordinate map.

	The linear maps $i_n$ have uniformly
	bounded norms, whereas
	\[
	\|i_{n+1}^{-1}\|
	\le C_S\bigl(1+\cot
	\angle(v_{n+1}^u,v_{n+1}^s)\bigr)
	\le C_Sb_{n+1}
	\]
	by \cref{H:lem:scale-recurrence}.
	Composing with $i_n$ and $i_{n+1}^{-1}$
	therefore gives
	\eqref{H:eq:coordinate-holder}.

	Finally, by \cref{H:lem:box-shape},
	every $z\in U_n$ satisfies
	\[
	\|z\|\le2r_n(1+\kappa_n).
	\]
	Since $DF_n(0)=DG_n(0)=0$,
	the H\"older estimate
	\eqref{H:eq:coordinate-holder}, with
	$C_{S,\beta}$ enlarged once if necessary,
	gives \eqref{H:eq:remainder-derivative}.
\end{proof}

\begin{lemma}\label{H:lem:four-errors}
	After increasing $D_0(S,\beta)$, for every
	$0\le n<L$,
	\begin{align}
		\varepsilon_n
		&\le \mathfrak D^{-3M}c_n,
		\label{H:eq:eps-one}\\
		\varepsilon_n(1+\kappa_n)
		&\le \mathfrak D^{-2M},
		\label{H:eq:eps-kappa}\\
		\frac{\varepsilon_n}{\widetilde\kappa_n}
		&\le \mathfrak D^{-2M},
		\label{H:eq:eps-over-tkappa}\\
		\varepsilon_n\widetilde\kappa_n
		&\le K^2\mathfrak D^{-4M}.
		\label{H:eq:eps-tkappa}
	\end{align}
\end{lemma}

\begin{proof}
	We estimate the nonlinear derivative error
	budget $\varepsilon_n$ defined in
	\eqref{H:eq:epsilon-def}.
	Since $\lambda_n^e\ge-2a$ and
	$\delta=\rho a$ with $\rho<1$, we have
	$\psi_p(\lambda_n^e-\delta)\ge-3a$.
	Using $c_n\le K$ in the recurrence
	\eqref{H:eq:c-recurrence} gives
	\[
	c_{n+1}\ge A^{-3}c_n,
	\qquad
	b_{n+1}\le A^3\bar b/c_n.
	\]

	Inserting this into \eqref{H:eq:epsilon-def}
	and using \eqref{H:eq:box-parameters},
	$p\beta=\beta+2$, and
	$(x+y)^\beta\le x^\beta+y^\beta$
	for $x,y\ge0$, we obtain
	\begin{align*}
		\varepsilon_n
		&\le C_{S,\beta}\mathfrak D^2 A^3
		\bar b\bar r^\beta
		\left(
		c_n^{p\beta-1}
		+\bar\kappa^\beta
		c_n^{(p-1)\beta-1}
		\right)\\
		&=C_{S,\beta}\mathfrak D^2 A^3
		\bar b\bar r^\beta
		\left(
		c_n^{\beta+1}
		+\bar\kappa^\beta c_n
		\right).
	\end{align*}

	Since $c_n\le K$, $A\le\mathfrak D$,
	$\bar\kappa\le1$, and
	\[
	\bar b\bar r^\beta
	=\mathfrak D^M\mathfrak D^{-5M}
	=\mathfrak D^{-4M},
	\]
	the right-hand side is bounded by
	\[
	C_{S,\beta}(K^\beta+1)
	\mathfrak D^{5-4M}c_n
	\le\mathfrak D^{-3M}c_n
	\]
	once $D_0$ is sufficiently large.
	This proves \eqref{H:eq:eps-one}.

	The remaining estimates follow directly:
	\[
	\begin{aligned}
		\varepsilon_n(1+\kappa_n)
		&\le\mathfrak D^{-3M}
		(c_n+\bar\kappa)
		\le(K+1)\mathfrak D^{-3M}
		\le\mathfrak D^{-2M},\\
		\frac{\varepsilon_n}{\widetilde\kappa_n}
		&\le\mathfrak D^{-3M}\bar\kappa^{-1}
		=\mathfrak D^{-2M},\\
		\varepsilon_n\widetilde\kappa_n
		&\le\mathfrak D^{-3M}\bar\kappa c_n^2
		\le K^2\mathfrak D^{-4M}.
	\end{aligned}
	\]
	The inverse cone estimate uses
	\eqref{H:eq:eps-over-tkappa}.
\end{proof}

\subsection{Cone and strip propagation}

\Needspace{6\baselineskip}
\begin{proposition}\label{H:prop:one-step}
	Under the standing conditions
	$\mathfrak D\ge A\ge A_*(S,\beta)$ and
	$\mathfrak D\ge D_0(S,\beta)$, the following
	hold for $0\le n<L$.
	\begin{enumerate}[label=\textnormal{(\roman*)}]
		\item For every $z\in U_n$,
		$Dg_n(z)\cC_n\subset\cC_{n+1}$ and
		$Dg_n(z)^{-1}\tcC_{n+1}\subset\tcC_n$.

		\item The image of every full
		$\kappa_n$-horizontal graph in $U_n$,
		intersected with $U_{n+1}$, is a full
		$\kappa_{n+1}$-horizontal graph.
		Moreover,
		\[
		H_{n+1}:=g_n(U_n)\cap U_{n+1}
		\]
		is a horizontal strip, and every arc
		$\Gamma=\{(\phi(w),w):w\in I\}\subset H_{n+1}$
		with $\operatorname{Lip}(\phi)
		\le\widetilde\kappa_{n+1}$
		joining its two horizontal sides pulls back
		to a full $\widetilde\kappa_n$-vertical
		graph in $U_n$.

		\item For $|w|\le\tau_n+\kappa_n r_n$,
		\[
		\pi_1g_n(r_n,w)>r_{n+1},\qquad
		\pi_1g_n(-r_n,w)<-r_{n+1}.
		\]
		Moreover,
		$g_n(U_n)\cap\partial_hU_{n+1}
		=\varnothing$.
	\end{enumerate}
	Consequently,
	\[
	V_n:=(g_n|_{U_n})^{-1}(H_{n+1})
	\]
	is a vertical strip of $U_n$, and
	$g_n:V_n\to H_{n+1}$ is regular.
\end{proposition}

\begin{proof}
Write $A_n=e^{\lambda_n^u}$, $B_n=e^{\lambda_n^s}$ and
$e_n=\varepsilon_n(1+\kappa_n)$.  By the remainder estimate,
\[
 Dg_n(z)=\begin{pmatrix}A_n&0\\0&B_n\end{pmatrix}+E_n(z),
 \qquad \|E_n(z)\|\le\varepsilon_n.
\]
All thresholds below are uniform in $n$;
we use freely $\delta\ge1$,
$A_n,B_n\ge\mathfrak D^{-1}$, and
\cref{H:lem:four-errors}.

\textit{1. Forward and inverse cones.}
For $\xi\in\overline{\cC_n}\setminus\{0\}$ we have
\begin{equation}\label{H:eq:forward-slope}
 \frac{|\pi_2Dg_n(z)\xi|}{|\pi_1Dg_n(z)\xi|}
 \le\frac{B_n\kappa_n+e_n}{A_n-e_n}.
\end{equation}
The numerator here is at most
$(B_n\kappa_n+e_n)|\xi_1|$ and the denominator at least
$(A_n-e_n)|\xi_1|$; since $e_n\le\mathfrak D^{-2M}$,
\[
A_n-e_n\ge A_ne^{-\delta/4}>0.
\]
Writing
\begin{equation}\label{H:eq:Xi}
 X_n:=\frac{B_n\kappa_n+e_n}{A_n-e_n},
\end{equation}
we obtain from \eqref{H:eq:kappa-rec} the bound
\[
X_n\le e^{-3\delta/4}\kappa_{n+1}
+\mathfrak D^{1+\rho/4-2M}.
\]
We choose $D_0$ so large that
$K\mathfrak D^{1+\rho/4-M}<1-e^{-3/4}$;
since $\kappa_{n+1}\ge\mathfrak D^{-M}/K$,
it follows that $X_n<\kappa_{n+1}$.

For the inverse cone, let $\bar\xi=Dg_n(z)\xi$ satisfy
$|\bar\xi_1|\le\widetilde\kappa_{n+1}|\bar\xi_2|$; eliminating
$\bar\xi$ gives
\begin{equation}\label{H:eq:inverse-slope}
 \frac{|\xi_1|}{|\xi_2|}
 \le\frac{\widetilde\kappa_{n+1}B_n+
                  \varepsilon_n(1+\widetilde\kappa_{n+1})}
              {A_n-\varepsilon_n(1+\widetilde\kappa_{n+1})}.
\end{equation}
Here the denominator is at least
$A_ne^{-\delta/4}$ for sufficiently large
$\mathfrak D$, as
$\varepsilon_n\le K\mathfrak D^{-3M}$ and
$\widetilde\kappa_{n+1}\le K\mathfrak D^{-M}$.
The main term is at most
$e^{-3\delta/4}\widetilde\kappa_n$
by \eqref{H:eq:tkappa-rec}, while
\eqref{H:eq:eps-over-tkappa} bounds the
error, divided by $\widetilde\kappa_n$, by
\[
\mathfrak D^{1+\rho/4-2M}
(1+K\mathfrak D^{-M})<1-e^{-3/4}.
\]
The right side of \eqref{H:eq:inverse-slope} is therefore strictly
smaller than $\widetilde\kappa_n$.  Both cone conclusions hold even for
nonzero vectors on the boundaries of the input cones.

\textit{2. Signed crossing and vertical containment.}
From the recurrence,
\begin{equation}\label{H:eq:kappa-upper}
 \kappa_{n+1}\le A^3\kappa_n.
\end{equation}
Since the box $U_n$ is star-shaped about $0$, we may integrate along
radial segments, and because the remainders vanish at $0$ this yields
\[
 |F_n(v,w)|,\ |G_n(v,w)|
 \le\varepsilon_n(|v|+|w|)
 \le\varepsilon_n((1+\kappa_n)|v|+\tau_n).
\]
Hence
\begin{align}
 |\pi_1g_n(z)|&\ge(A_n-e_n)|v|-\varepsilon_n\tau_n,
       \label{H:eq:first-coordinate-lower}\\
 |\pi_2g_n(z)|&\le(B_n\kappa_n+e_n)|v|+(B_n+\varepsilon_n)\tau_n.
       \label{H:eq:second-coordinate-upper}
\end{align}
Eliminating $|v|$ and using $X_n<\kappa_{n+1}$, we get
\[
 |\pi_2g_n(z)|\le\kappa_{n+1}|\pi_1g_n(z)|
       +(B_n+\varepsilon_n+\varepsilon_n\kappa_{n+1})\tau_n.
\]
The error stays below the available margin, since
\[
\varepsilon_n+\varepsilon_n\kappa_{n+1}
\le K\mathfrak D^{-3M}
+\mathfrak D^{3-2M}
<(e-1)\mathfrak D^{-1}
\le B_n(e^\delta-1).
\]
This gives
\begin{equation}\label{H:eq:vertical-containment}
 B_n+\varepsilon_n+\varepsilon_n\kappa_{n+1}<B_ne^\delta,
 \qquad
 |\pi_2g_n(z)|<\tau_{n+1}+\kappa_{n+1}|\pi_1g_n(z)|.
\end{equation}
For the crossing we need the sign of the first coordinate, which is
lost in \eqref{H:eq:first-coordinate-lower}.  Directly from
$\pi_1g_n(v,w)=A_nv+F_n(v,w)$ and $\tau_n\le r_n$ we have
\begin{align*}
 \pi_1g_n(r_n,w)&\ge(A_n-\varepsilon_n(2+\kappa_n))r_n
                 >A_ne^{-\delta/2}r_n>r_{n+1},\\
 \pi_1g_n(-r_n,w)&\le-(A_n-\varepsilon_n(2+\kappa_n))r_n
                 <-r_{n+1}.
\end{align*}
Here we used that
$\varepsilon_n(2+\kappa_n)
\le2\mathfrak D^{-2M}$
for sufficiently large $\mathfrak D$, while
\[
A_n(1-e^{-\delta/2})
\ge(1-e^{-1/2})\mathfrak D^{-1}.
\]

\textit{3. Full graphs and the strip between them.}
Let $w=\phi(v)$, $-r_n\le v\le r_n$,
be a full $\kappa_n$-horizontal graph in $U_n$.
  Then, almost everywhere,
\[
 \frac{d}{dv}\pi_1g_n(v,\phi(v))\ge A_n-e_n>0,
 \qquad
 \left|\frac{d\pi_2g_n}{d\pi_1g_n}\right|\le X_n<\kappa_{n+1}.
\]
Absolute continuity turns this into global monotonicity of the first
coordinate.  Together with the crossing of step~2 and
\eqref{H:eq:vertical-containment}, this shows that the part of the
image lying in $U_{n+1}$ is a full horizontal graph.

To describe the whole intersection, parameterize $U_n$ by
$Z(v,u)=(v,u(\tau_n+\kappa_n|v|))$ with $-1\le u\le1$.
For each $u$, the monotonicity just proved gives a unique $v=v(x,u)$
with $\pi_1g_n(Z(v,u))=x$ whenever $|x|\le r_{n+1}$, and
$v(x,u)$ is continuous.  For fixed $x$, the second
coordinate $u\mapsto\pi_2g_n(Z(v(x,u),u))$ is continuous, and injective
because $g_n$ is; its image is thus the interval between
the images of $u=-1$ and $u=1$.  In other words,
\[
 H_{n+1}:=g_n(U_n)\cap U_{n+1}
\]
is a single horizontal strip, bounded by the two full image graphs.

Since $g_n$ is a diffeomorphism on a
neighborhood of the connected box $U_n$
and $\det Dg_n(0)=A_nB_n>0$, we have
$\det Dg_n>0$ throughout $U_n$.

Since
$\kappa_{n+1}\widetilde\kappa_{n+1}
=\bar\kappa^2<1$,
a $\widetilde\kappa_{n+1}$-vertical graph
meets each horizontal side of $H_{n+1}$
at most once.

Now let $\Gamma=(\phi(w),w)\subset H_{n+1}$
be a $\widetilde\kappa_{n+1}$-vertical arc,
so that $|\phi'|\le\widetilde\kappa_{n+1}$
almost everywhere. In the inverse
matrix formula the vertical numerator is
$A_n+\partial_vF_n-\partial_vG_n\phi'>0$, by the same bound as for
the denominator in \eqref{H:eq:inverse-slope}, and $\det Dg_n>0$.
Hence the vertical coordinate of $g_n^{-1}(\Gamma)$ is strictly
increasing and its slope is at most $\widetilde\kappa_n$.
Since the two horizontal sides of $H_{n+1}$
are images of the corresponding subarcs of
$\partial_hU_n$, the endpoints of
$g_n^{-1}(\Gamma)$ lie on $\partial_hU_n$.
Applied to the two vertical sides of $H_{n+1}$, this gives two
inverse full graphs, which avoid $\partial_vU_n$ by the strict crossing
and bound
\[
 V_n:=(g_n|_{U_n})^{-1}(H_{n+1}).
\]
The restriction $g_n:V_n\to H_{n+1}$ is regular, and all the claims
are proved, every inverse map being taken on the branch indicated.
\end{proof}

The horizontal and vertical graph slopes satisfy
\begin{equation}\label{H:eq:cone-product}
 \kappa_n\widetilde\kappa_n=\bar\kappa^2<1.
\end{equation}
Let $w=f(v)$ be a full
$\kappa_n$-horizontal graph and
$v=a(w)$ a full
$\widetilde\kappa_n$-vertical graph in $U_n$,
with $w\in[w_-,w_+]$.
Put $h(w)=w-f(a(w))$.
Since the endpoints of the vertical graph
lie on the lower and upper horizontal
boundaries of $U_n$, respectively, we have
\[
h(w_-)\le0\le h(w_+).
\]
Moreover, for $w_1<w_2$,
\[
h(w_2)-h(w_1)
\ge(1-\kappa_n\widetilde\kappa_n)
(w_2-w_1)>0.
\]
Thus the two graphs meet exactly once.

Consequently, given two disjoint full
horizontal graphs and two disjoint full
vertical graphs, each horizontal--vertical
pair has exactly one intersection.
The subarcs joining these four intersection
points form a Jordan curve, whose enclosed
region is a topological quadrilateral with
the four intersections as its corners.

The next lemma concatenates the one-step
crossings into a regular branch from the
first box to the last, preserving both
cone inclusions.

\begin{lemma}\label{H:lem:strip-concatenation}
Let $G'=g_{L-1}\circ\cdots\circ g_0$
denote the local composition on its natural
domain. There is a
$\widetilde\kappa_0$-vertical strip $R'$ of $U_0$ such that $G'(R')$ is a
$\kappa_L$-horizontal strip $R$ of $U_L$ and $G'\colon R'\to R$ is
regular.  Moreover,
\begin{align}
 DG'(z)\cC_0&\subset\cC_L &&(z\in R'),
                                                        \label{H:eq:long-forward-cone}\\
 D(G'^{-1})_z\tcC_L&\subset\tcC_0 &&(z\in R).
                                                        \label{H:eq:long-inverse-cone}
\end{align}
\end{lemma}

\begin{proof}
	We run the concatenation argument of
	\cite[\S3.3.3, Corollary~D]{AFLXZ}
	with the one-step regular strip maps
	$g_n:V_n\to H_{n+1}$ of
	\cref{H:prop:one-step}.
	Write
	$\mathcal F_n=g_{n-1}\circ\cdots\circ g_0$;
	the calligraphic letter distinguishes this
	composition from the scalar remainder $F_n$
	in \eqref{H:eq:coordinate-map}.

	We construct inductively a
	$\widetilde\kappa_0$-vertical strip
	$V^{(n)}\subset U_0$ and a
	$\kappa_n$-horizontal strip
	$H^{(n)}\subset U_n$ such that
	\[
	\mathcal F_n:
	V^{(n)}\longrightarrow H^{(n)}
	\]
	is regular, every intermediate orbit point
	stays in its corresponding box, and
	\[
	\begin{aligned}
		D\mathcal F_n(z)\cC_0
		&\subset\cC_n
		&& (z\in V^{(n)}),\\
		D(\mathcal F_n^{-1})_y\tcC_n
		&\subset\tcC_0
		&& (y\in H^{(n)}).
	\end{aligned}
	\]
	We also retain preservation of the positive
	vertical half-cones by the inverse derivatives.
	For $n=1$ we take
	$V^{(1)}=V_0$ and $H^{(1)}=H_1$.

	Given the objects at time $n<L$, put
	\begin{equation}\label{eq:strip-induction}
		Q_n=H^{(n)}\cap V_n,\qquad
		V^{(n+1)}
		=(\mathcal F_n|_{V^{(n)}})^{-1}(Q_n),
		\qquad
		H^{(n+1)}=g_n(Q_n).
	\end{equation}
	By \eqref{H:eq:cone-product} and the
	full-graph intersection argument above,
	$Q_n$ is a single topological quadrilateral.
	Its horizontal sides are subarcs of
	$\partial_hH^{(n)}$, and its vertical sides
	are subarcs of $\partial_vV_n$.

	Each vertical side of $Q_n$ is a Lipschitz
	graph $\gamma(w)=(a(w),w)$ with
	$|a'(w)|\le\widetilde\kappa_n$ almost
	everywhere. Its positively oriented tangent
	vectors therefore lie in
	\[
	\{\xi\in\overline{\tcC_n}\setminus\{0\}:
	\xi_2>0\}.
	\]
	As shown by the inverse-matrix calculation
	in the proof of \cref{H:prop:one-step},
	each one-step inverse maps this closed
	positive vertical half-cone into the
	corresponding positive vertical cone at
	the preceding time. This property is
	preserved under composition.

	Consequently, the vertical coordinate
	of each preimage curve under
	$\mathcal F_n^{-1}$ is strictly increasing,
	and its slope is bounded by
	$\widetilde\kappa_0$. Thus the preimages
	are $\widetilde\kappa_0$-vertical graphs.

	The endpoints of these preimage graphs
	lie on opposite horizontal sides of
	$V^{(n)}$, hence on $\partial_hU_0$,
	by the regularity of $\mathcal F_n$.
	Moreover, the vertical sides of $Q_n$
	are disjoint from $\partial_vH^{(n)}$,
	since the vertical sides of $V_n$ are
	internal to $U_n$.
	Their preimages are therefore disjoint
	from $\partial_vV^{(n)}$ and hence from
	$\partial_vU_0$.
	Thus they bound a
	$\widetilde\kappa_0$-vertical strip
	$V^{(n+1)}$ of $U_0$.

	In the other direction, $g_n$ carries
	the horizontal sides of $Q_n$ to
	$\kappa_{n+1}$-horizontal graphs with
	strictly monotone first coordinate.
	By the regularity of $g_n:V_n\to H_{n+1}$,
	their endpoints lie on opposite vertical
	sides of $H_{n+1}$, hence on
	$\partial_vU_{n+1}$.
	Their images are disjoint from
	$\partial_hU_{n+1}$.
	Consequently they bound a
	$\kappa_{n+1}$-horizontal strip
	$H^{(n+1)}$ of $U_{n+1}$.

	The boundary correspondences show that
	\[
	\mathcal F_{n+1}:
	V^{(n+1)}\longrightarrow H^{(n+1)}
	\]
	is regular.
	For $z\in V^{(n+1)}$ we have
	$\mathcal F_n(z)\in Q_n\subset V_n$,
	so the next composition is defined.
	The chain rule then propagates both
	cone inclusions to time $n+1$.

	Finally, since $g_n(0)=0$ and
	$0\in U_n\cap U_{n+1}$, each one-step
	strip $V_n$ and $H_{n+1}$ contains $0$.
	The induction therefore gives
	\[
	0\in V^{(n)},\qquad
	0\in H^{(n)},\qquad
	\mathcal F_n(0)=0.
	\]

	At $n=L$ set
	$R'=V^{(L)}$ and $R=H^{(L)}$.
	Then $G'=\mathcal F_L$ on $R'$,
	and the iterated cone inclusions are
	\eqref{H:eq:long-forward-cone} and
	\eqref{H:eq:long-inverse-cone}.
\end{proof}

\Needspace{6\baselineskip}
A regular crossing of a vertical strip onto a
horizontal strip has a fixed point; together
with area preservation, strict projective cone
invariance makes each such fixed point
hyperbolic.

\begin{lemma}\label{H:lem:crossed-strip}
	Let $U\subset\mathbb R^2$ be a topological
	rectangle, let $V\subset U$ be a vertical strip,
	and let $H\subset U$ be a horizontal strip in
	the sense fixed above. If a homeomorphism
	$F\colon V\to H$ is regular, then $F$ has a
	fixed point.

	Suppose in addition that $F$ is differentiable
	at each fixed point and preserves a smooth
	area form there. Let $\mathcal C$ and
	$\widetilde{\mathcal C}$ be two open cones
	whose projectivizations
	\[
	I=\mathbb P(\mathcal C),\qquad
	J=\mathbb P(\widetilde{\mathcal C})
	\]
	are nondegenerate proper intervals in
	$\mathbb{RP}^1$ with disjoint closures.
	Assume that, at every fixed point $z$,
	\[
	\mathbb P(DF_z)(\overline I)\subset I,
	\qquad
	\mathbb P((DF_z)^{-1})(\overline J)\subset J,
	\]
	where $\mathbb P(B)$ denotes the projective
	map induced by an invertible linear map $B$.
	Then every fixed point of $F$ is hyperbolic.
\end{lemma}

\begin{proof}
	Choose a side-preserving homeomorphism
	$\chi:U\to[-1,1]^2$, and write
	\[
	\widehat V=\chi(V),\qquad
	\widehat H=\chi(H),\qquad
	\widehat F=\chi\circ F\circ\chi^{-1}.
	\]
	Choose a homeomorphism
	$a:[-1,1]^2\to\widehat V$ matching the four
	distinguished sides, with $u_2=-1$ and
	$u_2=1$ corresponding to the lower and
	upper horizontal sides, respectively.
	Set
	\[
	d(u)=\widehat F(a(u))-a(u).
	\]

By regularity, the two vertical sides of
$\widehat V$, which are disjoint from the
square's vertical boundary, are mapped
onto subarcs of the two opposite vertical
sides of the square. Hence, for some
$\epsilon\in\{-1,1\}$,
\[
\epsilon d_1(-1,u_2)<0<
\epsilon d_1(1,u_2).
\]
The horizontal sides of $\widehat V$
lie at heights $-1$ and $1$, whereas
$\widehat H$ avoids these heights. Thus
\[
d_2(u_1,-1)>0>d_2(u_1,1).
\]

	The Poincar\'e--Miranda theorem applied
	to $(\epsilon d_1,-d_2)$ gives
	$d(u_*)=0$ for some $u_*\in(-1,1)^2$.
	Consequently, $F$ has an interior
	fixed point.

	The strict boundary inequalities exclude
	boundary fixed points. Moreover, the
	boundary values of $d$ are homotopic
	through nonzero boundary maps to
	$u\mapsto(\epsilon u_1,-u_2)$, so
	\[
	\deg(d,(-1,1)^2,0)
	=-\epsilon\in\{-1,1\}.
	\]
	In particular, the degree need not be
	positive.

	For hyperbolicity, return to the original
	smooth coordinates and let $z$ be any
	fixed point. Put $B=DF_z$. Since $F$
	preserves a smooth area form at $z$
	and $F(z)=z$, we have $\det B=1$.

	The projective map $\mathbb P(B)$ sends
	$\overline I$ into $I$, while
	$\mathbb P(B^{-1})$ sends $\overline J$
	into $J$. The interval fixed-point
	theorem therefore gives an invariant
	line $E^u$ with projective direction
	in $I$ and an invariant line $E^s$
	with projective direction in $J$.
	Because $\overline I$ and $\overline J$
	are disjoint, these lines are distinct.

	Let $\lambda_u,\lambda_s$ be the
	corresponding real eigenvalues of $B$.
	In the eigenvector basis, use the
	projective coordinate $t$ in which
	$E^u$ corresponds to $0$ and $E^s$
	to infinity. Since $E^s\notin\overline I$,
	the first projective interval becomes
	a finite interval
	\[
	I=(t_-,t_+),\qquad t_-<0<t_+,
	\]
	on which the projective action is
	\[
	t\longmapsto
	\frac{\lambda_s}{\lambda_u}t.
	\]
	The strict inclusion
	\[
	\frac{\lambda_s}{\lambda_u}
	[t_-,t_+]\subset(t_-,t_+)
	\]
	implies
	\[
	\left|\frac{\lambda_s}{\lambda_u}\right|<1.
	\]
	Since $|\lambda_u\lambda_s|=1$, it follows
	that
	\[
	|\lambda_u|>1>|\lambda_s|.
	\]
	Thus every fixed point is hyperbolic.
	This proves the crossed-strip criterion
	used in \cite[Proposition~3.3]{AFLXZ}.
\end{proof}

\paragraph{The endpoint coordinate change.}
With $\Phi_n=\exp_{x_n}\circ i_n$, define
\[
I=\Phi_0^{-1}\circ\Phi_L
=i_0^{-1}\circ\exp_{x_0}^{-1}
\circ\exp_{x_L}\circ i_L.
\]
Then $\Phi_0(Iz)=\Phi_L(z)$; the map $I$
is a change of coordinates and involves
no iterate of $g$.

In a common smooth tangent trivialization,
the endpoint return condition
\eqref{H:eq:endpoint-return} and the
endpoint angle condition
\eqref{H:eq:endpoint-angle} give
\[
d(x_0,x_L)+\|i_L-i_0\|
\le C_Sq^{-\sigma},
\qquad
\|i_0^{-1}\|+\|i_L^{-1}\|
\le C_SA^3.
\]

Write
\[
T_{0L}=\exp_{x_0}^{-1}\circ\exp_{x_L},
\qquad
T_{L0}=\exp_{x_L}^{-1}\circ\exp_{x_0}
\]
in the common smooth tangent trivialization.
The transition maps depend smoothly on
their base points and reduce to the identity
when the base points coincide. Hence, on
uniform small exponential-coordinate
neighborhoods,
\[
\|T_{0L}-\mathrm{id}\|_{C^1}
+\|T_{L0}-\mathrm{id}\|_{C^1}
\le C_Sd(x_0,x_L).
\]
Moreover,
\[
I(z)-z
=i_0^{-1}\bigl[
(T_{0L}-\mathrm{id})(i_Lz)
+(i_L-i_0)z
\bigr].
\]

The chain rule, the endpoint frame bounds,
and the analogous formula for $I^{-1}$
therefore yield
\begin{equation}\label{H:eq:return-close}
	\|I-\mathrm{id}\|_{C^1(U_L)}
	+\|I^{-1}-\mathrm{id}\|_{C^1(I(U_L))}
	\le E_I:=C_{\rm ret}A^6q^{-\sigma}.
\end{equation}

The exponent $6$ is a convenient
non-optimal bound accounting for the
inverse-frame factors. The constant
$C_{\rm ret}$ depends only on the fixed
background geometry. In the application,
the final choice of $H_*$ ensures
\[
E_I\le\mathfrak D^{-20M/\beta},
\]
including the polynomial loss from
the boundary extension.

\Needspace{5\baselineskip}
The small endpoint coordinate change turns the long strip into a return branch in the initial box, with quantitative cone margins and expansion.

\begin{lemma}\label{H:lem:return-strip}
If $E_I\le\mathfrak D^{-20M/\beta}$,
the map $G=I\circ G'$ restricts to a
regular map from a vertical strip $R_1$
of $U_0$ onto a horizontal strip $R_2$
of $U_0$. Its internal vertical
boundary graphs have slope at most $\bar\kappa$, its internal horizontal
boundary graphs have slope at most $\bar\kappa/2$, and
$R_2\subset\{|w|<\bar r/3\}$.  It satisfies
\begin{equation}\label{H:eq:return-cones}
 DG_x\cC(\bar\kappa)\subset\cC(\bar\kappa/2)
       \quad(x\in R_1),
 \qquad
 D(G^{-1})_y\tcC(2\bar\kappa)\subset\tcC(\bar\kappa)
       \quad(y\in R_2).
\end{equation}
Moreover,
\begin{equation}\label{H:eq:central-return}
 0\in R_1,\qquad e:=G(0)\in R_2,\qquad |e|\le E_I.
\end{equation}
Put $\gamma=1-\eta-\rho$ and
\[
                         \Lambda=\tfrac12A^{\gamma L}.
\]
After increasing $A_*$ and $D_0$, we have $\Lambda\ge4$ and
\begin{align}
 |\pi_1DG_z\xi|&\ge\Lambda|\xi_1|
 &&(z\in R_1,\ \xi\in\cC(\bar\kappa)),
                                      \label{H:eq:return-horizontal-expansion}\\
 |\pi_2D(G^{-1})_z\zeta|&\ge\Lambda|\zeta_2|
 &&(z\in R_2,\ \zeta\in\tcC(2\bar\kappa)).
                                      \label{H:eq:return-vertical-expansion}
\end{align}
Every fixed point $z_*=(v_*,w_*)$ of $G$
satisfies
\begin{equation}\label{H:eq:fixed-point-localization}
	\|z_*\|_\infty\le\mu_L\bar r,
	\qquad
	\mu_L=C_{S,\beta}
	\left(
	K^pA^{-\gamma L}
	+K^{p-1}\mathfrak D^{-M}
	+\mathfrak D^{-15M/\beta}
	\right).
\end{equation}
In particular, $G$ has a hyperbolic fixed
point at distance $o(\bar r)$ from the origin,
uniformly as $A\to\infty$ over all
$\mathfrak D\ge A$ and $L\ge1$.
\end{lemma}

\begin{proof}
\textit{1. Returning the horizontal strip and pulling back its sides.}
By \eqref{H:eq:shape-bounds}--\eqref{H:eq:end-scales},
\[
 r_L=K^p\bar r,\qquad
 \tau_L+\kappa_Lr_L\le\bar r/10+K^{p-1}\bar\kappa\bar r<\bar r/5.
\]
The assumption on $E_I$ gives,
for sufficiently large $\mathfrak D$,
\[
 E_I<\bar r/10,\qquad K^p\bar r-E_I>2\bar r,\qquad
 E_I/\bar\kappa\longrightarrow0.
\]
The first two inequalities take care of positions; to keep track of
slopes we also need the $C^1$ part of \eqref{H:eq:return-close}.
If $B$ is a linear map with $\|B-\mathrm{Id}\|\le E_I$, a direction of
horizontal slope $h$ is sent to one of slope at most
\[
 \mathcal S(h)=\frac{h+E_I(1+h)}{1-E_I(1+h)},
\]
as long as the denominator is positive; the same formula applies to
vertical slopes.  It follows that
\[
 DI\cC(\bar\kappa/K)\subset\cC(\bar\kappa/2),\qquad
 DI^{-1}\tcC(2\bar\kappa)\subset\tcC(K\bar\kappa).
\]
Moreover, along a Lipschitz graph the
relevant horizontal or vertical projection
has derivative at least
$1-E_I(1+h)>0$ almost everywhere,
and hence is strictly monotone.

Let $G':R'\to R$ be the branch of \cref{H:lem:strip-concatenation}.
Each horizontal side of $I(R)$ is a
$\bar\kappa/2$-graph whose endpoints have first coordinates beyond
$\pm\bar r$, and
\[
 |\pi_2I(z)|\le\tau_L+\kappa_Lr_L+E_I<3\bar r/10\quad(z\in R).
\]
The two graphs are disjoint, and the
remaining two boundary arcs of $I(R)$
lie beyond the corresponding vertical
sides of $U_0$.
Since
\[
U_0=\{(v,w):|v|\le\bar r,\,
|w|\le\bar r+\bar\kappa|v|\}
\]
and $I(R)\subset\{|w|<3\bar r/10\}$,
the intersection with $U_0$ is exactly
the restriction to $|v|\le\bar r$.
Consequently,
\[
R_2=I(R)\cap U_0
\]
is the horizontal strip bounded by
these two restricted graphs, and
$R_2\subset\{|w|<\bar r/3\}$.
Let $\Gamma_\pm$ be its vertical sides.  By the inverse endpoint cone
inclusion, $I^{-1}\Gamma_\pm$ are vertical graphs joining the
horizontal sides of $R$.  Equation \eqref{H:eq:long-inverse-cone} and
regularity then pull them back to full $\bar\kappa$-vertical graphs
in $R'$, which miss its vertical boundary, since $I G'$ sends that
boundary outside $U_0$.  The strip between them is
\[
 R_1=(G'|_{R'})^{-1}(I^{-1}(R_2)),\qquad G=I G':R_1\to R_2,
\]
and the restriction is regular.  Composing the two endpoint cone
inclusions with
\eqref{H:eq:long-forward-cone}--\eqref{H:eq:long-inverse-cone}
gives \eqref{H:eq:return-cones}.

The induction in \cref{H:lem:strip-concatenation} gave $0\in R'$,
$0\in R$ and $G'(0)=0$, hence
\[
 G(0)=I(0)\in I(R),\qquad |I(0)|\le E_I<\bar r,
\]
so that $I(0)\in R_2$ and $0\in R_1$, which is
\eqref{H:eq:central-return}.

\textit{2. Coordinate expansion of the returned map.}
For $\xi\in\cC_n$ the one-step estimate reads
\[
 |\pi_1Dg_n\xi|
 \ge(A_n-\varepsilon_n(1+\kappa_n))|\xi_1|
 \ge e^{\lambda_n^u-\delta/4}|\xi_1|.
\]
In the other direction, given $\zeta\in\tcC_{n+1}$ put
$\xi=Dg_n^{-1}\zeta\in\tcC_n$; then
\[
 |\zeta_2|\le(B_n+\varepsilon_n(1+\widetilde\kappa_n))|\xi_2|
       \le B_ne^{\delta/4}|\xi_2|,
\]
because
\[
\varepsilon_n(1+\widetilde\kappa_n)
\le K\mathfrak D^{-3M}
(1+K\mathfrak D^{-M})
\le\mathfrak D^{-1}(e^{1/4}-1)
\]
for sufficiently large $\mathfrak D$.
Since $B_n=e^{\lambda_n^s}\ge A^{-1}
\ge\mathfrak D^{-1}$ and $\delta\ge1$,
this implies
\[
\varepsilon_n(1+\widetilde\kappa_n)
\le B_n(e^{\delta/4}-1).
\]
Consequently,
\[
|\pi_2Dg_n^{-1}\zeta|
\ge e^{-\lambda_n^s-\delta/4}
|\zeta_2|.
\]
Goodness at length $L$ together with
$\bar\lambda_n^e\le\min\{\lambda_n^u,-\lambda_n^s\}$ gives
\[
 \sum_{n<L}\lambda_n^u>(1-\eta)La,\qquad
 -\sum_{n<L}\lambda_n^s>(1-\eta)La.
\]
These one-step bounds may be iterated, since every intermediate point
and vector stays in the prescribed box and cone, and with
$\gamma=1-\eta-\rho>0$ we obtain
\begin{align*}
 |\pi_1DG'\xi|&\ge A^{(1-\eta-\rho/4)L}|\xi_1|
                   \ge A^{\gamma L}|\xi_1|,\\
 |\pi_2D(G'^{-1})\zeta|&\ge A^{\gamma L}|\zeta_2|.
\end{align*}
On the narrow cones in question, $DI$ and $DI^{-1}$ cost at most a
factor of $2$, because $1-E_I(1+2K\bar\kappa)\ge1/2$.  This proves
\eqref{H:eq:return-horizontal-expansion}--\eqref{H:eq:return-vertical-expansion}
with $\Lambda=\tfrac12A^{\gamma L}$, and choosing $A_*^\gamma\ge8$
makes $\Lambda\ge4$.

\textit{3. Localization of every fixed point of the branch.}
Let $G(z_*)=z_*$ with $z_*=(v_*,w_*)\in R_1$.  Then also
$z_*\in R_2\subset I(R)$, and since $\tau_L\le\bar r A^{-\gamma L}$,
\begin{equation}\label{H:eq:stable-location}
	\frac{|w_*|}{\bar r}
	\le\frac{\tau_L+\kappa_Lr_L+E_I}{\bar r}
	\le A^{-\gamma L}
	+K^{p-1}\mathfrak D^{-M}
	+\mathfrak D^{-15M/\beta}.
\end{equation}
For the first coordinate, write the vertical sides of $R'$ as
$v=\phi_\pm(w)$; they are defined on parameter intervals containing
$[-\bar r,\bar r]$ and have Lipschitz constant at most $\bar\kappa$.
The horizontal slice of $R'$ at $w=0$ is the interval
$[\phi_-(0),\phi_+(0)]$, which contains $0$.  On it we set
$f(v)=\pi_1G'(v,0)$.  Then
\[
 f(0)=0,\quad |f(\phi_\pm(0))|=r_L,\quad |f'(v)|\ge A^{\gamma L}.
\]
The derivative is continuous and never zero, so its sign is constant.
Integrating along each segment from $0$ to $\phi_\pm(0)$ gives
\[
 r_L=\left|\int_0^{\phi_\pm(0)}f'(v)\,dv\right|
     \ge A^{\gamma L}|\phi_\pm(0)|.
\]
Since $z_*\in R'$ with $|w_*|<\bar r/3$, where both boundary functions
are defined,
\[
 |v_*|\le\max_\pm|\phi_\pm(w_*)|
       \le r_LA^{-\gamma L}+\bar\kappa|w_*|.
\]
Together with $r_L=K^p\bar r$ and \eqref{H:eq:stable-location}, this
gives \eqref{H:eq:fixed-point-localization}.

Finally,
$G=\Phi_0^{-1}\widetilde g^L\Phi_0$
preserves the smooth area form
$\Phi_0^*\widehat\omega$
($\omega$ in the closed case).
We apply \cref{H:lem:crossed-strip}
with the cones
$\cC(\bar\kappa)$ and
$\tcC(2\bar\kappa)$.
Their product slope is
$2\bar\kappa^2<1$.
By continuity, the cone inclusions
\eqref{H:eq:return-cones} extend to
the closures of the input cones,
and the narrower target cones provide
strict projective containment.
Thus \cref{H:lem:crossed-strip} gives
existence and hyperbolicity of a fixed point.

For fixed $\beta$, since
$L\ge1$ and $\mathfrak D\ge A$,
the quantity $\mu_L$ in
\eqref{H:eq:fixed-point-localization}
tends to zero as $A\to\infty$,
uniformly over all
$\mathfrak D\ge A$ and $L\ge1$.
Hence the localization is $o(\bar r)$
with the asserted uniformity.
\end{proof}

\Needspace{5\baselineskip}
For the boundary argument we also need stable and unstable graphs crossing the whole rectangle through the returned fixed point.

\begin{lemma}\label{H:lem:full-cross}
	Let $U=[-r,r]\times[-t,t]$, where $r,t>0$.
	Fix positive slopes $h,v$ with $hv<1$,
	and define the closed cones
	\[
	\mathcal C_h=\{\xi:|\xi_2|\le h|\xi_1|\},
	\qquad
	\mathcal C_v=\{\zeta:|\zeta_1|\le v|\zeta_2|\}.
	\]
	Let $V\subset U$ be a full vertical strip
	whose internal vertical sides are
	$v$-Lipschitz graphs, and let $H\subset U$
	be a full horizontal strip whose internal
	horizontal sides are $h$-Lipschitz graphs.
	Suppose that a $C^1$ diffeomorphism $F$
	between neighborhoods of $V$ and $H$
	maps $V$ regularly onto $H$.

	Assume, for nonzero vectors, the strict
	cone inclusions
	\begin{equation}\label{H:eq:full-cross-cones}
		DF_z\mathcal C_h\Subset\mathcal C_h
		\quad(z\in V),
		\qquad
		D(F^{-1})_z\mathcal C_v\Subset\mathcal C_v
		\quad(z\in H),
	\end{equation}
	where compact containment is understood
	in projective directions.
	Assume also that there is $\Lambda>1$
	such that
	\begin{align*}
		|\pi_1DF_z\xi|
		&\ge\Lambda|\xi_1|
		&&(z\in V,\ \xi\in\mathcal C_h),\\
		|\pi_2D(F^{-1})_z\zeta|
		&\ge\Lambda|\zeta_2|
		&&(z\in H,\ \zeta\in\mathcal C_v).
	\end{align*}

	Then the branch $F:V\to H$ has a unique
	hyperbolic fixed point $p$.
	Through $p$ pass a full $C^1$ vertical
	graph $W_U^s(p)\subset V$ and a full
	$C^1$ horizontal graph $W_U^u(p)\subset H$,
	with slopes bounded by the corresponding
	cone slopes.

	For every $n\ge0$, they satisfy
	\[
	|F^n(x)-p|
	\le 2t\sqrt{1+v^2}\,\Lambda^{-n}
	\qquad (x\in W_U^s(p)),
	\]
	and
	\[
	|F^{-n}(x)-p|
	\le 2r\sqrt{1+h^2}\,\Lambda^{-n}
	\qquad (x\in W_U^u(p)).
	\]
	In particular, the two graphs meet
	exactly at $p$, and their union is
	a connected full cross in $U$:
	the stable graph joins the horizontal
	sides, while the unstable graph joins
	the vertical sides.
\end{lemma}

\begin{proof}
\textit{1. Finite strips and their widths.}
For $n\ge1$ set
\[
 V_n=\{z\in V:F^j(z)\in V\ (0\le j<n)\},\qquad H_n=F^n(V_n).
\]
The strip-intersection argument of \cref{H:lem:strip-concatenation}
applies to these sets with the same map at every step; explicitly,
\[
 V_1=V,\quad H_1=H,\qquad
 V_{n+1}=(F^n|_{V_n})^{-1}(H_n\cap V),\quad
 H_{n+1}=F(H_n\cap V).
\]
Since $H_n$ is a full horizontal strip
and $V$ is a full vertical strip,
the full-graph boundary conditions
give the four corner intersections.
Each is unique because $hv<1$.
Regularity and the inverse cones pull the vertical sides back to full
vertical graphs, and the forward cones carry the horizontal sides to
full horizontal graphs.  This makes each $V_{n+1}\subset V_n$ a full
vertical strip and $F^n:V_n\to H_n$ regular.  Note also that every
derivative in an $n$-fold composition is evaluated at a point of its
required domain.

Denote the two boundary functions by
$\varphi_n^-\le\varphi_n^+$ on $[-t,t]$, and for fixed $y$ put
$J_n(y)=[\varphi_n^-(y),\varphi_n^+(y)]$ and
$f_{n,y}(x)=\pi_1F^n(x,y)$.  Cone invariance and the expansion give
\[
 |f_{n,y}'(x)|\ge\Lambda^n\quad(x\in J_n(y)).
\]
This derivative is continuous and never vanishes, hence of constant
sign, and since $F^n(V_n)\subset U$,
\begin{equation}\label{eq:cross-strip-width}
 \Lambda^n(\varphi_n^+(y)-\varphi_n^-(y))
 \le\int_{J_n(y)}|f_{n,y}'(x)|\,dx
 =|f_{n,y}(\varphi_n^+(y))-f_{n,y}(\varphi_n^-(y))|
 \le2r.
\end{equation}
The left and right boundary functions increase and decrease,
respectively; they are $v$-Lipschitz, and their gap is uniformly at
most $2r\Lambda^{-n}$.  Their common limit $\varphi$ is again
$v$-Lipschitz, and
\[
 W_U^s:=\bigcap_{n\ge1}V_n
       =\{(\varphi(y),y):-t\le y\le t\}.
\]
For the inverse branch we set
$\widehat H_n=\{z\in H:F^{-j}(z)\in H\ (0\le j<n)\}$; the same
argument with the coordinates interchanged gives
\[
 W_U^u:=\bigcap_{n\ge1}\widehat H_n
       =\{(x,\psi(x)):-r\le x\le r\},\qquad \Lip(\psi)\le h,
\]
and the bound $2t\Lambda^{-n}$ on the vertical widths of
$\widehat H_n$.

\textit{2. Identification of the fixed point and the invariant sets.}
Since $y\mapsto\psi(\varphi(y))$ maps $[-t,t]$ into itself with
Lipschitz constant at most $hv<1$, the two limiting graphs meet in a
single point.  Now $F$ has a fixed point by
\cref{H:lem:crossed-strip}, and every fixed point belongs to all $V_n$
and all $\widehat H_n$; so the intersection point is the unique fixed
point $p$.  The orbit conditions defining the strips give the
stronger identities
\begin{equation}\label{eq:cross-invariance-exact}
 F(W_U^s)=W_U^s\cap H,\qquad
 F^{-1}(W_U^u)=W_U^u\cap V.
\end{equation}
For instance, a point $z\in W_U^s\cap H$ has $F^{-1}z\in V$, with
all subsequent iterates in $V$, so $F^{-1}z\in W_U^s$.

\textit{3. Contraction along connected subarcs.}
Write the horizontal boundaries of $H$ as $y=\psi_\pm(x)$.  For
$y_2>y_1$ we have
\begin{align*}
 &[y_2-\psi_\pm(\varphi(y_2))]
       -[y_1-\psi_\pm(\varphi(y_1))]\\
 &\hspace{30mm}\ge(1-hv)(y_2-y_1)>0.
\end{align*}
It follows that $W_U^s\cap H$ is a connected vertical subarc.  Given
$x\in W_U^s$ and $k\ge1$, the subarc joining $F^k(x)$ to $p$ lies in
$H$, and by \eqref{eq:cross-invariance-exact} its preimage lies in
$W_U^s$.  Parameterizing it by $\zeta(y)=(\varphi(y),y)$, we have
almost everywhere
\[
 \zeta'(y)\in \mathcal C_v,\qquad
 \left|\frac{d}{dy}\pi_2F^{-1}(\zeta(y))\right|\ge\Lambda.
\]
The scalar function in the second
relation is absolutely continuous.
By \eqref{eq:cross-invariance-exact},
the image under $F^{-1}$ of this
subarc is a connected subarc of
the vertical graph $W_U^s$.
Consequently, its second-coordinate
projection is injective and hence
strictly monotone. Integrating
without cancellation gives
\[
|\pi_2(F^k(x)-p)|
\le\Lambda^{-1}
|\pi_2(F^{k-1}(x)-p)|.
\]
Consequently
\begin{align}
 |F^n(x)-p|&\le2t\sqrt{1+v^2}\,\Lambda^{-n}
          &&(x\in W_U^s),\label{eq:cross-stable-rate}\\
 |F^{-n}(x)-p|&\le2r\sqrt{1+h^2}\,\Lambda^{-n}
          &&(x\in W_U^u).\label{eq:cross-unstable-rate}
\end{align}
The second estimate uses the connected horizontal subarc $W_U^u\cap V$
and forward expansion in the same way.  This confirms that the
limiting graphs, defined as geometric limits of strips, are the
stable and unstable sets of the branch.

\textit{4. Hyperbolicity and $C^1$
	regularity up to the sides.}
The strict cone conditions at $p$
produce two transverse invariant lines
$E^u\subset\mathcal C_h$ and
$E^s\subset\mathcal C_v$.
Writing $\lambda_u,\lambda_s$ for
the corresponding eigenvalues of
$DF_p$, the coordinate expansions
applied to their nonzero eigenvectors give
\[
|\lambda_u|\ge\Lambda>1,\qquad
|\lambda_s|^{-1}\ge\Lambda>1.
\]
Hence $p$ is hyperbolic; area
preservation plays no role here. The
$C^1$ stable-manifold theorem gives local $C^1$ graphs at $p$.
Choose a local stable neighborhood $N$ for the usual orbit
characterization of the local stable manifold; by the uniform
estimate \eqref{eq:cross-stable-rate} there is $N_0$ with
\[
 F^{N_0+j}(W_U^s)\subset N\quad(j\ge0),
 \qquad F^{N_0}(W_U^s)\subset W^s_{\rm loc}(p).
\]
Near any $x\in W_U^s$, the image of a small relative graph interval
is an arc of $W^s_{\rm loc}(p)$, and pulling it back by the local
$C^1$ inverse of $F^{N_0}$ gives $C^1$ regularity near $x$.  At the
endpoints on $\partial U$ the same argument applies one-sidedly.  The
tangent lies in $\mathcal C_v$, either by the Lipschitz bound or by
inverse cone invariance.  Backward convergence gives the corresponding
assertions for $W_U^u$, and the full-cross conclusion follows.
\end{proof}

\Needspace{6\baselineskip}
We apply the full-cross lemma to a rectangular restriction of the
returned moving-box branch.  The geometric fact below forces a
sufficiently flat separating curve close to the crossing point to
meet one of its stable or unstable branches.

\begin{lemma}\label{H:lem:flat-seam}
	There are universal constants
	$\epsilon_0,\kappa_{\rm fs},\theta_0>0$
	with the following property.

	Let $R>0$, let $p\in\mathbb R^2$, and put
	\[
	Q_R=p+[-R,R]^2,\qquad
	2Q_R=p+[-2R,2R]^2.
	\]
	Suppose that $W^v,W^h\subset Q_R$ are,
	respectively, a full vertical and a full
	horizontal $C^1$ graph through $p$,
	both of slope at most $\kappa_{\rm fs}$.

	Let $N$ be an open neighborhood of $2Q_R$,
	and let $\Gamma\subset N$ be a relatively
	closed, embedded $C^1$ one-dimensional
	submanifold without boundary.
	Assume that
	\[
	N\setminus\Gamma=N_+\sqcup N_-,
	\]
	where $N_+$ and $N_-$ are disjoint open
	subsets of $N$.
	Suppose that $p\in N_+$ and that
	there exists $z\in N_-$ satisfying
	\[
	|z-p|\le\epsilon_0R.
	\]
	If
	\[
	\operatorname{diam}_{\mathbb{RP}^1}
	\{T_x\Gamma:x\in\Gamma\cap Q_R\}
	\le\theta_0,
	\]
	then
	\[
	\Gamma\cap(W^v\cup W^h)\ne\varnothing.
	\]
\end{lemma}

\begin{proof}
	Fix universal constants with
	\begin{equation}\label{eq:flat-constants}
		0<\kappa_{\rm fs}<\tfrac18,\qquad
		\theta_0+4\kappa_{\rm fs}<\tfrac\pi8,\qquad
		c_*:=\cos(3\pi/8)>0,\qquad
		(1+\kappa_{\rm fs})(1+c_*^{-1})
		\epsilon_0<\tfrac18.
	\end{equation}
	Translate $p$ to $0$ and suppose, for a contradiction, that
	$\Gamma\cap(W^v\cup W^h)=\varnothing$.  Write
	\[
	W^v=\{(\phi(y),y):|y|\le R\},\qquad
	W^h=\{(x,\psi(x)):|x|\le R\},
	\]
	with $\phi(0)=\psi(0)=0$ and
	$\|\phi'\|_0,\|\psi'\|_0\le\kappa_{\rm fs}$.

	\textit{1. Straightening and a quantitative common domain.}
	On $Q_R$ define $\Phi(x,y)=(x-\phi(y),y-\psi(x))$.  For
	$u,u'\in Q_R$ we have
	\begin{equation}\label{eq:flat-bilip}
		(1-\kappa_{\rm fs})|u-u'|
		\le|\Phi(u)-\Phi(u')|
		\le(1+\kappa_{\rm fs})|u-u'|,
		\qquad
		\|D\Phi-\mathrm{Id}\|\le\kappa_{\rm fs}.
	\end{equation}
	Together with $\det D\Phi=1-\phi'(y)\psi'(x)\ge1-\kappa_{\rm fs}^2>0$ this
	makes $\Phi$ a $C^1$ diffeomorphism onto its image, understood
	one-sidedly at the boundary of the square.  For
	$|(X,Y)|_\infty\le R/2$ the map
	\[
	(x,y)\longmapsto(X+\phi(y),Y+\psi(x))
	\]
	is a $\kappa_{\rm fs}$-contraction of $Q_R$ into its interior, and its fixed
	point is the unique preimage of $(X,Y)$.  Since also
	$\|\Phi(u)\|_\infty\le(1+\kappa_{\rm fs})R$ on $Q_R$, we obtain
	\begin{equation}\label{eq:flat-square-domains}
		[-R/2,R/2]^2\Subset\Phi(\operatorname{int}Q_R),\qquad
		\Phi(Q_R)\subset[-2R,2R]^2.
	\end{equation}
	Here compact containment is meant in the ordinary Euclidean sense.

	The vertical graph is sent to
	$\{(0,y-\psi(\phi(y))):|y|\le R\}$, and for $y_2>y_1$
	\[
	[y_2-\psi(\phi(y_2))]-[y_1-\psi(\phi(y_1))]
	\ge(1-\kappa_{\rm fs}^2)(y_2-y_1).
	\]
	The endpoint values therefore straddle
	$[-(1-\kappa_{\rm fs}^2)R,(1-\kappa_{\rm fs}^2)R]$, and with the analogous
	horizontal estimate we get
	\begin{equation}\label{eq:flat-axes-covered}
		\bigl(\{0\}\times[-R/2,R/2]\bigr)
		\cup\bigl([-R/2,R/2]\times\{0\}\bigr)
		\subset\Phi(W^v\cup W^h).
	\end{equation}
	It is here that the full lengths of the two graphs are used.

	\textit{2. A nearby boundary point and its tangent direction.}
	For a linear map $B$ with
	$\|B-\mathrm{Id}\|\le\kappa_{\rm fs}<1/8$
	we have
	\[
	d_{\mathbb{RP}^1}(B\ell,\ell)
	\le\arcsin\frac{\kappa_{\rm fs}}{1-\kappa_{\rm fs}}\le2\kappa_{\rm fs}.
	\]
	Put $Q'=[-R/2,R/2]^2$ and consider
	$\widehat\Gamma=\Phi(\Gamma\cap\operatorname{int}Q_R)$ on a
	neighborhood of $Q'$.  Since every preimage of a point of $Q'$ lies
	in $Q_R$, the tangent hypothesis and the triangle inequality give
	\begin{equation}\label{eq:flat-tangent-bound}
		\operatorname{diam}_{\mathbb{RP}^1}
		\{T_u\widehat\Gamma:u\in\widehat\Gamma\cap Q'\}
		\le\theta_0+4\kappa_{\rm fs}<\pi/8.
	\end{equation}
	The segment $[0,\Phi(z)]$ stays inside $Q'$ by
	\eqref{eq:flat-constants}--\eqref{eq:flat-bilip}, and it joins opposite
	sides of $\widehat\Gamma$, so it contains a point $y\in\widehat\Gamma$
	with
	\begin{equation}\label{eq:flat-nearby-point}
		|y|\le|\Phi(z)|\le(1+\kappa_{\rm fs})\epsilon_0R.
	\end{equation}
	If $y$ lies on one of the axes, \eqref{eq:flat-axes-covered} already
	contradicts the disjointness assumption.

	\textit{3. Hitting an axis before leaving the domain.}
	Otherwise, let $\nu$ be a unit tangent at $y$, choose $i\in\{1,2\}$
	with $|\nu_i|\ge1/\sqrt2$, and orient it so that
	$\operatorname{sign}(y_i)\nu_i<0$.  Let $\gamma(s)$ be the
	corresponding unit-speed parameterization, with $\gamma(0)=y$.  The
	continuous lift of its tangent to the unit circle stays in the
	component containing $\nu$, since the other lift is separated from it
	by nearly $\pi$.  So as long as $\gamma$ stays in $Q'$,
	\eqref{eq:flat-tangent-bound} yields
	\begin{equation}\label{eq:flat-directed-speed}
		\operatorname{sign}(y_i)\gamma_i'(s)\le-c_*.
	\end{equation}
	Set $s_0:=|y_i|/c_*$.
	For every $s$ in the interval of existence
	of $\gamma$ with $0\le s\le s_0$,
	unit speed and \eqref{eq:flat-nearby-point}
	give
	\[
	|\gamma(s)|\le|y|+s
	\le(1+c_*^{-1})(1+\kappa_{\rm fs})
	\epsilon_0R<R/8.
	\]
	Thus the curve stays in a compact subset
	of $\operatorname{int}Q'$.
	Since $\widehat\Gamma$ is relatively
	closed in $\Phi(\operatorname{int}Q_R)$
	and has no boundary there, its
	unit-speed parametrization cannot
	terminate at a finite time while
	remaining in this compact subset.
	Consequently, $\gamma$ is defined
	through time $s_0$.

	Integrating \eqref{eq:flat-directed-speed},
	we obtain
	\[
	\operatorname{sign}(y_i)\gamma_i(s_0)
	\le |y_i|-c_*s_0=0.
	\]
	By continuity there exists $s_*\le s_0$
	such that $\gamma_i(s_*)=0$.
	This point lies in $Q'$ and hence,
	by \eqref{eq:flat-axes-covered},
	belongs to $\Phi(W^v\cup W^h)$.
	Injectivity of $\Phi$ gives the
	required contradiction.
\end{proof}

\begin{proof}[Proof of \cref{thm:finite-closing}]
	\textit{1. Constants and domains.}
	The background geometry and $\beta$ are fixed throughout.  Take
	$p,\eta,\rho$ as in \eqref{H:eq:small-constants} and
	$\vartheta_*=\vartheta_+(\eta)$ from \cref{H:lem:good-segment}.
	Choose $A_*$ above the finitely many thresholds of the extension,
	coordinate, scale and strip lemmas, and so large that
	\[
	A_*\ge\max\{e^{1/\rho},K^{2p},(10K)^{1/(M-2)}\}.
	\]
	For $\mathfrak D\ge A\ge A_*$ we then have
	\[
	\delta\ge1,\qquad A_*^{1/(2p)}\ge K,\qquad
	\mathfrak D^M/K\ge A^M/K\ge10A^2.
	\]
	The smallness requirements below for the universal constants of
	\cref{H:lem:flat-seam} are imposed now as well.
	Finally we fix $H_*\ge H_0$ so large that, with $a_*=20M/\beta$,
	\[
	\sigma H_*\ge6+a_*m_{\rm ext}+1+
	\frac{\log^+(C_{\rm ret}C_{\rm ext}^{a_*})}{\log A_*},
	\qquad \log^+x=\max\{0,\log x\}.
	\]
	Then $q\ge D^{H_*}$ and $A\le D$ give $q\ge A^{H_0}$ and
	\begin{equation}\label{H:eq:extended-return-small}
		C_{\rm ret}A^6q^{-\sigma}
		\le C_{\rm ret}D^{6-\sigma H_*}
		\le C_{\rm ext}^{-a_*}D^{-a_*m_{\rm ext}}
		\le\mathfrak D^{-a_*}
		\quad\bigl(\mathfrak D\le C_{\rm ext}D^{m_{\rm ext}}\bigr).
	\end{equation}
	This choice does not depend on $g,A,D,q$.

	Now apply \cref{H:lem:good-segment} to $g$.  Because of its strict
	exponent reserve, the initial point, and with it the whole rebased orbit
	$x_0,\ldots,x_L$, can be taken in $\operatorname{int}S$.
	When $S$ is closed we put $\widetilde g=g$ and $\mathfrak D=D$;
	otherwise we invoke \cref{H:lem:symplectic-extension} and put
	$\mathfrak D=\Delta=C_{\rm ext}D^{m_{\rm ext}}$.
	We then run the local construction with
	$\mathfrak D$ as its norm parameter:
	\[
	\bar r=\mathfrak D^{-5M/\beta},
	\qquad
	\bar\kappa=\mathfrak D^{-M},
	\qquad
	\bar b=\mathfrak D^M.
	\]
	Since $\widetilde g=g$ along the selected
	orbit segment,
	\[
	D\widetilde g(x_n)=Dg(x_n)
	\qquad(0\le n<L).
	\]
	The rates therefore still satisfy
	$|\lambda_n^{u,s}|\le\log A$.
	The physical source boxes have radius at most $C_{S,\beta}K^p\bar r$
	and their images radius at most $C_{S,\beta}K^p\mathfrak D\bar r$, so
	in the boundary case the collar margin they require follows from
	\begin{equation}\label{eq:closing-collar-margin}
		\frac{C_{S,\beta}K^p\mathfrak D\bar r}{\mathfrak D^{-1}}
		=C_{S,\beta}K^p\mathfrak D^{2-5M/\beta}
		\longrightarrow0.
	\end{equation}
	We choose the thresholds so that this ratio is below $1/100$.
	Every source and image set then lies in the extension and chart
	domains, whether or not it comes near the boundary, and
	\eqref{H:eq:extended-return-small} supplies the hypothesis of
	\cref{H:lem:return-strip}.

	\textit{2. A rectangular returned branch.}
	Let $z_*$ be a fixed point of $G:R_1\to R_2$ and
	$P=\Phi_0(z_*)$ the corresponding point of the ambient surface.
	Since $G=\Phi_0^{-1}\widetilde g^L\Phi_0$ is hyperbolic at $z_*$, the
	return $\widetilde g^L(P)=P$ is hyperbolic.  If $S$ is closed, or if
	$P\in S$, every intermediate iterate is an iterate of $g$ inside $S$,
	and we have found the desired periodic point.  Suppose then, for a
	contradiction, that $P\in\widehat S\setminus S$.

	Set
	\[
	B=[-\bar r,\bar r]^2\subset U_0,\qquad
	V=R_1\cap B,\qquad H=G(V).
	\]
	The internal vertical sides of $R_1$ are graphs
	$v=\varphi_\pm(w)$ with $\Lip(\varphi_\pm)\le\bar\kappa$, and their
	endpoints on $\partial_hU_0$ sit at heights of absolute value at least
	$\bar r$, so that
	\[
	V=\{(v,w):|w|\le\bar r,\ \varphi_-(w)\le v\le\varphi_+(w)\}
	\]
	is a full vertical strip of $B$.
	Along the horizontal cuts
	$\ell_\pm(v)=(v,\pm\bar r)$,
	where
	$\varphi_-(\pm\bar r)\le v
	\le\varphi_+(\pm\bar r)$,
	we have
	\[
	|\partial_v\pi_1G(\ell_\pm(v))|\ge\Lambda>1,\qquad
	|\partial_v\pi_2G(\ell_\pm(v))|
	\le(\bar\kappa/2)|\partial_v\pi_1G(\ell_\pm(v))|.
	\]
	The first of these derivatives is continuous and never vanishes, hence
	of constant sign.  Regularity sends the two endpoints to opposite
	vertical sides of $R_2$, at $v=\pm\bar r$, so each image is a full
	$\bar\kappa/2$-horizontal graph.  As
	$H\subset R_2\subset\{|w|<\bar r/3\}$, these graphs stay strictly
	between the horizontal sides of $B$.  By injectivity and planar
	separation, these two images and the two subarcs of the vertical
	sides bound exactly the region $H$; hence $G:V\to H$ is regular.

	The cone estimates of \cref{H:lem:return-strip} are stated for open
	cones, but they pass to the closures by continuity.  The factor of two
	in the slopes gives strict projective containment for the closed cones
	with
	\[
	h=\bar\kappa,\qquad v=2\bar\kappa,\qquad hv=2\bar\kappa^2<1.
	\]
	The coordinate expansion estimates likewise extend to the boundary
	directions, and $\Lambda\ge4$.  For large $A_*$ we have
	$\|z_*\|_\infty<\bar r$, so $z_*\in V$, and \cref{H:lem:full-cross}
	produces full $C^1$ stable and unstable graphs $W^s,W^u$ across $B$
	through the same point $z_*$.

	\textit{3. Localization and a full cross
		at the boundary scale.}
	By \eqref{H:eq:fixed-point-localization},
	\[
	\|z_*\|_\infty\le\mu_L\bar r,
	\qquad
	|z_*|\le\sqrt2\,\mu_L\bar r,
	\]
	where
	\[
	\mu_L=C_{S,\beta}
	\bigl(
	K^pA^{-\gamma L}
	+K^{p-1}\mathfrak D^{-M}
	+\mathfrak D^{-15M/\beta}
	\bigr).
	\]
	Since $L\ge1$ and $\mathfrak D\ge A$,
	\[
	\mu_L\le C_{S,\beta}
	\bigl(
	K^pA^{-\gamma}
	+K^{p-1}A^{-M}
	+A^{-15M/\beta}
	\bigr)\longrightarrow0
	\]
	uniformly in $L$.
	We arrange that
	\[
	\mu_L<
	\min\left\{\frac16,
	\frac{\epsilon_0}{3\sqrt2}\right\},
	\qquad
	2\bar\kappa<
	\min\{1,\kappa_{\rm fs}\}.
	\]
	With $R=\bar r/3$ and $Q_R=z_*+[-R,R]^2$ this gives
	\[
	\|z_*\|_\infty+2R\le(\mu_L+2/3)\bar r<5\bar r/6,
	\qquad 2Q_R\Subset B.
	\]
	Write $W^s=\{(\varphi(w),w)\}$, where
	$\varphi(w_*)=v_*$ and $\Lip(\varphi)\le2\bar\kappa$.  For
	$|w-w_*|\le R$ we get
	\[
	|\varphi(w)-v_*|\le2\bar\kappa|w-w_*|<R.
	\]
	Hence $W^s\cap Q_R$ joins the two horizontal sides and misses the
	vertical ones, and the same estimate with slope $\bar\kappa$ handles
	$W^u$.  The two restrictions therefore form a full cross in $Q_R$.

	\textit{4. Quantitative flatness of the physical boundary.}
	Recall that
	$\Phi_0=\exp_{x_0}\circ i_0$.
	Choose an open coordinate neighborhood
	$N$ containing $2Q_R$ such that
	$\Phi_0(N)$ lies in a single smooth
	boundary chart of $\widehat S$.
	Set
	\[
	\Gamma=\Phi_0^{-1}(\partial S)\cap N.
	\]
	Since $0$ represents
	$x_0\in\operatorname{int}S$ and
	$z_*$ represents $P\in\widehat S\setminus S$,
	the two points lie on opposite sides
	of $\Gamma$.

	The trace $\Gamma$ is relatively closed
	in $N$ and is an embedded $C^1$
	one-manifold without boundary.
	The signed boundary defining function
	separates $N\setminus\Gamma$ into
	the interior and exterior open sides.
	The boundary tangent field is Lipschitz
	in this fixed chart.
	Connectedness of $\Gamma\cap Q_R$
	is not required.
	For physical boundary points $y,y'$ in this neighborhood, the uniform
	$C^2$ bounds on the exponential charts give, in $T_{x_0}\widehat S$,
	\[
	d_{\mathbb{RP}^1}\bigl(D\exp_{x_0}^{-1}|_y(T_y\partial S),
	D\exp_{x_0}^{-1}|_{y'}(T_{y'}\partial S)\bigr)
	\le C_Sd(y,y').
	\]
	For any invertible linear map $B_0$, the sine formula for two unit
	vectors gives
	\[
	\Lip(\mathbb P B_0)\le C\cond(B_0)^2.
	\]
	Here $\cond(i_0)\le C_SA^3$ by \eqref{H:eq:endpoint-angle}, while
	$\|i_0\|\le\sqrt2$ gives
	$\operatorname{diam}\Phi_0(Q_R)\le C_SR$.  Combining these estimates
	we obtain
	\begin{equation}\label{H:eq:boundary-tangent-flatness}
		\begin{split}
			\operatorname{diam}_{\mathbb{RP}^1}
			\{T_z\Gamma:z\in\Gamma\cap Q_R\}
			&\le C_S\cond(i_0)^2R\le C_SA^6\bar r\\
			&\le C_{S,\beta}A^6\bar r^\beta
			=C_{S,\beta}A^6\mathfrak D^{-5M}
			\le C_{S,\beta}\mathfrak D^{6-5M}.
		\end{split}
	\end{equation}
	The right-hand side tends to zero, so after enlarging the initial
	threshold it is at most $\theta_0$.
	Since also
	\[
	|z_*|\le3\sqrt2\,\mu_LR
	\le\epsilon_0R,
	\]
	every hypothesis of
	\cref{H:lem:flat-seam} is satisfied
	with center $p=z_*$ and
	opposite-side point $z=0$, and we conclude that
	$\Gamma\cap(W^s\cup W^u)\ne\varnothing$.

	\textit{5. Boundary invariance excludes an exterior saddle.}
	Suppose first that $y\in\Gamma\cap W^s$.  Then every forward iterate
	under the branch is defined and stays in the branch, and
	$G^n(y)\to z_*$.  Invariance of the boundary and the coordinate
	conjugacy give
	\[
	\Phi_0(G^n(y))=\widetilde g^{nL}(\Phi_0(y))\in\partial S,
	\qquad \Phi_0(G^n(y))\longrightarrow P.
	\]
	Since $\partial S$ is closed in $\widehat S$, this forces
	$P\in\partial S$.  If instead $y\in\Gamma\cap W^u$, the same argument
	applies with $G^{-n}$.  In either case we contradict $P\notin S$.
	The limit here is taken in the physical set $\partial S$, which is
	closed, and only relative closedness of its coordinate trace is used.

	We have shown that $P\in S$ and $g^L(P)=P$.  If $m$ is its least
	positive period, then $m\mid L$ and
	\[
	Dg^L(P)=\bigl(Dg^m(P)\bigr)^{L/m}.
	\]
	An eigenvalue of modulus one for $Dg^m(P)$ would therefore produce one
	for $Dg^L(P)$, which is impossible since the return is hyperbolic.
	Hence the derivative at the least period is hyperbolic as well.  If
	$g$ has no periodic point on $\partial S$, then necessarily
	$P\in\operatorname{int}S$.
\end{proof}

\begin{figure}[H]
	\centering
	\begin{tikzpicture}[x=1.45cm,y=1.45cm,>=Latex,font=\small]
		\begin{scope}
			\filldraw[fill=gray!4,draw=black!65,line width=.6pt]
			(-1.3,-1.508)--(0,-1.3)--(1.3,-1.508)--(1.3,1.508)
			--(0,1.3)--(-1.3,1.508)--cycle;
			\filldraw[fill=purple!13,draw=purple!65!black,line width=.8pt]
			(-.23,-1.3368)--(.25,-1.34)--(.25,1.34)--(0,1.3)
			--(-.23,1.3368)--cycle;
			\filldraw[fill=cyan!14,draw=cyan!55!black,line width=.8pt]
			(-1.3,-.20)--(1.3,-.14)--(1.3,.16)--(-1.3,.12)--cycle;
			\node[anchor=south west] at (-1.29,1.53) {$U_0$};
			\node[purple!65!black] at (.02,.97) {$R_1$};
			\node[cyan!45!black] at (1.00,-.01) {$R_2$};
			\draw[->,line width=.8pt] (.27,.90) .. controls (.86,.90) and (1.13,.55) .. (1.14,.23)
			node[pos=.45,above] {$G$};
			\fill (0,0) circle (1.3pt);
			\fill[blue!70!black] (.10,.045) circle (1.3pt);
			\fill (-.085,-.08) circle (1.3pt);
			\node[anchor=east] at (-.37,.38) {$0$};
			\draw[thin,gray] (-.33,.31)--(-.015,.015);
			\node[blue!70!black,anchor=west] at (.38,.34) {$e$};
			\draw[thin,blue!45] (.38,.27)--(.12,.06);
			\node[anchor=east] at (-.37,-.48) {$z_*$};
			\draw[thin,gray] (-.34,-.40)--(-.10,-.10);
			\node[font=\footnotesize] at (0,-1.78) {(a) moving box and returned strips};
		\end{scope}
		\begin{scope}[xshift=6.8cm]
			\fill[orange!10] (-1.3,-1.3) rectangle (1.3,1.3);
			\fill[blue!8] (-1.3,.12)--(1.3,.38)--(1.3,1.3)--(-1.3,1.3)--cycle;
			\draw[black,line width=.8pt] (-1.3,.12)--(1.3,.38);
			\draw[green!50!black,line width=1pt]
			(-.04,-1.3) .. controls (.07,-.45) and (-.07,.45) .. (.04,1.3);
			\draw[red!65!black,line width=1pt,densely dashed]
			(-1.3,-.03) .. controls (-.45,.05) and (.45,-.05) .. (1.3,.03);
			\draw[black!65,line width=.6pt] (-1.3,-1.3) rectangle (1.3,1.3);
			\fill (0,0) circle (1.5pt);
			\fill[blue!70!black] (.26,.43) circle (1.4pt);
			\fill (-.015,.249) circle (1.4pt);
			\node[anchor=south east] at (1.3,1.34) {$Q_R$};
			\node[font=\footnotesize,blue!50!black,align=center] at (.76,.98) {physical\\[-1pt]side};
			\node[font=\footnotesize,orange!70!black] at (.73,-.96) {added side};
			\node[anchor=east] at (-.18,-.22) {$z_*$};
			\node[blue!70!black,anchor=west] at (.32,.53) {$0$};
			\node[anchor=east] at (-.20,.56) {$y$};
			\draw[thin,gray] (-.18,.49)--(-.025,.27);
			\node at (1.0,.55) {$\Gamma$};
			\node[green!50!black] at (-.39,-.88) {$W^s$};
			\node[red!65!black] at (.96,-.22) {$W^u$};
			\node[font=\footnotesize] at (0,-1.78) {(b) exclusion of an exterior saddle};
		\end{scope}
	\end{tikzpicture}
	\caption{Returned strips and the boundary argument.  Panel (a) uses the
		actual V-shaped upper and reflected lower sides of $U_0$; relative
		widths and displacements are enlarged for visibility.  The map sends
		$R_1$ regularly onto $R_2$ and $e=G(0)$ is close to the origin.
		In (b), an exterior saddle would be separated from the physical point
		$0$ by $\Gamma$.  The full stable--unstable cross forces an intersection
		$y$, and invariance of the physical boundary together with convergence to $z_*$ gives
		a contradiction.}
	\label{H:fig:return-boundary}
\end{figure}
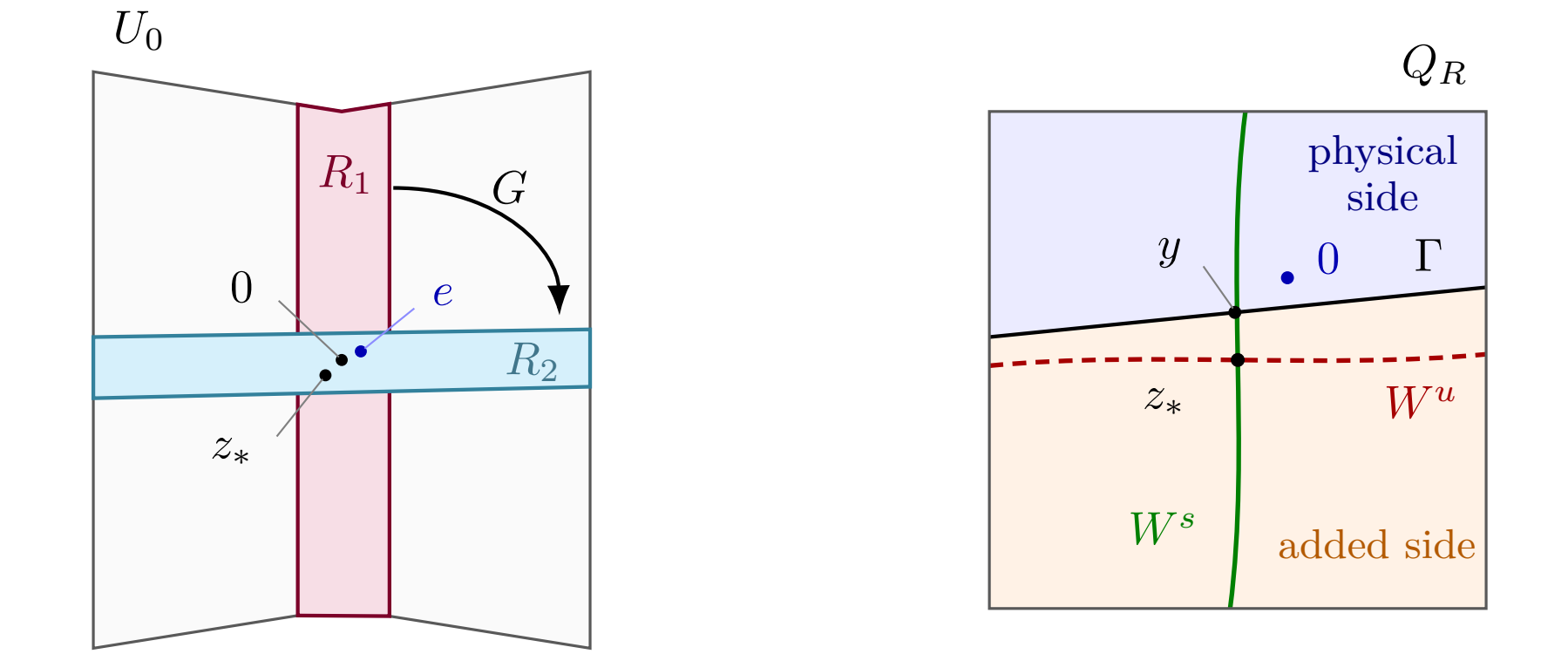

	\Needspace{23\baselineskip}
The geometric argument in Steps 4 and 5 is illustrated
in \Cref{H:fig:return-boundary}.

\section{Derivative growth and non-Brjuno rigidity}\label{H:sec:derivative-growth}

\subsection{A lacunary derivative-growth gap}

In this subsection, we use \cref{thm:finite-closing} to estimate the derivative growth. Explicitly, if the map has no hyperbolic periodic points, then its derivative cannot grow too fast along some particular subsequences of iterates.

\begin{proposition}\label{H:prop:lacunary-gap}
Let $0<\beta\le1$ and let
$f\in\Diff^{1,\beta}_{\mathrm{vol}}(S)$.  Assume that $f$ has no
hyperbolic periodic point.  There are constants
\[
 C_f\ge1,\qquad \Theta\in(0,1),\qquad H_f>1
\]
such that every sequence of integers satisfying
\begin{equation}\label{H:eq:tower}
                    Q_0\ge H_f,\qquad
                    Q_j\ge H_f^{Q_{j-1}}\quad(j\ge1)
\end{equation}
also satisfies
\begin{equation}\label{H:eq:lacunary-gap}
              \log\|Df^{Q_j}\|_0\le C_f\Theta^jQ_j
              \qquad(j\ge0).
\end{equation}
\end{proposition}

\begin{proof}
In Euclidean coordinates the product rule gives
\[
 [D(u\circ v)]_\beta
 \le[Du]_\beta\Lip(v)^\beta\|Dv\|_0
       +\|Du\|_0[Dv]_\beta.
\]
For the Euclidean calculation on a convex coordinate domain, write
$L_f=\max\{2,\|Df\|_0\}$, $B_f=[Df]_\beta$ and
$a_m=[D(f^m)]_\beta$; then
\begin{align*}
 a_{m+1}&\le B_fL_f^{(1+\beta)m}+L_fa_m,\\
 a_m&\le B_f\sum_{j=0}^{m-1}
          L_f^{m-1-j+(1+\beta)j}
       \le B_fmL_f^{(1+\beta)(m-1)}.
\end{align*}
On the surface, we measure
$a_m=[D(f^m)]_\beta$ using the fixed
compatible atlas of \S\ref{H:sec:ambient}.
The chain rule first gives
\[
\|Df^m\|_0\le L_f^m.
\]

For pairs of source points at distance
at most $cL_f^{-m}$, where $c>0$
depends only on the fixed atlas,
the images under $f^m$ lie in a
common buffered coordinate chart.
For more distant pairs, the uniform
derivative bound gives a H\"older
quotient bounded by
$C L_f^{(1+\beta)(m+1)}$.

The preceding composition estimate,
together with the uniform chart
comparisons and the near/far-pair
argument used in
\eqref{eq:atlas-change}, gives
constants $C_1,C_2>0$, independent
of $m$, such that
\[
a_{m+1}
\le C_1a_m+
C_2L_f^{(1+\beta)(m+1)}.
\]
Indeed, the principal composition
term is bounded by
$C B_fL_f^{(1+\beta)m}$,
while the additional coordinate-change
terms are bounded by
$C L_f^{(1+\beta)(m+1)}$.
Iterating this recurrence gives an
exponential bound for $a_m$.
Combining it with
$\|Df^m\|_0\le L_f^m$ and enlarging
$C_f\ge1$, we obtain
\begin{equation}\label{H:eq:iterate-holder}
	\|Df^m\|_0+[D(f^m)]_\beta
	\le e^{C_fm}\qquad(m\ge1).
\end{equation}
 Let $A_*,\vartheta_*,H_*$ be the constants of
\cref{thm:finite-closing} and fix
\[
                         \Theta\in(\vartheta_*,1).
\]

Now choose $H_f$ so large that
\begin{equation}\label{H:eq:Hf-choice}
 H_f\Theta>2,\qquad
 \log H_f>2C_fH_*,\qquad H_f\ge A_*,
\end{equation}
and that, for every $m\ge H_f$,
\begin{equation}\label{H:eq:Hf-absorption}
 \frac{H_f^m}{2m}\ge e^{C_fH_*m},
 \qquad
 (\Theta-\vartheta_*)\frac{(H_f\Theta)^m}{2m}\ge1.
\end{equation}
Suppose that \eqref{H:eq:lacunary-gap} fails for some index, and let
$j$ be the first one.  Since
$\|Df^{Q_0}\|_0\le\|Df\|_0^{Q_0}\le e^{C_fQ_0}$, the index $j=0$
cannot fail, so $j\ge1$.  Put
\[
 m=Q_{j-1},\qquad N=Q_j,\qquad
 q=\left\lfloor\frac Nm\right\rfloor,\qquad
 r=N-qm,\quad 0\le r<m,
\]
and
\[
 g=f^m,\qquad A=e^{C_f\Theta^{j-1}m},\qquad D=e^{C_fm}.
\]
By the minimality of $j$ and \eqref{H:eq:iterate-holder},
\[
                 \|Dg\|_0\le A,\qquad [Dg]_\beta\le D.
\]
The tower condition \eqref{H:eq:tower} gives $m=Q_{j-1}\ge H_f^j$,
hence $\Theta^{j-1}m\ge H_f(H_f\Theta)^{j-1}\ge H_f$, and therefore
\[
                  A\ge A_*,\qquad q\ge D^{H_*}.
\]
For the second of these, $N\ge H_f^m$ and $H_f^m/m\ge2$ give
$\lfloor N/m\rfloor\ge N/(2m)\ge H_f^m/(2m)$, and
\eqref{H:eq:Hf-choice} together with the first inequality in
\eqref{H:eq:Hf-absorption} does the rest.

Since $f^N=f^r\circ g^q$,
\begin{align*}
 \log\|Dg^q\|_0
 &\ge\log\|Df^N\|_0-\log\|Df^r\|_0\\
 &>C_f\Theta^jN-C_fr\\
 &\ge C_f\Theta^jqm-C_fm.
\end{align*}
(When $r=0$ the remainder term vanishes, since $\|Df^0\|_0=1$.)
On the other hand,
\[
 \vartheta_*q\log A
   =\vartheta_*C_f\Theta^{j-1}qm.
\]
Hence \eqref{H:eq:large-growth} follows as soon as
\begin{equation}\label{H:eq:remainder-absorption}
           (\Theta-\vartheta_*)\Theta^{j-1}q\ge1.
\end{equation}
Because $j-1\le m$ and $0<\Theta<1$, we have
$\Theta^{j-1}\ge\Theta^m$, and the left side of
\eqref{H:eq:remainder-absorption} is therefore at least
\[
 (\Theta-\vartheta_*)\frac{(H_f\Theta)^m}{2m}\ge1
\]
by \eqref{H:eq:Hf-absorption}.  Now \cref{thm:finite-closing} applies
and yields a hyperbolic periodic point of $g=f^m$.  That point is also
periodic for $f$, and its derivative at the least $f$-period is
hyperbolic, since a positive power of it is the hyperbolic return for
$g$.  This contradicts the hypothesis on $f$.
\end{proof}

\subsection{Non-Brjuno arithmetic and rigidity}
In this section, we combine the derivative-growth estimate with the arithmetic properties of a non-Brjuno rotation number. We will use two results from \cite{AFLXZ}. 

\begin{lemma}\label{H:lem:arithmetic-tower}
Let $\alpha\in\T\setminus\Q$ be non-Brjuno.  For every $H>1$ there is a
sequence of positive integers $Q_j$ satisfying
\[
                   Q_0\ge H,\qquad Q_{j+1}\ge H^{Q_j},
\]
and, for infinitely many $j$,
\begin{equation}\label{H:eq:small-rotation}
                 \normT{Q_j\alpha}<e^{-Q_j/(j+1)^2}.
\end{equation}
\end{lemma}

This is the zero-indexed version of
\cite[Lemma~2.2]{AFLXZ}.
The proof partitions the continued-fraction
denominators into first-passage blocks
for the threshold $H^q$ and identifies
large jumps in blocks carrying sufficient
Brjuno mass.
Selecting blocks of one parity preserves
the tower inequality between successive
selected denominators.

To ensure $Q_0\ge H$, one may restart
the first-passage construction at a
sufficiently large convergent denominator,
since every tail of the divergent Brjuno
series still diverges.
The block indices are then restarted
as well, and indexing the selected
sequence from $j=0$ gives the
denominator $(j+1)^2$ in
\eqref{H:eq:small-rotation}.

The second ingredient is the free-disk displacement estimate.

\begin{lemma}\label{H:lem:displacement}
	There exists a universal constant
	$C_0>0$ such that every $C^1$
	area-preserving disk pseudo-rotation
	$h$, with
	$\varepsilon=\normT{\rot(h)}$,
	satisfies
	\begin{equation}\label{H:eq:displacement}
		\|h-\id\|_{C^0}
		\le C_0(1+\|Dh\|_0)\sqrt{\varepsilon}.
	\end{equation}
\end{lemma}

For the representative
$+\varepsilon$ this is \cite[Corollary~B]{AFLXZ}, the constant
absorbing the area normalization.  The proof there bounds the area of
a free disk by Franks' returning-disk theorem \cite{Franks} and Kac
recurrence, then applies that bound to disks controlled by the
derivative.  For the representative $-\varepsilon$ we apply the same
result to $h^{-1}$, which gives the same estimate because
\[
 \|h-\id\|_{C^0}=\|h^{-1}-\id\|_{C^0},
 \qquad \|Dh^{-1}\|_0=\|Dh\|_0
\]
for a two-dimensional area-preserving map in the area-compatible
metric.

\begin{proof}[Proof of \cref{thm:rigidity}]
Let $f$ satisfy the hypotheses of the theorem and put $\alpha=\rot(f)$.  The hypothesis of \cref{H:prop:lacunary-gap}
holds here: by the pseudo-rotation property, for every $n\ge1$ the
origin is the unique fixed point of $f^n$ and there is no fixed point
on the boundary, so the Lefschetz--Hopf formula gives
\[
                         \operatorname{ind}(f^n,0)=1.
\]
If $Df(0)$ had positive real reciprocal hyperbolic eigenvalues, then
$\operatorname{sign}\det(I-Df(0))=-1$, against the formula for $n=1$.
If both eigenvalues were negative, the same contradiction would arise
for $f^2$.  So the origin is nonhyperbolic, and $f$ has no hyperbolic
periodic point.

Let $C_f,\Theta,H_f$ be the constants of \cref{H:prop:lacunary-gap}
for the disk, and apply \cref{H:lem:arithmetic-tower} with $H=H_f$; the
resulting tower satisfies
\begin{equation}\label{H:eq:disk-derivative-gap}
                         \|Df^{Q_j}\|_0
                  \le e^{C_f\Theta^jQ_j}.
\end{equation}
For the infinitely many indices satisfying \eqref{H:eq:small-rotation},
$f^{Q_j}$ is a pseudo-rotation with rotation number at distance
$\varepsilon_j=\normT{Q_j\alpha}$ from the integers, and
\cref{H:lem:displacement,H:eq:disk-derivative-gap} give
\[
 \|f^{Q_j}-\id\|_{C^0}
 \le C_0\left(1+e^{C_f\Theta^jQ_j}\right)
       e^{-Q_j/(2(j+1)^2)}.
\]
For all sufficiently large $j$ we have
\[
C_f\Theta^j<\frac{1}{4(j+1)^2}.
\]
Put $a_j=Q_j/(j+1)^2$.
For those indices satisfying
\eqref{H:eq:small-rotation},
the preceding estimate therefore gives
\[
\|f^{Q_j}-\id\|_{C^0}
\le C_0\left(e^{-a_j/2}+e^{-a_j/4}\right).
\]

By \eqref{H:eq:Hf-choice},
$H_f\Theta>2$, hence $H_f>2$.
The tower condition
\eqref{H:eq:tower} then implies,
by induction, that
$Q_j\ge H_f^{j+1}$.
Consequently,
\[
a_j=\frac{Q_j}{(j+1)^2}
\ge\frac{H_f^{j+1}}{(j+1)^2}
\longrightarrow\infty.
\]
Since \eqref{H:eq:small-rotation}
holds for infinitely many indices,
we may choose an increasing sequence
$j_k$ of such indices.
Then $Q_{j_k}\to\infty$ and
\[
\|f^{Q_{j_k}}-\id\|_{C^0}
\longrightarrow0.
\]
Thus $f$ is $C^0$-rigid, as asserted
in \cref{thm:rigidity}.
\end{proof}

\subsection{The \texorpdfstring{$C^1$}{C1} endpoint}
We now explain why the argument above does not directly extend to the \(C^1\) case. The main difficulty is that continuity of the derivative alone does not provide the quantitative estimates needed in the moving-box construction.

The quantitative moving-box argument in
Part~I relies on the H\"older estimate
\eqref{H:eq:coordinate-holder}, which gives
\[
\|Dg_n(z)-Dg_n(0)\|
\le C_{S,\beta}\mathfrak D^2
b_{n+1}|z|^\beta
\qquad (z\in U_n,\ 0\le n<L).
\]

For a map that is only $C^1$, the
oscillation of its derivative on a
compact coordinate domain is controlled
instead by a modulus of continuity.
Fix a coordinate cutoff $r_{\rm mod}>0$,
and let
\[
\varpi:[0,r_{\rm mod}]\to[0,\infty)
\]
be continuous and nondecreasing, with
$\varpi(0)=0$, such that
\[
\|Df(x)-Df(y)\|
\le\varpi(|x-y|)
\]
for pairs of points in a common
compatible source--target coordinate
domain.

For $\varepsilon>0$, define the
generalized inverse scale
\[
r_\varpi(\varepsilon)
:=\sup\{r\in[0,r_{\rm mod}]:
\varpi(r)\le\varepsilon\}.
\]
Continuity and monotonicity imply that
$r_\varpi(\varepsilon)>0$ and
\[
\varpi(r)\le\varepsilon
\qquad
(0\le r\le r_\varpi(\varepsilon)).
\]

However, $C^1$ regularity alone gives
no fixed power-law lower bound on
$r_\varpi(\varepsilon)$ as
$\varepsilon\downarrow0$.
For example, the modulus
\[
\varpi(r)=\frac{1}{\log(e/r)}
\]
for sufficiently small $r>0$,
extended by $\varpi(0)=0$, gives
\[
r_\varpi(\varepsilon)
=e^{1-1/\varepsilon}
\]
for all sufficiently small
$\varepsilon$.
Thus a tolerance of size $B^{-a}$
can require a radius of order
$e^{-B^a}$, smaller than every fixed
negative power of $B$.

By contrast, the $C^{1,\beta}$
construction uses the power-law
scale
\[
\bar r=\mathfrak D^{-5M/\beta},
\]
and \cref{thm:finite-closing}
requires only the polynomial
threshold $q\ge D^{H_*}$.
The lacunary argument of
\cref{H:prop:lacunary-gap}
is designed to absorb such a fixed
power, but it provides no uniform
control over arbitrary super-polynomial
thresholds arising from a general
$C^1$ modulus.

Therefore, the estimates become weaker as \(\beta\) approaches zero, and the present argument does not give a version of \cref{thm:finite-closing} under \(C^1\) regularity alone. The qualitative and
effective-hyperbolicity results
cited here do not by themselves
supply the missing uniform
power-law estimate; compare
\cite{Katok,CP,Pugh}.

\part{Non-rigidity at the \texorpdfstring{$C^1$}{C1} endpoint}\label{part:construction}
The construction of non-rigid examples is based on the method developed in \cite[\S\S2--7]{LWY}. We keep the original local squeeze models and their Hamiltonian realization, but modify the way in which the perturbations are selected. The main new ingredients are an adaptive choice of arithmetic parameters, a sparse phase profile, and a recursive selection of pairwise disjoint supports in the original coordinates. We will recall the necessary estimates and explain the modifications in the following sections.

\section{Adaptive arithmetic and the circle construction}\label{B:sec:circle-construction}

\paragraph{Conventions for the construction.}
Throughout Sections~\ref{B:sec:circle-construction}--\ref{B:sec:completion},
we fix the rotation number $\alpha$ of
\cref{thm:nonrigidity}, the one-dimensional squeeze data of
\cref{B:lem:protected-squeeze}, and the cutoff of \S\ref{sec:hamiltonian}.
Constants may depend on these data but are uniform in the stage $n$,
the micro-scale $r\in[a_n,b_n)$ of \eqref{B:eq:stage-scales}, the tower
heights in \eqref{B:eq:tower-height}, the selected edges, and the
starting index $n_0$ chosen later.  Comparison notation is as in
\S\ref{H:sec:analytic-conventions}.

Differences of circle and annulus maps are measured in compatible
degree-one lifts: the lift of a flow starts at the identity, cumulative
lifts are formed by composition, and the rotation is lifted to
$x\mapsto x+\alpha$. Norms of periodic functions are taken over one
fundamental domain, using the sum $C^1$ norm of \S\ref{sec:conventions};
$\|u\|_\infty=\|u\|_0$ is the supremum norm.  The
\emph{relative closed support} of a function $u$ on $E$ is
\[
 \supp_Eu=\operatorname{cl}_E\{x\in E:u(x)\ne0\},
\]
with $\operatorname{cl}_E$ the closure in $E$; an overline likewise
denotes closure, with the ambient space named when it matters.  For a
map $M$ we write $\supp M$ for $\supp(M-\id)$ in local coordinates, and
arcs in $\T$ are measured with the normalization $|\T|=1$.

\FloatBarrier

\subsection{Adaptive arithmetic checkpoints}\label{B:sec:arithmetic}
In this subsection, we choose arithmetic checkpoints that divide the construction into stages. The resulting stage budgets tend to zero, while their sum diverges. These are the two properties needed in the subsequent construction.

Let $\alpha$ be as in \cref{thm:nonrigidity}, and write
\[
 \theta_n:=\frac{\log q_{n+1}}{q_n},
 \qquad
 \theta_n\longrightarrow0,
 \qquad
 \sum_n\theta_n=\infty.
\]
We shall take the starting index as large as the argument requires.

\Needspace{5\baselineskip}
We begin with denominator estimates that will control the cumulative tower height and the arithmetic correction terms.

\begin{lemma}\label{B:lem:arithmetic}
For every irrational rotation number,
\begin{equation}\label{B:eq:sum-q}
	q_{n+2}\ge2q_n,\qquad
	\sum_{i=0}^n q_i\le4q_n,\qquad
	n=O(\log q_n).
\end{equation}
Moreover
\begin{equation}\label{B:eq:log-q-summable}
 \sum_n\frac{\log q_n}{q_n}<\infty.
\end{equation}
\end{lemma}

\begin{proof}
	Let $\mathfrak a_k\ge1$ denote the
	continued-fraction partial quotients,
	with $q_{-1}=0$ and $q_0=1$.
	The denominator recurrence is
	\[
	q_{n+1}
	=\mathfrak a_{n+1}q_n+q_{n-1}.
	\]
	Since $q_{n+1}\ge q_n$, it follows that
	\[
	q_{n+2}
	\ge q_{n+1}+q_n
	\ge2q_n.
	\]

	Iterating this inequality along the
	even and odd backward subsequences gives
	\[
	\sum_{i=0}^n q_i
	\le 2q_n+2q_{n-1}
	\le4q_n,
	\]
	where $q_{-1}=0$ handles $n=0$.

	The same two-step growth yields
	\[
	q_n\ge2^{\lfloor n/2\rfloor},
	\]
	and hence $n=O(\log q_n)$.

	Finally, $t\mapsto(\log t)/t$ is
	decreasing for $t>e$.
	For all sufficiently large $n$,
	\[
	0\le\frac{\log q_n}{q_n}
	\le C(n+1)2^{-n/2}.
	\]
	The right-hand side is summable,
	which proves
	\eqref{B:eq:log-q-summable}.
\end{proof}

Fix a constant $C_{\rm chk}>0$, to be chosen sufficiently large, and write
\[
 \nu_n:=\frac{\log q_n}{\log 2}.
\]
The factor $1/\log 2$ converts the natural logarithm to the dyadic
scale $\delta_r=2^{-r}$ of \S\ref{B:sec:micro}.  The checkpoints are
\begin{equation}\label{B:eq:checkpoints}
 k_n:=\left\lceil
 \nu_n+\frac{2\log(\nu_n+2)+2\log(n+1)}{\log 2}
           +C_{\rm chk}\right\rceil+n.
\end{equation}
At stage $n$ we put
\begin{equation}\label{B:eq:stage-scales}
 a_n=k_n,\qquad b_n=k_{n+1},\qquad w_n=b_n-a_n,
 \qquad e_n=\frac{w_n}{q_n}.
\end{equation}
The micro-scales of this stage are the integers $a_n\le r<b_n$.  Since $b_n=a_{n+1}$, every integer micro-scale $r\ge a_{n_0}$ belongs
to exactly one stage $n\ge n_0$, and the term
$+n$ guarantees that no stage is empty, even when consecutive
denominators are of comparable size.

\Needspace{6\baselineskip}
The checkpoints make the stage budgets $e_n$ tend to zero while
preserving $\sum_n e_n=\infty$.

\begin{lemma}\label{B:lem:width}
For the prescribed $\alpha$,
\begin{equation}\label{B:eq:width-asymp}
 w_n\ge1,\qquad e_n\longrightarrow0,\qquad
 \sum_ne_n=\infty.
\end{equation}
There is a constant $C_{\log}$,
depending on $C_{\rm chk}$ but not
on $n$, such that
\begin{equation}\label{B:eq:checkpoint-small}
	2^{-a_n}\le
	\frac{2^{-C_{\rm chk}}2^{-n}}
	{q_n(\nu_n+2)^2(n+1)^2},
	\qquad
	a_n+1\le C_{\log}(\nu_n+2).
\end{equation}
Moreover, $C_{\log}$ can be chosen
so that
\[
C_{\log}\le C_*(1+C_{\rm chk}),
\]
where $C_*>0$ is an absolute constant.
In particular,
\begin{equation}\label{B:eq:checkpoint-support}
	q_na_n^22^{-a_n}
	\le C_{\log}^2\,2^{-C_{\rm chk}}
	\frac{2^{-n}}{(n+1)^2}.
\end{equation}
\end{lemma}

\begin{proof}
Write $f_n$ for the expression inside the ceiling in
\eqref{B:eq:checkpoints}.  It is strictly increasing, so that
$w_n=\lceil f_{n+1}\rceil-\lceil f_n\rceil+1\ge1$ and
\[
 w_n>f_{n+1}-f_n
 \ge\nu_{n+1}-\nu_n=\frac{\log q_{n+1}-\log q_n}{\log 2}.
\]
In the other direction, \cref{B:lem:arithmetic} controls $n$ and the
logarithmic correction terms and gives
$w_n\le k_{n+1}\le C(\log q_{n+1}+1)$, whence
$e_n\le C\theta_n+C/q_n\to0$.
Combining the lower bound above with
\eqref{B:eq:log-q-summable}, we obtain
\[
\sum_{n=0}^{J}e_n
>
\frac1{\log2}
\left(
\sum_{n=0}^{J}
\frac{\log q_{n+1}}{q_n}
-
\sum_{n=0}^{J}
\frac{\log q_n}{q_n}
\right).
\]
The first series on the right diverges
by the non-Brjuno assumption, whereas
the second converges by
\eqref{B:eq:log-q-summable}.
Letting $J\to\infty$ gives
$\sum_n e_n=\infty$.

Since $a_n=k_n$ and
$\lceil x\rceil\ge x$, we have
\[
a_n\ge
\nu_n+
\frac{2\log(\nu_n+2)+2\log(n+1)}
{\log2}
+C_{\rm chk}+n.
\]
Using $2^{-\nu_n}=q_n^{-1}$
gives
\[
2^{-a_n}
\le
\frac{2^{-C_{\rm chk}}2^{-n}}
{q_n(\nu_n+2)^2(n+1)^2},
\]
which is the first estimate in
\eqref{B:eq:checkpoint-small}.

Moreover, by \cref{B:lem:arithmetic},
$n\le2\nu_n+1$.
The definition of $k_n$ and the
elementary logarithmic bounds
therefore give
\[
a_n+1
\le C_*(1+C_{\rm chk})(\nu_n+2)
\]
for an absolute constant $C_*>0$.
Thus we may choose
$C_{\log}\le C_*(1+C_{\rm chk})$,
proving the second estimate in
\eqref{B:eq:checkpoint-small}.

Finally, multiplying the first
estimate in \eqref{B:eq:checkpoint-small}
by $q_na_n^2$ and using
$a_n\le C_{\log}(\nu_n+2)$ gives
\eqref{B:eq:checkpoint-support}.
\end{proof}

\subsection{The protected micro-squeeze}
\label{B:sec:micro}
\begin{lemma}\label{B:lem:protected-squeeze}
Set $\kappa_{\rm sq}=e^{-3}$.  There are fixed constants
$\varepsilon_{\rm adm}>0$, $C>0$, and $C_{\rm sq,\ell}>0$,
$0\le\ell\le5$, with the following property.  For
$0<\varepsilon\le\varepsilon_{\rm adm}$, put
\[
 A_\varepsilon=C\varepsilon^{-1}+3>4.
\]
There is a smooth odd vector field $V_\varepsilon$, vanishing near the
ends of $(-A_\varepsilon,A_\varepsilon)$ and extended by zero to $\R$,
whose flow $\Phi_\varepsilon^t$ has the following properties.  With
$P=[-2,2]$, $P^+=[-1,1]$, and $0\le t\le1$,
\begin{align}
 \Phi_\varepsilon^1(P)&=P^+,\label{B:eq:S1}\\
 \kappa_{\rm sq}\le D\Phi_\varepsilon^t&\le e^{\varepsilon/2},\label{B:eq:S2}\\
 D\Phi_\varepsilon^t(z)<1&\Longrightarrow
             \Phi_\varepsilon^t(z)\in P^+.\label{B:eq:S3}
\end{align}
Moreover,
\[
 \|V_\varepsilon\|_{C^\ell(\R)}\le C_{\rm sq,\ell}
       \qquad(1\le\ell\le5),
\]
and for every signed time $|u|\le1$,
\begin{align}
 \|\Phi_\varepsilon^u-\id\|_{C^0}
       &\le C_{\rm sq,0}|u|,\label{B:eq:S40}\\
 \|D\Phi_\varepsilon^u-1\|_{C^0}
       &\le C_{\rm sq,1}|u|,\label{B:eq:S41}\\
 \|D^\ell\Phi_\varepsilon^u\|_{C^0}
       &\le C_{\rm sq,\ell}|u|\qquad(2\le\ell\le5).
       \label{B:eq:S4h}
\end{align}
The constants are independent of $\varepsilon$, and the fixed inner
model also satisfies $V_\varepsilon(0)=0$ and $V_\varepsilon'(0)=-1$.
\end{lemma}

We use exactly the inner model of \cite[proof of Lemma~2.1]{LWY},
namely $-x$ near $0$; hence the two identities at $0$ require no
modification.  The protection implication \eqref{B:eq:S3} is used
only for nonnegative times.

Following \cite[\S2.1]{LWY}, we choose $K_{\min}$ so large that
$\bar\lambda:=e^{1/(2K_{\min}^2)}<3/2$ and that every integer
$r\ge K_{\min}$ satisfies
\begin{equation}\label{B:eq:micro-admissibility}
 \begin{gathered}
 r^{-2}\le\varepsilon_{\rm adm},\qquad
 2(Cr^2+3)2^{-r}<1,\\
 (C(r+1)^2+3)2^{-(r+1)}<(Cr^2+3)2^{-r}.
 \end{gathered}
\end{equation}
We put
\[
 \delta_r=2^{-r},\qquad \varepsilon_r=r^{-2},\qquad
 \lambda_r=e^{\varepsilon_r/2},\qquad A_r=Cr^2+3.
\]
Fix an angular base point $\theta_*\in\T$, written $\theta_0$ in
\cite{LWY}; the scale-$r$ squeeze interval and protected core are
\[
 J_r=\theta_*+\delta_r(-A_r,A_r),\qquad
 P_r=\theta_*+\delta_rP.
\]
By \eqref{B:eq:micro-admissibility}, $J_r$ is an embedded arc of the
circle and $\overline{J_{r+1}}\subset J_r$.  Moreover
$\theta_*+\delta_rP^+=P_{r+1}$.  In the lift of $J_r$ centered at
$\theta_*$ we set
\[
 L_{r,t}(\theta_*+\delta_r z)
       =\theta_*+\delta_r\Phi_{\varepsilon_r}^t(z),
\]
and extend by the identity outside $J_r$.  The protection estimates
rescale accordingly: for $0\le t\le1$,
\begin{align}
 L_{r,1}(P_r)&=P_{r+1},\label{B:eq:micro-full}\\
 \kappa_{\rm sq}\le L_{r,t}'&\le\lambda_r,\label{B:eq:micro-der}\\
 L_{r,t}'(x)<1&\Longrightarrow L_{r,t}(x)\in P_{r+1}.
       \label{B:eq:micro-protect}
\end{align}
For $|u|\le1$ we have
\begin{align}
 \|L_{r,u}-\id\|_{C^0}&\le C_{\rm sq,0}\delta_r|u|,
       \label{B:eq:micro-C0}\\
 \|L_{r,u}'-1\|_{C^0}&\le C_{\rm sq,1}|u|,
       \label{B:eq:micro-C1}\\
 \|L_{r,u}^{(\ell)}\|_{C^0}&\le C_{\rm sq,\ell}
             \delta_r^{1-\ell}|u|\quad(2\le\ell\le5).
       \label{B:eq:micro-high}
\end{align}
These are \cite[(4)--(9)]{LWY}; in particular
$\prod_{r\ge K_{\min}}\lambda_r<\infty$.  We now increase $n_0$ so that
the starting scale $K:=a_{n_0}$ is at least $K_{\min}$; since
$b_n=a_{n+1}$, every scale $r\ge K$ then occurs in exactly one stage.

For integers $K_{\min}\le a<b$ we keep the clipped phase coordinates of
\cite[(10) and Lemma~2.2]{LWY}:
\begin{equation}\label{B:eq:clipped-phase}
 \tau_{a,b,r}(s)=\min\{1,\max\{0,s-(b-1-r)\}\},
 \quad 0\le s\le b-a,\quad a\le r<b.
\end{equation}
The block composite built from them is
\begin{equation}\label{B:eq:block-composite}
 G_{a,b}(s)=L_{b-1,\tau_{a,b,b-1}(s)}\circ\cdots
                         \circ L_{a,\tau_{a,b,a}(s)}.
\end{equation}
Its scale-$r$ phase is the interval $[b-1-r,b-r]$; there, writing
$t_r=s-(b-1-r)$, the composite takes the equivalent right-insertion
form
\begin{equation}\label{B:eq:right-insertion}
 G_{a,b}(s)=L_{b-1,1}\circ\cdots\circ L_{r+1,1}\circ L_{r,t_r}.
\end{equation}
When $r=b-1$ the product preceding $L_{r,t_r}$ is the identity.
The definitions agree at shared phase endpoints, and
\[
 G_{a,b}(0)=\id,\qquad G_{a,b}(b-a)(P_a)=P_b.
\]
From \cite[Lemma~2.2]{LWY} we need two conclusions: if $s,t$ lie in
the same scale-$r$ phase, then
\begin{equation}\label{B:eq:phase-cancel}
 G_{a,b}(t)^{-1}G_{a,b}(s)=L_{r,s-t},
\end{equation}
and, for $0\le s\le b-a$,
\begin{equation}\label{B:eq:block-protect}
 G_{a,b}'(s;x)<1\quad\Longrightarrow\quad G_{a,b}(s)(x)\in P_b.
\end{equation}
The prime denotes differentiation in $x$, the phase parameter being
held fixed.  Since these are the original block composites, both
properties hold as proved there, phase endpoints included.

\subsection{Sparse fronts and the circle construction}\label{B:sec:fronts}
Fix once and for all
\begin{equation}\label{B:eq:c-theta}
 \mathfrak c=\frac1{100},
 \qquad
 \vartheta=\frac1{10000}.
\end{equation}
The construction begins at a sufficiently large continued-fraction index $n_0$; the precise threshold will be raised finitely many times below.

We follow the tower-height convention of \cite[\S3]{LWY}: put
$T_{n_0-1}=0$ and, for $n\ge n_0$,
\begin{equation}\label{B:eq:tower-height}
 N_n:=T_{n-1},\qquad L_n:=\lfloor\mathfrak c q_n\rfloor,
 \qquad T_n:=N_n+L_n.
\end{equation}
In this notation $N_n$ is the old height and $T_n$ the new one.  In
\cite{LWY} a single quantity $D_m=T_m-N_m$ counts both the fresh edges
and the active edges on each front; here the two roles come apart.  The
difference $L_n=T_n-N_n$ counts candidates, while $D_n$ in
\eqref{B:eq:active-count} will count the edges actually selected to be
active.

\begin{definition}\label{B:def:tower-data}
The stage-$n$ rotation tower is
\[
 \mathcal T_n=\{R_\alpha^jJ_{a_n}:j\in\Z,\ |j|\le T_n+2\}.
\]
Its level $j$ is $R_\alpha^jJ_{a_n}$.  As in
\cite[Definition~3.1]{LWY}, an edge is the ordered pair $(j,j+1)$, with
source and target intervals
\[
 J_{n,j}^{\rm src}=R_\alpha^jJ_{a_n},\qquad
 J_{n,j}^{\rm tgt}=R_\alpha^{j+1}J_{a_n}=R_\alpha(J_{n,j}^{\rm src}).
\]
The positive and negative candidate fronts have index sets
\begin{equation}\label{B:eq:candidate-fronts}
 \mathcal J_n^+=\{N_n+1,\ldots,T_n\},\qquad
 \mathcal J_n^-=\{-T_n-1,\ldots,-N_n-2\}.
\end{equation}
Each front contains exactly $L_n$
candidate edges.  We call the collection
of levels $|j|\le N_n+1$ the
\emph{old plateau}.

The one-level offset between the
positive and negative candidate fronts
comes from the orientation
$(j,j+1)$ of an edge.  In particular,
the source and target levels of every
candidate edge lie outside
$|j|\le N_n$, the index range of the
protected block introduced below.

A front is specified by edge indices,
not by a connected arc of the circle.
In contrast with \cite{LWY}, only
selected candidate edges will be active.
\end{definition}

\Needspace{5\baselineskip}
The checkpoint correction makes the tower levels disjoint even after the padding required by the local squeezes.

\begin{lemma}\label{B:lem:padded-separation}
After taking $C_{\rm chk}$ sufficiently large and then $n_0$ sufficiently large,
\begin{equation}\label{B:eq:padded-separation}
 \min_{1\le m\le2T_n+4}\normT{m\alpha}
 >2A_{a_n}\delta_{a_n}
\end{equation}
for every $n\ge n_0$.  Consequently the closures of the intervals
\[
 R_\alpha^jJ_{a_n},\qquad |j|\le T_n+2,
\]
are pairwise disjoint.
\end{lemma}

\begin{proof}
From \eqref{B:eq:sum-q},
\[
T_n\le\mathfrak c\sum_{i=0}^nq_i
\le4\mathfrak c q_n.
\]
Since $8\mathfrak c<1$, increasing $n_0$ gives $2T_n+4<q_n$, and for $1\le m\le2T_n+4$ the best-approximation property yields
\[
 \normT{m\alpha}\ge\normT{q_{n-1}\alpha}
 >\frac1{q_{n-1}+q_n}>\frac1{2q_n}.
\]

On the other hand, by
\eqref{B:eq:checkpoint-small} and
$A_r=Cr^2+3$,
\[
\begin{aligned}
	2A_{a_n}\delta_{a_n}
	&=2(Ca_n^2+3)2^{-a_n}\\
	&\le
	2(C+3)C_{\log}^2\,2^{-C_{\rm chk}}
	\frac{2^{-n}}{q_n(n+1)^2}\\
	&\le
	2(C+3)C_*^2(1+C_{\rm chk})^2
	2^{-C_{\rm chk}}\frac1{q_n}.
\end{aligned}
\]
Since
$(1+C_{\rm chk})^2 2^{-C_{\rm chk}}
\to0$ as $C_{\rm chk}\to\infty$,
we may fix $C_{\rm chk}$ sufficiently
large that
\[
2A_{a_n}\delta_{a_n}
<\frac1{2q_n}
\]
uniformly in $n$.

By \eqref{B:eq:micro-admissibility},
$A_r\delta_r$ is decreasing for
$r\ge K_{\min}$.
After increasing $n_0$ so that
$a_n\ge K_{\min}$ for all
$n\ge n_0$, we obtain
\[
2A_r\delta_r
\le2A_{a_n}\delta_{a_n}
<\frac1{2q_n}
\qquad(r\ge a_n).
\]
Together with the best-approximation
lower bound above, this proves
\eqref{B:eq:padded-separation}.
\end{proof}

Now set
\begin{equation}\label{B:eq:active-count}
 D_n:=\lfloor\vartheta q_n\rfloor.
\end{equation}
By \cref{B:lem:width}, $D_n/w_n\to\infty$, and after a further increase
of $n_0$ we have $D_n\ge2w_n$.
We use a balanced variant of \cite[(14)]{LWY}, distributing
the $D_n$ selected edges on each front separately
over the $w_n$ scales as evenly as
possible: each scale receives
$\lfloor D_n/w_n\rfloor$ edges per
front, and the first
$D_n-w_n\lfloor D_n/w_n\rfloor$
scales receive one additional edge
per front. The
resulting positive integers $d_{n,r}$, $a_n\le r<b_n$, satisfy
\begin{equation}\label{B:eq:balanced-allocation}
 \sum_{r=a_n}^{b_n-1}d_{n,r}=D_n,
 \qquad
 \left|d_{n,r}-\frac{D_n}{w_n}\right|\le1.
\end{equation}
For all sufficiently large $n$ the floor in $D_n$ and the rounding in
\eqref{B:eq:balanced-allocation} give
\begin{equation}\label{B:eq:d-asymp}
 \frac{\vartheta q_n}{4w_n}
 \le d_{n,r}\le\frac{2\vartheta q_n}{w_n},
 \qquad
 \frac{e_n}{2\vartheta}\le\frac1{d_{n,r}}
 \le\frac{4e_n}{\vartheta}.
\end{equation}
Every phase jump therefore tends to zero, uniformly in the scale.

\Needspace{5\baselineskip}

\begin{definition}\label{B:def:sparse-front-profile}
Fix $n\ge n_0$ and suppose that $D_n$ candidate edges have been selected
on each front; that such selections exist and meet the avoidance
requirement is proved recursively in \S\ref{B:sec:selection}.
Enumerate the selected edges in outward order,
\[
 j_{n,1}^+<\cdots<j_{n,D_n}^+\quad\hbox{in }\mathcal J_n^+,
 \qquad
 j_{n,1}^->\cdots>j_{n,D_n}^-\quad\hbox{in }\mathcal J_n^-.
\]
Using \eqref{B:eq:balanced-allocation}, assign the first $d_{n,a_n}$
edges of each list to scale $a_n$, the next $d_{n,a_n+1}$ to scale
$a_n+1$, and so on.  Write $\mathcal I_{n,r}^\pm$ for the selected
scale-$r$ indices on the two fronts and put
\begin{equation}\label{B:eq:active-index-sets}
 \begin{gathered}
 \mathcal I_{n,r}=\mathcal I_{n,r}^+\sqcup\mathcal I_{n,r}^-,\qquad
 \mathcal I_n^\pm=\bigsqcup_{r=a_n}^{b_n-1}\mathcal I_{n,r}^\pm
                                  \subset\mathcal J_n^\pm,\\
 \mathcal I_n=\mathcal I_n^+\sqcup\mathcal I_n^-,\qquad
 |\mathcal I_{n,r}^\pm|=d_{n,r},\qquad |\mathcal I_n^\pm|=D_n.
 \end{gathered}
\end{equation}

The phase profile $\sigma_{n,j}\in[0,w_n]$ is specified by
\[
 \sigma_{n,j}=w_n\quad(|j|\le N_n+1),\qquad
 \sigma_{n,j}=0\quad(|j|\ge T_n+1),
\]
and by the signed jumps
\begin{equation}\label{B:eq:profile-signed-jumps}
 \sigma_{n,j+1}-\sigma_{n,j}=
 \begin{cases}
  -d_{n,r}^{-1},&j\in\mathcal I_{n,r}^+,\\
  +d_{n,r}^{-1},&j\in\mathcal I_{n,r}^-,\\
  0,&j\text{ is not selected}.
 \end{cases}
\end{equation}
As in \cite[Definition~3.1]{LWY}, an edge is \emph{active} when
$\sigma_{n,j}\ne\sigma_{n,j+1}$.  An active edge is \emph{active at
scale $r$} when
\begin{equation}\label{B:eq:active-phase-definition}
 \bigl(\min\{\sigma_{n,j},\sigma_{n,j+1}\},
       \max\{\sigma_{n,j},\sigma_{n,j+1}\}\bigr)
               \subset (b_n-1-r,b_n-r).
\end{equation}
By the alignment checked below, these are exactly the edges indexed by
$\mathcal I_{n,r}$.  On such an edge,
\begin{equation}\label{B:eq:active-jump}
 |\sigma_{n,j+1}-\sigma_{n,j}|=d_{n,r}^{-1}.
\end{equation}
\end{definition}

The definition is consistent: on each front and at each scale,
$\sum_{j\in\mathcal I_{n,r}^\pm}d_{n,r}^{-1}=1$,
so the total decrease in the outward direction is $b_n-a_n=w_n$.  This
matches the prescribed endpoint values, and the profile is determined
uniquely.  Moreover the decrease assigned to scale $r$ takes the phase
from $b_n-r$ down to $b_n-1-r$.  Both phase values of an active
scale-$r$ edge therefore lie in $[b_n-1-r,b_n-r]$, exactly the phase
appearing in \eqref{B:eq:right-insertion}; this gives
\eqref{B:eq:active-phase-definition}, and with it the phase alignment
on which the cancellation rests.  The signed time on an edge will be
\begin{equation}\label{B:eq:signed-edge-time}
 u_{n,j}:=\sigma_{n,j}-\sigma_{n,j+1}
 =\begin{cases}
   d_{n,r}^{-1},&j\in\mathcal I_{n,r}^+,\\
   -d_{n,r}^{-1},&j\in\mathcal I_{n,r}^-,\\
   0,&j\notin\mathcal I_n.
  \end{cases}
\end{equation}
In \cite[Definition~3.1]{LWY}, every fresh edge is active, and
two equivalences available there do not survive for the sparse profile.
The plateau implication $|j|\le N_n+1\Longrightarrow\sigma_{n,j}=w_n$ is
not an equivalence, since an unselected edge may carry a constant phase
into a candidate front, and likewise the intervening levels need not all
have phase strictly between $0$ and $w_n$.

As in the staircase definition just recalled, the active source and
target unions are
\begin{equation}\label{B:eq:circle-active-unions}
 E_n^{\rm src}:=\bigcup_{j\in\mathcal I_n}J_{n,j}^{\rm src},\qquad
 E_n^{\rm tgt}:=R_\alpha(E_n^{\rm src}).
\end{equation}
The circle diffeomorphism $M_n$ is given, as in \cite[(17)--(18)]{LWY}, by
\begin{equation}\label{B:eq:Mn-circle}
 M_n|_{R_\alpha^jJ_{a_n}}
 =R_\alpha^jG_{a_n,b_n}(\sigma_{n,j})R_\alpha^{-j},
 \qquad |j|\le T_n+2,
\end{equation}
and by the identity elsewhere.  The level closures are disjoint by
\cref{B:lem:padded-separation}, and stationarity near the support
endpoints makes the resulting map smooth.  Put
\begin{equation}\label{B:eq:hn-gn}
 h_{n_0}=\id,\qquad h_{n+1}=M_nh_n,\qquad
 g_n=h_n^{-1}R_\alpha h_n.
\end{equation}
The angle $x$ stays fixed throughout the construction; $y=h_n(x)$ is
the stage-$n$ angle, the coordinate in which $g_n$ is a rotation.

\begin{figure}[H]
\centering
\begin{tikzpicture}[x=0.82cm,y=0.90cm,>=Stealth,font=\small]
 \draw[->] (-0.4,0)--(10,0) node[right]{outward level index};
 \draw[->] (0,-0.2)--(0,3.3) node[above]{$\sigma_{n,j}$};
 \draw[densely dashed,gray] (0,2.7)--(9.4,2.7);
 \node[left] at (0,2.7) {$w_n$};
 \draw[thick,blue!65!black]
 (0.15,2.7)--(2,2.7)--(2,2.1)--(4.2,2.1)--(4.2,1.4)
 --(6.5,1.4)--(6.5,0.7)--(8.5,0.7)--(8.5,0)--(9.5,0);
 \foreach \x in {0.5,1.2,2.8,3.6,5.0,5.7,7.3,7.9,9.2}
   \fill[gray] (\x,0) circle (1.6pt);
 \foreach \x in {2,4.2,6.5,8.5}
   \fill[blue!65!black] (\x,0) circle (2.3pt);
 \draw[<->,blue!65!black] (2.22,2.12)--(2.22,2.68)
   node[midway,right,inner xsep=3pt]{$1/d_{n,r}$};
 \node[align=center] at (1.0,-0.55) {old\\plateau};
 \node[blue!65!black] at (5.0,-1.23) {selected edges: small jumps};
 \node[gray!80!black] at (5.0,-1.73)
 {unselected edges: zero phase jump};
 \node[align=center] at (9.15,-0.55) {outside\\the tower};
\end{tikzpicture}
\caption{One sparse front, shown
	schematically in the outward order
	of edge indices (increasing on the
	positive front and decreasing on
	the negative front), not in angular
	order. The phase decreases only
	across selected edges. Each active
	scale-$r$ edge has jump size
	$1/d_{n,r}$, and the $d_{n,r}$
	jumps assigned to that scale have
	total size one. Across an unselected
	edge the adjacent phases agree, so
	the moving-front quotient cancels
	exactly. The drawing suppresses
	many edges and phase steps and is
	not to scale.}
\label{B:fig:sparse-front}
\end{figure}
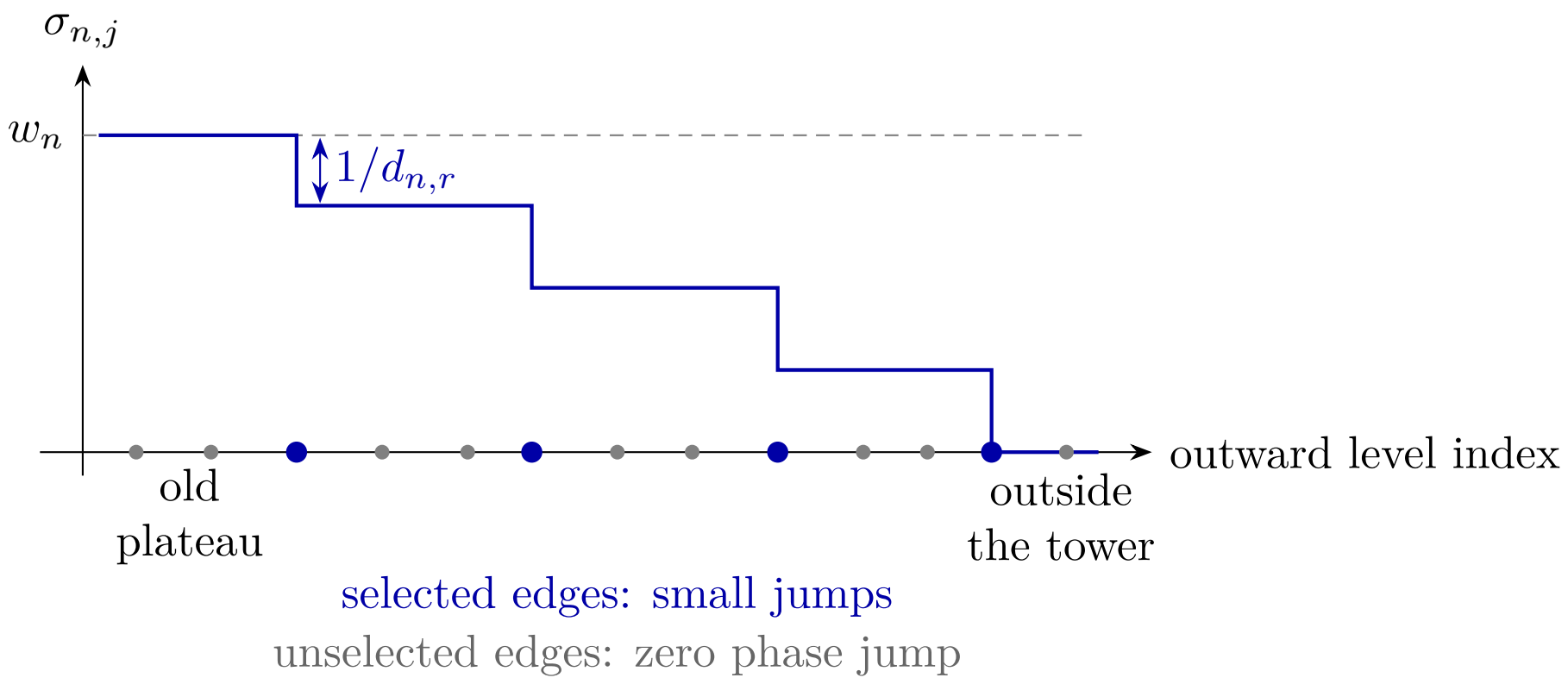

\Needspace{18\baselineskip}
The profile of \cref{B:def:sparse-front-profile} is drawn schematically
in \cref{B:fig:sparse-front}; across an unselected edge the jump is zero.

\Needspace{5\baselineskip}
\begin{lemma}\label{B:lem:sparse-moving-front}
Let $S_n=M_n^{-1}R_\alpha M_n$.  For $-T_n-1\le j\le T_n$,
on the source level $J_{n,j}^{\rm src}=R_\alpha^jJ_{a_n}$,
\begin{equation}\label{B:eq:moving-front}
 S_n=
 R_\alpha^{j+1}
 G_{a_n,b_n}(\sigma_{n,j+1})^{-1}
 G_{a_n,b_n}(\sigma_{n,j})
 R_\alpha^{-j}.
\end{equation}
For $j\in\mathcal I_{n,r}$ the middle quotient is $L_{r,u_{n,j}}$;
for $j\notin\mathcal I_n$ it is the identity.  In particular,
$S_n-R_\alpha$ vanishes outside the selected fine-scale source intervals,
and
\[
 S_n(R_\alpha^jJ_r)=R_\alpha^{j+1}J_r
                  =R_\alpha(R_\alpha^jJ_r)\quad(j\in\mathcal I_{n,r}).
\]
For every $a_n\le r<b_n$ and
$j\in\mathcal I_{n,r}$, the following
estimates hold on the corresponding
coarse source level $J_{n,j}^{\rm src}$:
\begin{align}
	\|S_n-R_\alpha\|_{C^0(J_{n,j}^{\rm src})}
	&\le C\frac{\delta_r}{d_{n,r}},
	\label{B:eq:circle-C0}\\
	\|DS_n-1\|_{C^0(J_{n,j}^{\rm src})}
	&\le C\frac1{d_{n,r}}.
	\label{B:eq:circle-C1}
\end{align}
\end{lemma}

\begin{proof}
This is the local computation of \cite[Lemma~3.1]{LWY}.  The levelwise
formula for $M_n$ gives \eqref{B:eq:moving-front}.  Since our profile
satisfies the same phase-alignment condition
\eqref{B:eq:active-phase-definition}, the cancellation
\eqref{B:eq:phase-cancel} reduces the middle quotient on a selected edge
to $L_{r,u_{n,j}}$.  Equal adjacent phases cancel exactly, in particular
across the unselected edges inside either fresh front.  As $L_{r,u}$ is
supported inside $J_r$ and maps that interval onto itself, the
assertions about support and image follow.  The estimates are
\eqref{B:eq:micro-C0}--\eqref{B:eq:micro-C1}.  In other words, the local
computation never used that every fresh edge is active.
\end{proof}

\subsection{Protected density and the boundary Denjoy mechanism}
In this subsection, we study the conjugacies in the rotation coordinate. As in \cite[\S4]{LWY}, we define the density function and the protected block. 

For an orientation-preserving circle
diffeomorphism $h$ let
\[
 \rho_h(y)=h'(h^{-1}(y)).
\]
Write $\rho_n=\rho_{h_n}$.  The protected block is
\begin{equation}\label{B:eq:protected-block}
 \mathcal P_n=\bigcup_{|j|\le N_n}R_\alpha^jP_{a_n}.
\end{equation}
This is a union of protected cores; the old plateau, by contrast,
consists of the full levels with $|j|\le N_n+1$.

\begin{proposition}\label{B:prop:protected-density}
For every $n\ge n_0$,
\begin{equation}\label{B:eq:rho-lower}
 \rho_n(y)\ge1
 \qquad(y\in\T\setminus\mathcal P_n).
\end{equation}
Every active source level and its target level are disjoint from $\mathcal P_n$.
\end{proposition}

\begin{proof}
The argument is the induction of \cite[Proposition~4.1]{LWY}, whose
inputs are the full plateau, block protection, and the separation of
the active source and target levels from the protected block.  All three
remain valid for the sparse profile.  At the
start $h_{n_0}=\id$ and the density is one.  Assume the claim at stage
$n$ and write $y=M_n(z)$, so that
\begin{equation}\label{B:eq:rho-rec}
 \rho_{n+1}(y)=M_n'(z)\rho_n(z).
\end{equation}
If $z\in\mathcal P_n$, its level lies in the full plateau, where the full block composite carries $P_{a_n}$ onto $P_{b_n}=P_{a_{n+1}}$; hence $M_n(z)\in\mathcal P_{n+1}$.  So if $y\notin\mathcal P_{n+1}$, then $z\notin\mathcal P_n$ and $\rho_n(z)\ge1$.

Suppose in addition that
$M_n'(z)<1$.
Then $z$ lies in a tower level
$R_\alpha^jJ_{a_n}$ with
$\sigma_{n,j}>0$, since $M_n$
is the identity outside the tower
and $G_{a_n,b_n}(0)=\id$.
By the definition of the phase
profile, $|j|\le T_n=N_{n+1}$.
Block protection
\eqref{B:eq:block-protect}
then gives
\[
M_n(z)\in R_\alpha^jP_{b_n}
\subset\mathcal P_{n+1},
\]
contrary to
$y\notin\mathcal P_{n+1}$.
It follows that $M_n'(z)\ge1$,
and \eqref{B:eq:rho-rec} gives
\eqref{B:eq:rho-lower}
at the next stage.

Finally, every active edge belongs to a fresh front, so both its source index and its target index lie outside $|j|\le N_n$, and by padded separation the corresponding levels are disjoint from \eqref{B:eq:protected-block}.
\end{proof}

\Needspace{5\baselineskip}

\begin{lemma}\label{B:lem:stabilized}
For each $j\in\Z$, choose $n$ so large that $|j|\le N_n$ and define
\begin{equation}\label{B:eq:Ij}
 I_j=h_n^{-1}(R_\alpha^jP_{a_n}).
\end{equation}
Then $I_j$ is independent of every sufficiently large choice of $n$.  The intervals $I_j$ are pairwise disjoint and nondegenerate, and
\begin{equation}\label{B:eq:shift-Ij}
 g_n(I_j)=I_{j+1}
\end{equation}
for all sufficiently large $n$ depending on $j$.
\end{lemma}

\begin{proof}
The stabilization argument of \cite[Lemma~4.1]{LWY} applies, because
the full-plateau action it requires is unchanged: on every protected old
level the coefficient is $w_n$, so
\[
 M_n(R_\alpha^jP_{a_n})
 =R_\alpha^jP_{b_n}
 =R_\alpha^jP_{a_{n+1}}.
\]
Taking inverse images under $h_{n+1}=M_nh_n$ gives the stabilization.  Pairwise disjointness and nondegeneracy are checked at one common late stage, using padded separation.  Once $j$ and $j+1$ are both protected,
\[
 g_n(I_j)
 =h_n^{-1}R_\alpha(R_\alpha^jP_{a_n})
 =I_{j+1}.
\qedhere
\]
\end{proof}

\section{Hamiltonian realization and sparse support selection}\label{B:sec:hamiltonian-selection}

\subsection{Hamiltonian realization and derivative estimates}
\label{sec:hamiltonian}
In this subsection, we extend the above circle construction to an area-preserving Hamiltonian flow on the annulus as in \cite[\S\S5--6]{LWY}. We also establish the derivative estimates needed in the proof.

Let
\[
 \mathbb A=\T\times[0,1],
 \qquad
 \omega=dx\wedge ds,
\]
where $s=0$ is the outer boundary.  We use the flat product metric,
write $\mathcal R_\alpha(x,s)=(x+\alpha,s)$, and compute all periodic
norms on $[0,1]\times[0,1]$ in compatible lifts.  The cutoff is that
of \cite[(3), \S5]{LWY}: we fix $\mathfrak b=7/8$ and a smooth
$\chi:[0,\infty)\to[0,1]$ which equals one on $[0,\mathfrak b]$ and
vanishes near $[1,\infty)$.  Then $\bar\lambda/2<\mathfrak b<1$.

\Needspace{5\baselineskip}
We distinguish the full support collar from the smaller region where the flow is exactly a cotangent lift.

\begin{definition}
\label{B:def:collar-data}
For $r\ge K_{\min}$, the scale-$r$ boundary collar and its local
support rectangle are, as in \cite[Definition~5.1]{LWY},
\[
 \T\times[0,\delta_r],\qquad Q_r:=J_r\times[0,\delta_r].
\]
The exact cotangent zone is $\T\times[0,\mathfrak b\delta_r]$; its
intersection with $Q_r$ is $J_r\times[0,\mathfrak b\delta_r]$.
The cotangent formula below applies along an orbit segment only when
the entire segment stays in that zone.
\end{definition}

At scale $r$ we put
\[
 V_{r,\delta_r}(x)=
 \delta_rV_{\varepsilon_r}\!\left(\frac{x-\theta_*}{\delta_r}\right)
\]
on $J_r$, extended by zero, and define
\begin{equation}\label{B:eq:Hamiltonian}
 K_r(x,s)=sV_{r,\delta_r}(x)\chi(s/\delta_r).
\end{equation}
We write $\Lambda_{r,t}$ for its Hamiltonian flow, with the sign
convention $X_H=(\partial_sH,-\partial_xH)$.

\Needspace{5\baselineskip}
The cutoff Hamiltonian realizes the circle squeeze by an area-preserving flow, with uniform estimates at the physical scale.

\begin{lemma}\label{B:lem:ham-squeeze}
The flow $\Lambda_{r,t}$ is smooth, area-preserving, and supported in $J_r\times[0,\delta_r]$.  On the outer boundary,
\[
 \Lambda_{r,t}(x,0)=(L_{r,t}(x),0).
\]
Whenever the whole orbit segment stays in $0\le s\le\mathfrak b\delta_r$,
\begin{equation}\label{B:eq:cotangent-local}
 \Lambda_{r,t}(x,s)
 =\left(L_{r,t}(x),\frac{s}{L_{r,t}'(x)}\right).
\end{equation}
Moreover, for $|u|\le1$ and spatial orders $\ell=0,1,2,3$,
\begin{equation}\label{B:eq:ham-scale}
 \|D^\ell(\Lambda_{r,u}-\id)\|_{C^0}
 \le C_\ell\delta_r^{1-\ell}|u|,
\end{equation}
where $D^0$ denotes the map itself.  The one-dimensional derivatives
also satisfy, for $0\le u\le t\le1$,
\begin{equation}\label{B:eq:micro-derivative-ratio}
 \frac{L_{r,t}'(x)}{L_{r,u}'(x)}\le\lambda_r.
\end{equation}
\end{lemma}

This is the Hamiltonian realization of \cite[Lemma~5.1]{LWY}, with the
same cutoff and sign convention.  As there, the exact cotangent formula
will only be invoked along orbit segments that stay entirely inside the
exact zone.

For $a_n\le r<b_n$ and $j\in\Z$ we keep the phase coefficient of
\cite[(32)]{LWY}, namely
\begin{equation}\label{B:eq:ham-phase-coefficient}
 c_{n,r,j}:=\tau_{a_n,b_n,r}(\sigma_{n,j})\in[0,1].
\end{equation}
Let $\mathcal N_{n,r}$ be the map which equals
$\mathcal R_\alpha^j\Lambda_{r,c_{n,r,j}}\mathcal R_\alpha^{-j}$ on
$\mathcal R_\alpha^jQ_r$ for $|j|\le T_n+2$, and the identity
elsewhere.  At fixed $r$ these rectangles have disjoint closures, and
in fact $\mathcal N_{n,r}$ is the time-one map of the smooth
Hamiltonian
\[
 \sum_{|j|\le T_n+2}c_{n,r,j}K_r\circ\mathcal R_\alpha^{-j},
\]
whose summands have disjoint supports.  Following \cite[(33)--(34)]{LWY}
we set
\begin{equation}\label{B:eq:macro-Ham}
 \mathcal M_n=\mathcal N_{n,b_n-1}\circ\cdots\circ\mathcal N_{n,a_n},
 \qquad \mathcal H_{n_0}=\id,\qquad
 \mathcal H_{n+1}=\mathcal M_n\mathcal H_n,
\end{equation}
and
\begin{equation}\label{B:eq:Fn}
 F_n=\mathcal H_n^{-1}\mathcal R_\alpha\mathcal H_n.
\end{equation}
On the boundary, $\mathcal M_n$, $\mathcal H_n$, and $F_n$ restrict to
$M_n$, $h_n$, and $g_n$.  Since the Hamiltonian maps are produced
directly by this finite factorization, the question whether a block
path is differentiable across the phase endpoints does not arise.

Besides the fixed physical coordinates we use the coordinates in which
$F_n$ is a rotation,
\begin{equation}\label{B:eq:physical-rotation-coordinates}
 (y,v)=\mathcal H_n(x,s),\qquad
 \mathcal H_n\circ F_n\circ\mathcal H_n^{-1}(y,v)=(y+\alpha,v).
\end{equation}
So $(x,s)$ are the physical coordinates and $(y,v)$ the rotation
coordinates at stage $n$, and on the outer boundary $y=h_n(x)$.

\Needspace{5\baselineskip}
To estimate the conjugated defects, we record both the active rectangles in rotation coordinates and their preimages in physical coordinates.

\begin{definition}
\label{B:def:active-rectangles}
For $j\in\mathcal I_{n,r}$, the active source and target rectangles are
those of \cite[Definition~5.1]{LWY}:
\[
 Q_{n,r,j}^{\rm src}:=\mathcal R_\alpha^jQ_r,\qquad
 Q_{n,r,j}^{\rm tgt}:=\mathcal R_\alpha^{j+1}Q_r.
\]
They lie in rotation coordinates.  Write
\[
 \mathcal E_n^{\rm src}:=\bigcup_{r=a_n}^{b_n-1}
                \bigcup_{j\in\mathcal I_{n,r}}Q_{n,r,j}^{\rm src},
 \qquad
 \mathcal E_n^{\rm tgt}:=\mathcal R_\alpha(\mathcal E_n^{\rm src}).
\]
Their physical components are
$\mathcal H_n^{-1}(Q_{n,r,j}^{\rm src})$ and
$\mathcal H_n^{-1}(Q_{n,r,j}^{\rm tgt})$, respectively.
\end{definition}

We also keep the notation of \cite[Lemma~5.2]{LWY} for the map
\begin{equation}\label{B:eq:ham-moving-map}
 \mathcal S_n:=\mathcal M_n^{-1}\mathcal R_\alpha\mathcal M_n.
\end{equation}
The map $\mathcal S_n$ is a diffeomorphism; $\mathcal S_n-\mathcal R_\alpha$
denotes its defect in compatible lifts.
For $K_{\min}\le a<b$ and $0\le s\le b-a$, the two-dimensional block
composite of that lemma is
\[
 \mathcal G_{a,b}(s)
 :=\Lambda_{b-1,\tau_{a,b,b-1}(s)}\circ\cdots
                                  \circ\Lambda_{a,\tau_{a,b,a}(s)}.
\]
On the coarse source level $\mathcal R_\alpha^jQ_{a_n}$,
$-T_n-1\le j\le T_n$, the same levelwise computation as there gives
\begin{equation}\label{B:eq:ham-moving-front}
 \mathcal S_n
 =\mathcal R_\alpha^{j+1}
   \mathcal G_{a_n,b_n}(\sigma_{n,j+1})^{-1}
   \mathcal G_{a_n,b_n}(\sigma_{n,j})\mathcal R_\alpha^{-j}.
\end{equation}
For $j\in\mathcal I_{n,r}$, phase alignment reduces the middle quotient
to $\Lambda_{r,u_{n,j}}$, while for an inactive edge the quotient is
the identity.  Since this flow preserves $Q_r$,
\begin{equation}\label{B:eq:active-rectangle-onto}
 \mathcal S_n(Q_{n,r,j}^{\rm src})
 =\mathcal R_\alpha(Q_{n,r,j}^{\rm src})=Q_{n,r,j}^{\rm tgt}.
\end{equation}
Thus \cite[Lemma~5.2]{LWY} applies
component by component.
The relative closed support of
$\mathcal S_n-\mathcal R_\alpha$
lies in $\mathcal E_n^{\rm src}$.
For every $a_n\le r<b_n$ and
$j\in\mathcal I_{n,r}$, the following
estimates hold on the corresponding
active source rectangle
$Q_{n,r,j}^{\rm src}$:
\begin{align}
	\|\mathcal S_n-\mathcal R_\alpha\|
	_{C^0(Q_{n,r,j}^{\rm src})}
	&\le C\frac{\delta_r}{d_{n,r}},
	\label{B:eq:S2d0}\\
	\|D\mathcal S_n-I\|
	_{C^0(Q_{n,r,j}^{\rm src})}
	&\le C\frac1{d_{n,r}},
	\label{B:eq:S2d1}\\
	\|D^2\mathcal S_n\|
	_{C^0(Q_{n,r,j}^{\rm src})}
	&\le C\frac{\delta_r^{-1}}{d_{n,r}}.
	\label{B:eq:S2d2}
\end{align}
Here $I$ is the identity matrix.  The target rectangle in
\eqref{B:eq:active-rectangle-onto} is also where the joining segments
of the derivative estimate below will lie.

\Needspace{5\baselineskip}
All earlier factors remain inside their exact cotangent zones when traced backward from the current collar.

\begin{lemma}\label{B:lem:nested-cotangent}
On the full rotation-coordinate strip $0\le v\le\delta_{a_n}$,
\begin{equation}\label{B:eq:H-cotangent}
 \mathcal H_n^{-1}(y,v)
 =\bigl(h_n^{-1}(y),v\rho_n(y)\bigr),
\end{equation}
and equivalently
\begin{equation}\label{B:eq:H-cotangent-forward}
 \mathcal H_n(x,s)
 =\left(h_n(x),\frac{s}{h_n'(x)}\right)
\end{equation}
on the physical preimage of that strip.  Thus the inverse formula holds
on every active rotation-coordinate rectangle, and the forward formula
on its physical preimage.
\end{lemma}

\begin{proof}
	We follow the proof of
	\cite[Lemma~5.3]{LWY}.
	Write $a=a_n$.
	Each micro-scale $K\le r<a$
	used so far occurs once.
	We trace the inverse factors
	in decreasing order of scale,
	$r=a-1,\ldots,K$.

	Assume inductively that, upon
	entering the inverse scale-$r$
	factor, the normal height is at most
	\[
	s_{\rm in}\le
	\delta_a\prod_{u=r+1}^{a-1}\lambda_u.
	\]
	For $r=a-1$ this holds because
	the product is empty.
	Since $\lambda_u\le\bar\lambda$
	and $\delta_a=2^{-(a-r)}\delta_r$,
	\[
	s_{\rm in}
	\le\frac{\delta_r}{2}
	\left(\frac{\bar\lambda}{2}
	\right)^{a-r-1}
	\le\frac12\delta_r
	<\mathfrak b\delta_r.
	\]
	Thus the backward segment starts
	strictly inside the exact
	cotangent zone.

	By \eqref{B:eq:micro-derivative-ratio},
	every partial-time inverse factor
	expands the normal direction by
	at most $\lambda_r$ as long as
	the segment stays in that zone.
	Hence, up to a possible first exit,
	its height is bounded by
	\[
s\le\lambda_r s_{\rm in}\le\delta_a\prod_{u=r}^{a-1}\lambda_u\le\delta_r\left(\frac{\bar\lambda}{2}\right)^{a-r}\le\frac{\bar\lambda}{2}\delta_r<\mathfrak b\delta_r.
\]
	The strict inequality rules out
	a first exit. Consequently, the
	whole inverse scale-$r$ factor
	stays in its exact cotangent zone,
	and its final height satisfies
	the induction bound for the
	next, coarser scale.

	Outside its angular support
	the factor is the identity,
	so the argument is independent
	of which edges were selected.
	All earlier factors are therefore
	exact cotangent lifts along the
	orbit in question.
	Composing them gives
	\eqref{B:eq:H-cotangent}, and
	inverting gives
	\eqref{B:eq:H-cotangent-forward}
	on the physical preimage of
	the strip.
\end{proof}

The derivative bounds for the boundary maps come from the auxiliary
micro-factorization of \cite[\S6.1]{LWY}.  Let $N_{n,r}:\T\to\T$ be
given by
\[
 \mathcal N_{n,r}(x,0)=(N_{n,r}(x),0).
\]
The double index distinguishes it from the integer old height $N_n$.
If $n$ is the unique stage containing $r$, set
\begin{equation}\label{B:eq:auxiliary-micro-factorization}
 \widehat h_K=\id,\qquad
 \widehat h_{r+1}=N_{n,r}\widehat h_r\quad(a_n\le r<b_n).
\end{equation}
Because the scale blocks are consecutive this is unambiguous, and
$h_n=\widehat h_{a_n}$, $h_{n+1}=\widehat h_{b_n}$.  Each $N_{n,r}$
agrees with a translate of $L_{r,c_{n,r,j}}$ on disjoint levels and
with the identity elsewhere, so
\[
 \|N_{n,r}'\|_\infty\le\lambda_r,\qquad
 \|N_{n,r}^{(\ell)}\|_\infty\le C_\ell\delta_r^{1-\ell}
                         \quad(2\le\ell\le4).
\]
The micro-scale induction of \cite[Lemma~6.1]{LWY} uses nothing else;
where the active levels sit, and how many of them there are, plays no
role.  At $r=a_n$ that lemma gives
\begin{equation}\label{B:eq:Bder}
 \|h_n'\|_\infty\le B_{\rm der}
       :=\exp\!\left(\frac12\sum_{r\ge1}r^{-2}\right),
\end{equation}
and, for $\ell=2,3,4$,
\begin{equation}\label{B:eq:h-high}
 \|h_n^{(\ell)}\|_\infty\le C_\ell\delta_{a_n}^{1-\ell}.
\end{equation}
Both bounds are uniform in the tower heights, in the sparse selection,
and in later increases of $n_0$.

On the angular projections of the active source and target rectangles,
\cref{B:prop:protected-density} together with \eqref{B:eq:Bder} gives
$1\le\rho_n\le B_{\rm der}$, and by the inverse cotangent formula
$s=v\rho_n(y)\le B_{\rm der}\delta_{a_n}$.
Thus \eqref{B:eq:Bder}--\eqref{B:eq:h-high}
supply the local cotangent derivative
bounds of \cite[Lemma~6.2]{LWY}.
For any active source or target
rectangle $Q$, put
$\Omega=\mathcal H_n^{-1}(Q)$.
Then
\begin{align}
	\|D\mathcal H_n\|_{C^0(\Omega)}
	+\|D\mathcal H_n^{-1}\|_{C^0(Q)}
	&\le C,\label{B:eq:cot1}\\
	\|D^2\mathcal H_n\|_{C^0(\Omega)}
	+\|D^2\mathcal H_n^{-1}\|_{C^0(Q)}
	&\le C\delta_{a_n}^{-1},
	\label{B:eq:cot2}\\
	\|D^3\mathcal H_n\|_{C^0(\Omega)}
	+\|D^3\mathcal H_n^{-1}\|_{C^0(Q)}
	&\le C\delta_{a_n}^{-2}.
	\label{B:eq:cot3}
\end{align}
Since a lifted active target
rectangle is convex, the inverse
estimates also hold along every
segment contained in it.

Let
\begin{equation}\label{B:eq:increment-definition}
 \Delta_n:=F_{n+1}-F_n
\end{equation}
in compatible degree-one annulus lifts.  Differentiating
\[
 F_{n+1}=\mathcal H_n^{-1}\mathcal S_n\mathcal H_n,
 \qquad F_n=\mathcal H_n^{-1}\mathcal R_\alpha\mathcal H_n,
\]
and combining \eqref{B:eq:S2d0}--\eqref{B:eq:S2d2} with
\eqref{B:eq:cot1}--\eqref{B:eq:cot3}, we arrive at the following bounds
on the increments.

\begin{proposition}\label{B:prop:ledger}
	There are constants $C_0,C_1,C_2>0$,
	independent of $n$, $r$, and the sparse
	selection, with the following property.

	For every $a_n\le r<b_n$ and
	$j\in\mathcal I_{n,r}$, put
	\[
	\Omega=\mathcal H_n^{-1}
	(Q_{n,r,j}^{\rm src}).
	\]
	Then
	\begin{align}
		\|\Delta_n\|_{C^0(\Omega)}
		&\le C_0\frac{\delta_r}{d_{n,r}},
		\label{B:eq:Delta0}\\
		\|D\Delta_n\|_{C^0(\Omega)}
		&\le C_1\frac1{d_{n,r}},
		\label{B:eq:Delta1}\\
		\|D^2\Delta_n\|_{C^0(\Omega)}
		&\le C_2\frac{\delta_r^{-1}}{d_{n,r}}.
		\label{B:eq:Delta2}
	\end{align}
	The defect vanishes in a relative
	neighborhood of the lateral and inner
	edges of each physical source component,
	so each component admits a smooth
	zero extension.
\end{proposition}

\begin{proof}
We run the local calculation from the proof of
\cite[Proposition~6.1]{LWY} on each selected active component.  That
calculation uses only the single-flow defect bounds
\eqref{B:eq:S2d0}--\eqref{B:eq:S2d2}, the cotangent estimates
\eqref{B:eq:cot1}--\eqref{B:eq:cot3} on the active region, and
$\delta_r\le\delta_{a_n}$, all of which have been established above.
Given a physical source point $z$, let $w=\mathcal H_n(z)$.  By
\eqref{B:eq:active-rectangle-onto} the points $\mathcal S_n(w)$ and
$\mathcal R_\alpha(w)$ lie in a common lifted target rectangle, so the
inverse estimates are available along the segment joining them.

For clarity, set $U=\mathcal H_n^{-1}$ and write out the
first-derivative estimate.  Since $D\mathcal R_\alpha=I$ in these
coordinates, the chain rule gives
\begin{align*}
 D\Delta_n(z)
 &=\Bigl(DU(\mathcal S_n(w))D\mathcal S_n(w)
             -DU(\mathcal R_\alpha(w))\Bigr)D\mathcal H_n(z)\\
 &=\Bigl(DU(\mathcal S_n(w))-DU(\mathcal R_\alpha(w))\Bigr)
          D\mathcal H_n(z)\\
 &\quad+DU(\mathcal S_n(w))(D\mathcal S_n(w)-I)
          D\mathcal H_n(z).
\end{align*}
The joining segment remains inside the lifted target rectangle, so
\eqref{B:eq:cot2} and \eqref{B:eq:S2d0} imply
\[
 \|DU(\mathcal S_n(w))-DU(\mathcal R_\alpha(w))\|
 \le C\delta_{a_n}^{-1}
             |\mathcal S_n(w)-\mathcal R_\alpha(w)|
 \le C\frac{\delta_r}{\delta_{a_n}d_{n,r}}.
\]
Together with \eqref{B:eq:cot1} and \eqref{B:eq:S2d1}, this yields
\[
 \|D\Delta_n(z)\|
 \le C\left(\frac{\delta_r}{\delta_{a_n}}+1\right)
                    \frac1{d_{n,r}}
 \le\frac{2C}{d_{n,r}},
\]
because $\delta_r\le\delta_{a_n}$.  This proves
\eqref{B:eq:Delta1} on the component using only the local cotangent
bounds.

The zeroth-order bound follows
from the mean-value estimate along
the joining segment:
\[
|\Delta_n(z)|
\le C|\mathcal S_n(w)
-\mathcal R_\alpha(w)|
\le C\frac{\delta_r}{d_{n,r}}.
\]
This proves \eqref{B:eq:Delta0}.

For the second derivative, put
\[
\begin{gathered}
	w_1=\mathcal S_n(w),\qquad
	w_0=\mathcal R_\alpha(w),\\
	P=D\mathcal H_n(z),\qquad
	Q=D^2\mathcal H_n(z),\\
	E=D\mathcal S_n(w)-I,\qquad
	Z=D^2\mathcal S_n(w).
\end{gathered}
\]
The second-order chain rule gives
\[
\begin{aligned}
	D^2\Delta_n(z)
	={}&D^2U(w_1)[(I+E)P,(I+E)P]\\
	&-D^2U(w_0)[P,P]\\
	&+DU(w_1)Z[P,P]\\
	&+\bigl(DU(w_1)(I+E)
	-DU(w_0)\bigr)Q.
\end{aligned}
\]
The bilinear notation is understood
componentwise.

Write
\[
a=\delta_{a_n},\qquad
b=\delta_r,\qquad d=d_{n,r}.
\]
By \eqref{B:eq:cot1}--\eqref{B:eq:cot3}
and \eqref{B:eq:S2d0}--\eqref{B:eq:S2d2},
\[
\begin{gathered}
	\|P\|\le C,\quad
	\|Q\|\le Ca^{-1},\\
	\|E\|\le Cd^{-1},\quad
	\|Z\|\le Cb^{-1}d^{-1}.
\end{gathered}
\]
Since the joining segment lies
inside the lifted target rectangle,
the same estimates imply
\[
\begin{aligned}
	\|DU(w_1)-DU(w_0)\|
	&\le C\frac{b}{ad},\\
	\|D^2U(w_1)-D^2U(w_0)\|
	&\le C\frac{b}{a^2d}.
\end{aligned}
\]
Expanding the bilinear terms and
combining these inequalities yields
\[
\|D^2\Delta_n(z)\|
\le C\left(
\frac1{bd}
+\frac1{ad}
+\frac{b}{a^2d}
+\frac1{ad^2}
\right).
\]
Because $b\le a$ and $d\ge1$,
the right-hand side is bounded by
$C/(bd)$, proving
\eqref{B:eq:Delta2}.

The calculation uses only local
cotangent bounds and single-flow
defect estimates. It requires
neither summability nor contiguity
of the active edges, nor any global
derivative bound for the cumulative
conjugacy.

Finally, the single-flow defect
vanishes near the lateral and inner
edges of the source rectangle.
Conjugating by the smooth
diffeomorphism $\mathcal H_n$
preserves this relative
stationarity, giving the asserted
smooth zero extensions.
\end{proof}

Since $\Delta_n$ vanishes outside
the physical active source components,
the local bounds
\eqref{B:eq:Delta0}--\eqref{B:eq:Delta1}
and \eqref{B:eq:d-asymp} give
\begin{equation}\label{B:eq:Delta-C1-vanish}
	\|\Delta_n\|_{C^1}
	\le C\max_{a_n\le r<b_n}
	\frac1{d_{n,r}}
	\le Ce_n\longrightarrow0.
\end{equation}

\subsection{Sparse support selection}\label{B:sec:selection}
In this subsection, we study edges at each stage, and we will choose the edges in \cref{B:def:sparse-front-profile} so that the increments from different stages have disjoint supports in the fixed physical coordinates. This section plays an important role in \cref{B:sec:completion} when we prove \(C^1\) convergence. 

Recall the coordinate relation
\eqref{B:eq:physical-rotation-coordinates}.  A selected scale-$r$ edge
$j$ has source rectangle
$Q_{n,r,j}^{\rm src}=\mathcal R_\alpha^j(J_r\times[0,\delta_r])$ in
rotation coordinates, and the cotangent formula gives
\[
 \mathcal H_n^{-1}(Q_{n,r,j}^{\rm src})
 \subset h_n^{-1}(R_\alpha^j\overline{J_r})\times[0,1].
\]
An angular preimage under $h_n$ therefore covers the physical support
of the corresponding increment.  It is these covers, rather than the
stage-dependent rotation intervals, that must be kept disjoint;
together with \eqref{B:eq:Delta-C1-vanish}, this will give $C^1$
convergence in \S\ref{B:sec:completion}.

\begin{definition}\label{B:def:physical-support-covers}
At a stage $i$ whose active edges have already been selected, let
\[
 \mathcal A_i=\{(r,j):a_i\le r<b_i,\
                    j\in\mathcal I_{i,r}\},
 \qquad |\mathcal A_i|=2\sum_{r=a_i}^{b_i-1}d_{i,r}=2D_i.
\]
Here $j$ always denotes the source index of the edge $(j,j+1)$.
For $(r,j)\in\mathcal A_i$, define the closed physical angular arc
\[
 K_{i,r,j}:=h_i^{-1}(R_\alpha^j\overline{J_r}).
\]
Set
\begin{equation}\label{B:eq:rotation-support-enlarged}
 E_i^*:=\bigcup_{(r,j)\in\mathcal A_i}R_\alpha^j\overline{J_r},
\end{equation}
\begin{equation}\label{B:eq:Omega-i}
 \Omega_i:=h_i^{-1}(E_i^*)=\bigcup_{(r,j)\in\mathcal A_i}K_{i,r,j},
 \qquad \mathrm{Len}_i:=\sum_{(r,j)\in\mathcal A_i}|K_{i,r,j}|.
\end{equation}
Here $\Omega_i$ is a closed angular cover in the fixed physical
coordinate, and $\mathrm{Len}_i$ is its summed arc length, overlaps
counted separately; it is not a minimal-component length.  The
refined closed union $E_i^*$ differs from the coarse source union
$E_i^{\rm src}$ of
\eqref{B:eq:circle-active-unions}, but satisfies
$E_i^*\subset\operatorname{cl}_{\T}E_i^{\rm src}$.
\end{definition}

The covers are additional data attached to the sparse selection and
do not alter the source--target components of
\cite[Definition~5.1]{LWY}.
By the support localization following
\eqref{B:eq:ham-moving-front}
and \cref{B:lem:nested-cotangent,B:prop:ledger},
\begin{equation}\label{B:eq:physical-support-cover}
 \supp\Delta_i\subset\Omega_i\times[0,1].
\end{equation}
We first bound the lengths of the covers, then show that each new
front still has enough candidates avoiding all the earlier ones.

\begin{lemma}\label{B:lem:support-size}
The $2D_i$ closed arcs in \eqref{B:eq:Omega-i} satisfy
\begin{equation}\label{B:eq:Omega-size}
 |\Omega_i|\le \mathrm{Len}_i\le C_\Omega\frac{2^{-i}}{(i+1)^2}
\end{equation}
uniformly over all admissible sparse selections at sufficiently late stages.
\end{lemma}
\begin{proof}
	For each $(r,j)\in\mathcal A_i$,
	the closed source arc
	$R_\alpha^j\overline{J_r}$
	is disjoint from $\mathcal P_i$,
	by \cref{B:prop:protected-density}
	and padded separation.
	Thus \eqref{B:eq:rho-lower}
	and the change-of-variables
	formula give
	\[
	|K_{i,r,j}|
	=\int_{R_\alpha^j\overline{J_r}}
	\rho_i(y)^{-1}\,dy
	\le |J_r|.
	\]

	By \eqref{B:eq:micro-admissibility},
	$A_r\delta_r$ is decreasing
	for $r\ge K_{\min}$.
	Consequently, for $r\ge a_i$,
	\[
	|J_r|
	\le |J_{a_i}|
	=2(Ca_i^2+3)2^{-a_i}
	\le C_1a_i^22^{-a_i}.
	\]
	Since there are $d_{i,r}$
	selected edges at scale $r$
	on each of the two fronts,
	\[
\mathrm{Len}_i\le2\sum_{r=a_i}^{b_i-1}d_{i,r}|J_r|\le C_2D_i a_i^22^{-a_i}\le C_2\vartheta q_i a_i^22^{-a_i}\le C_\Omega\frac{2^{-i}}{(i+1)^2}.
\]
	where the final inequality
	follows from
	\eqref{B:eq:checkpoint-support}.
	The constants are independent
	of the sparse selection.
	Finally,
	$|\Omega_i|\le\mathrm{Len}_i$
	by subadditivity of arc length.
\end{proof}

Summing over the stages, we obtain
\begin{equation}\label{B:eq:support-tail-sum}
 \sum_{i\ge n_0}|\Omega_i|\le\sum_{i\ge n_0}\mathrm{Len}_i
       \le C_\Omega\sum_{i\ge n_0}\frac{2^{-i}}{(i+1)^2},
\end{equation}
and the right-hand side tends to zero as $n_0\to\infty$.

\Needspace{5\baselineskip}

\begin{lemma}\label{B:lem:packing}
Let $\delta>0$ and let $X\subset\T$ be $\delta$-separated, and let
$E\subset\bigcup_{\ell=1}^m I_\ell$, where the $I_\ell$ are closed arcs
which may overlap.  Put $S_E=\sum_\ell|I_\ell|$.
Then
\begin{equation}\label{B:eq:packing}
 \#(X\cap E)\le\delta^{-1}S_E+m.
\end{equation}
For $a\ge0$, let
\[
E^{+a}:=\{x\in\T:
\operatorname{dist}_{\T}(x,E)\le a\}
\]
denote the closed $a$-neighborhood
in the circle metric.
\begin{equation}\label{B:eq:thickening}
 E^{+a}\subset\bigcup_{\ell=1}^m I_\ell^{+a},\qquad
 \sum_{\ell=1}^m|I_\ell^{+a}|\le S_E+2ma.
\end{equation}
A thickened arc that fills the circle is assigned length $1$.
\end{lemma}
\begin{proof}
An arc of length $L$ contains at most $L/\delta+1$ separated points,
and summing over the arcs of the cover gives the first claim;
overcounting is harmless.  Each thickened arc has length at most $L+2a$.
\end{proof}

\Needspace{5\baselineskip}
Combining this packing bound with the small total length of the earlier supports leaves enough unused candidates on both fronts.

\begin{proposition}\label{B:prop:avoidance}
After increasing $n_0$, each of the two candidate fronts admits exactly
$D_n$ selected edges, and the resulting closed sets $\Omega_n$ are
pairwise disjoint.  All preceding estimates remain valid.
\end{proposition}

\begin{proof}
We select the edges inductively, starting at $n=n_0$ with the empty
old cover.
Suppose that the stages $n_0,\ldots,n-1$ have been selected, with
disjoint physical angular covers, and put
\[
 \Omega_{<n}:=\bigcup_{i=n_0}^{n-1}\Omega_i,\qquad
 B_n:=h_n(\Omega_{<n})
     =\bigcup_{i=n_0}^{n-1}\ \bigcup_{(r,j)\in\mathcal A_i}h_n(K_{i,r,j}).
\]
We call $B_n$ the \emph{forbidden set in the current rotation
coordinates}: it is the old physical covers transported into the
coordinate system in which the new candidate intervals live.  A
candidate edge $j\in\mathcal J_n^\pm$ is \emph{good} when
\[
 R_\alpha^j\overline{J_{a_n}}\cap B_n=\varnothing,
\]
and \emph{bad} otherwise.  Since
$\overline{J_r}\subset\overline{J_{a_n}}$ for $r\ge a_n$, every finer
source interval assigned to a good edge avoids $B_n$ as well.
Selecting $D_n$ good edges on each front therefore yields a new
physical cover disjoint from $\Omega_{<n}$.  To see that enough good
candidates remain, we bound the number of bad ones.

\Needspace{15\baselineskip}
The change of coordinates is illustrated in \Cref{B:fig:physical-avoidance}.

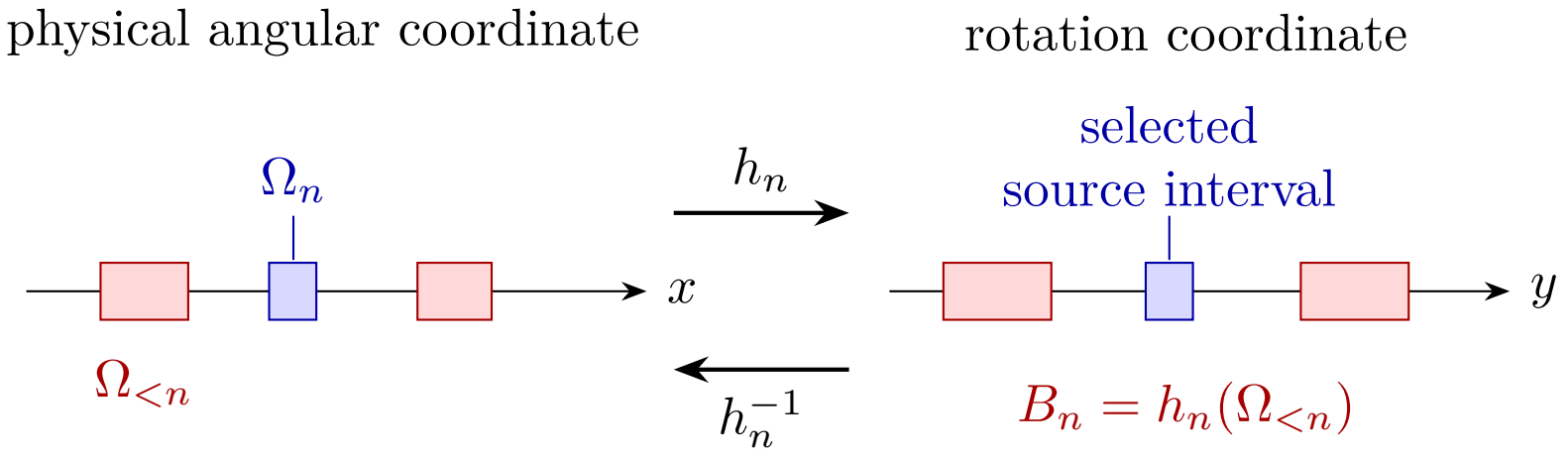
\begin{figure}[H]
\centering
\begin{tikzpicture}[x=0.90cm,y=0.95cm,>=Stealth,font=\small]
 \draw[->] (0,0)--(4.6,0) node[right]{$x$};
 \draw[->] (6.4,0)--(11,0) node[right]{$y$};
 \node[align=center] at (2.2,1.82) {physical angular coordinate};
 \node[align=center] at (8.6,1.82) {rotation coordinate};
 \filldraw[fill=red!15,draw=red!65!black] (0.55,-0.2) rectangle (1.2,0.2);
 \filldraw[fill=red!15,draw=red!65!black] (2.9,-0.2) rectangle (3.45,0.2);
 \filldraw[fill=blue!15,draw=blue!65!black] (1.8,-0.2) rectangle (2.15,0.2);
 \node[red!65!black] at (0.87,-0.65) {$\Omega_{<n}$};
 \node[blue!65!black] at (1.98,0.78) {$\Omega_n$};
 \draw[blue!65!black,thin] (1.98,0.53)--(1.98,0.22);
 \filldraw[fill=red!15,draw=red!65!black] (6.8,-0.2) rectangle (7.6,0.2);
 \filldraw[fill=red!15,draw=red!65!black] (9.45,-0.2) rectangle (10.25,0.2);
 \filldraw[fill=blue!15,draw=blue!65!black] (8.3,-0.2) rectangle (8.65,0.2);
 \node[red!65!black] at (8.6,-0.82) {$B_n=h_n(\Omega_{<n})$};
 \node[blue!65!black,align=center] at (8.475,0.95) {selected\\source interval};
 \draw[blue!65!black,thin] (8.475,0.53)--(8.475,0.22);
 \draw[->,thick] (4.8,0.55)--(6.1,0.55) node[midway,above]{$h_n$};
 \draw[->,thick] (6.1,-0.55)--(4.8,-0.55) node[midway,below]{$h_n^{-1}$};
\end{tikzpicture}
\caption{Physical and rotation angular coordinates in recursive
support selection.  The current conjugacy $h_n$ sends the old physical
cover $\Omega_{<n}$ to the forbidden set $B_n$.
A good candidate avoids $B_n$ with its entire coarsest source interval;
every finer interval assigned to it pulls back disjointly from
$\Omega_{<n}$.  The arcs are schematic and are not drawn to scale.}
\label{B:fig:physical-avoidance}
\end{figure}

We work with the displayed cover of the forbidden set, indexed by the
original arcs.  Its number of arcs and summed length satisfy
\begin{equation}\label{B:eq:component-count}
 m_n=2\sum_{i=n_0}^{n-1}D_i,
\end{equation}
\begin{equation}\label{B:eq:B-small}
 |B_n|\le \mathrm{Len}_n^{\rm old}:=
 \sum_{i=n_0}^{n-1}\sum_{(r,j)\in\mathcal A_i}|h_n(K_{i,r,j})|
 \le B_{\rm der}C_\Omega\sum_{i\ge n_0}\frac{2^{-i}}{(i+1)^2}.
\end{equation}
Neither a minimal component count nor disjointness of the transported
cover is needed here.  By \eqref{B:eq:sum-q},
\begin{equation}\label{B:eq:component-ratio}
 m_n\le8\vartheta q_{n-1}\le8\vartheta q_n.
\end{equation}
With $\ell_n=A_{a_n}\delta_{a_n}$,
the quantitative checkpoint bound
used in the proof of
\cref{B:lem:padded-separation}
gives
\begin{equation}\label{B:eq:ell-small}
	\ell_n<(4q_n)^{-1}.
\end{equation}
A candidate edge $j$ is bad if and only if
\[
 R_\alpha^j\overline{J_{a_n}}\cap B_n\ne\varnothing
 \quad\Longleftrightarrow\quad
 \theta_*+j\alpha\in B_n^{+\ell_n}.
\]
On either front the centers are $(2q_n)^{-1}$-separated, since all
nonzero index differences have absolute value $<q_n$ while
$\|q_{n-1}\alpha\|_\T>(q_{n-1}+q_n)^{-1}>(2q_n)^{-1}$.
Applying \cref{B:lem:packing} to the explicitly thickened cover, we
find
\begin{align*}
 \#\{\text{bad candidates on one front}\}
 &\le2q_n(\mathrm{Len}_n^{\rm old}+2m_n\ell_n)+m_n\\
 &\le2q_n\mathrm{Len}_n^{\rm old}+2m_n.
\end{align*}
Now increase $n_0$ so that $\mathrm{Len}_n^{\rm old}<\mathfrak c/8$ for all $n\ge n_0$;
with $\mathfrak c=1/100$ and $\vartheta=1/10000$ we then have
\begin{equation}\label{B:eq:good-margin}
 2\mathrm{Len}_n^{\rm old}+2m_n/q_n
 <\mathfrak c/4+16\vartheta=0.0041
 <(\mathfrak c-\vartheta)/2.
\end{equation}
As $L_n\ge\mathfrak c q_n-1$ and $D_n\le\vartheta q_n$, this gives
\[
 L_n-\#\{\text{bad}\}-D_n
 >\tfrac12(\mathfrak c-\vartheta)q_n-1>0
\]
at all sufficiently late stages.  We select $D_n$ good edges on each
front and allocate them by \eqref{B:eq:balanced-allocation}.

Since $\overline{J_r}\subseteq\overline{J_{a_n}}$ for each selected
scale $r\ge a_n$,
\[
 E_n^*\cap h_n(\Omega_{<n})=\varnothing
 \quad\Longrightarrow\quad\Omega_n\cap\Omega_{<n}=\varnothing.
\]
All the sets involved are closed.  Finally, the dynamical estimates
proved earlier involve the allocated counts but not the positions of
the selected edges, so they remain valid under this induction.  See
\cref{B:fig:physical-avoidance}.
\end{proof}

\section{\texorpdfstring{$C^1$}{C1} convergence and disk completion}\label{B:sec:completion}
In this section, we first prove that the sequence of annulus maps converges in the \(C^1\) topology. We then show that the limit has the required properties and extend it to the whole disk, following the construction in \cite[\S7]{LWY}.

\subsection{Disjoint-support \texorpdfstring{$C^1$}{C1} convergence}

\Needspace{5\baselineskip}
The main point is that the supports of different increments are pairwise disjoint. Therefore, at each point, at most one increment contributes to the difference between two maps. As a result, the \(C^1\) norm of the tail is controlled by the largest remaining increment, and we do not need to assume that the increments are summable.

\begin{lemma}
\label{B:lem:localized-Cauchy}
On $E=\R\times[0,1]$ let $f_n:E\to\R^2$ be $C^1$ degree-one
lifts, so $f_n(x+1,s)=f_n(x,s)+(1,0)$.  Put $u_n=f_{n+1}-f_n$.
Norms are suprema over $[0,1]\times[0,1]$, with the sum convention for $C^1$
fixed in \cref{sec:conventions}.
If the relative closed supports of the periodic functions $u_n$ are
pairwise disjoint and $\|u_n\|_{C^1}\to0$, then $f_n$ is Cauchy in
this norm and
\begin{equation}\label{B:eq:tail-sup}
 \|f_N-f_M\|_{C^1}\le 2\max_{M\le n<N}\|u_n\|_{C^1}.
\end{equation}
\end{lemma}
\begin{proof}
Write $f_N-f_M=\sum_{n=M}^{N-1}u_n$ in the fixed normed coordinate
space.  Since the supports are pairwise disjoint, at most one summand
is nonzero at any given point.  Off $\supp_Eu_n$ the function $u_n$
vanishes in a relative neighborhood of each point (one-sided at
$s=0,1$), so $Du_n=0$ there as well.  We therefore have the two
estimates
\[
 \|f_N-f_M\|_0\le\max_{M\le n<N}\|u_n\|_0,\qquad
 \|D(f_N-f_M)\|_0\le\max_{M\le n<N}\|Du_n\|_0.
\]
Each maximum is at most
$\max_{M\le n<N}\|u_n\|_{C^1}$,
so adding the two inequalities
gives \eqref{B:eq:tail-sup}.
Since
\[
\sup_{n\ge M}\|u_n\|_{C^1}
\longrightarrow0
\qquad(M\to\infty),
\]
the estimate proves that $(f_n)$
is Cauchy in $C^1$, uniformly
for $N>M$.
Completeness on the fixed
fundamental domain gives a
$C^1$ limit, and the degree-one
relation passes to this limit.
\end{proof}

\Needspace{5\baselineskip}
Applying this criterion to the stage increments gives the annulus limit and preserves the stabilized wandering intervals on its outer boundary.

\begin{proposition}\label{B:prop:annulus-limit}
The maps $F_n$ converge in $C^1$ to an area-preserving $C^1$ diffeomorphism $F:\mathbb A\to\mathbb A$.  Its outer boundary restriction $g$ is a $C^1$ Denjoy diffeomorphism of rotation number $\alpha$.
\end{proposition}

\begin{proof}
By \cref{B:prop:avoidance} and \eqref{B:eq:physical-support-cover}
the supports $\supp\Delta_n\subset\Omega_n\times[0,1]$ are pairwise
disjoint; together with $\|\Delta_n\|_{C^1}\to0$, the localized Cauchy
criterion gives $F_n\to F$ in $C^1$ in compatible annulus lifts.

From here we follow the local-to-global
argument of
\cite[Proposition~7.1]{LWY}.
Since every $F_n$ preserves area
and $F_n\to F$ in $C^1$,
\[
\det DF=\lim_n\det DF_n=1
\]
throughout the closed annulus.
In particular, $DF$ is invertible
at every point.

The limit maps $\mathbb A$ into itself
and maps each boundary circle into itself.
No interior point can be mapped to
the boundary: otherwise the normal
component of $F$ would attain
an interior local extremum, forcing
its differential to vanish and
contradicting $\det DF=1$.

At either boundary circle, the
derivative in local product
coordinates has the form
\[
DF=
\begin{pmatrix}
	\partial_xF_1&\partial_sF_1\\
	0&\partial_sF_2
\end{pmatrix}.
\]
The tangential derivative is
nonnegative as the $C^1$ limit
of those of orientation-preserving
boundary diffeomorphisms.
Since $\det DF=1$, both
$\partial_xF_1$ and
$\partial_sF_2$ are strictly
positive.
Consequently, $F$ is a local
$C^1$ diffeomorphism also in
boundary charts.

It follows that $F(\mathbb A)$
is relatively open in $\mathbb A$.
It is also closed, by compactness.
As $\mathbb A$ is connected,
$F(\mathbb A)=\mathbb A$.
Thus $F$ is a proper surjective
local diffeomorphism and hence a
finite-sheeted covering.

Each boundary restriction has
degree one, inherited from
the corresponding restrictions
of the maps $F_n$.
Since these restrictions are
orientation-preserving local
diffeomorphisms, each has
exactly one sheet.
The covering degree of $F$
is therefore one.
Hence $F$ is a global
area-preserving $C^1$
diffeomorphism of $\mathbb A$.

On the outer boundary, $g_n\to g=F|_{s=0}$ in $C^1$, and continuity
of the rotation number gives $\rot(g)=\alpha$.  For each $j$ the
stabilized intervals of \cref{B:lem:stabilized} eventually satisfy
$g_n(I_j)=I_{j+1}$, whence $g(I_j)=I_{j+1}$.  Since
their interiors are nonempty and pairwise disjoint, these interiors
form a wandering orbit, and $g$ is a Denjoy diffeomorphism.
\end{proof}

\Needspace{5\baselineskip}

\begin{proposition}\label{B:prop:sharp-increment}
For some constants $c,C>0$ and all sufficiently large $n$,
\begin{equation}\label{B:eq:Delta-sharp}
 c e_n\le\|F_{n+1}-F_n\|_{C^1}\le C e_n.
\end{equation}
In particular the increment norms are not summable; the lower bound
plays no part in the disjoint-support convergence criterion.
The boundary limit $g$ is not $C^{1,\beta}$ for any $\beta>0$.
\end{proposition}

\begin{proof}
We first draw a consequence of support avoidance.  We have $F_{n_0}=\mathcal R_\alpha$ and
\[
 F_n-\mathcal R_\alpha=\sum_{i=n_0}^{n-1}\Delta_i,
 \qquad \supp\Delta_i\subset\Omega_i\times[0,1],
\]
where the finitely many closed sets $\Omega_i$, $i<n$, are disjoint from the compact set $\Omega_n$.  It follows that $F_n$ coincides with $\mathcal R_\alpha$ on an angular neighborhood of $\Omega_n\times[0,1]$; in particular $g_n=R_\alpha$ and $g_n'=1$ near every current physical source interval.

Now take a selected edge $j$ on the positive front, assigned to the highest micro-scale
\[
 r_n:=b_n-1,
 \qquad d_n:=d_{n,r_n},
 \qquad u_n:=d_n^{-1}>0.
\]
Such an edge exists because $d_n\ge1$.  Across this positive edge the parameters decrease, so the moving-front quotient is $L_{r_n,u_n}$.  Put
\[
 y_n=\theta_*+j\alpha,\qquad x_n=h_n^{-1}(y_n),\qquad
 z_n=h_n^{-1}(y_n+\alpha).
\]
A micro-flow fixes its angular center, where its derivative is
\[
 L_{r_n,u_n}'(\theta_*)=
 \exp\bigl(u_nV_{\varepsilon_{r_n}}'(0)\bigr)=e^{-u_n}.
\]
Combining $g_n=h_n^{-1}R_\alpha h_n$ with $g_n'=1$ at $x_n$ gives
\[
 \frac{h_n'(x_n)}{h_n'(z_n)}=1.
\]
At the center, the moving-front
identity gives
\[
S_n(y_n)=y_n+\alpha,
\qquad S_n'(y_n)=e^{-u_n}.
\]
Consequently, the conjugacy formula gives
\begin{equation}\label{B:eq:exact-center}
 g_{n+1}(x_n)=g_n(x_n)=R_\alpha(x_n),
 \qquad g_{n+1}'(x_n)=e^{-u_n}.
\end{equation}
As \eqref{B:eq:d-asymp} gives
$u_n\ge e_n/(2\vartheta)$ and $u_n\le1$,
\[
 |g_{n+1}'(x_n)-g_n'(x_n)|
 =1-e^{-u_n}\ge \tfrac12u_n
 \ge\frac{e_n}{4\vartheta}.
\]
This is a tangential derivative of the annulus increment, so it
proves the lower bound in \eqref{B:eq:Delta-sharp}; the upper bound is
\eqref{B:eq:Delta-C1-vanish}.  Together with the divergence
$\sum_ne_n=\infty$ of \cref{B:lem:width}, the lower bound shows that
the increment norms are not summable.

For the H\"older statement, let
\[
y_n^+=y_n+A_{r_n}\delta_{r_n},
\qquad
x_n^+=h_n^{-1}(y_n^+),
\]
where $y_n^+$ is the right endpoint
of the translated closed scale-$r_n$
interval $R_\alpha^j\overline{J_{r_n}}$,
under the compatible lift.

Since $V_{\varepsilon_{r_n}}$
vanishes near the endpoints of
$(-A_{r_n},A_{r_n})$, the local
flow is the identity near $y_n^+$.
Consequently,
$g_{n+1}'(x_n^+)=g_n'(x_n^+)=1$.
Both $x_n$ and $x_n^+$ lie in the closed set $\Omega_n$.  By support avoidance every later increment vanishes there together with its first derivative, so passing to the $C^1$ limit gives
\begin{equation}\label{B:eq:limit-oscillation}
 g'(x_n)=e^{-u_n},\qquad g'(x_n^+)=1.
\end{equation}
Moreover $\rho_n\ge1$ throughout the selected source interval, so the positively oriented arc from $x_n$ to $x_n^+$ has length
\begin{equation}\label{B:eq:physical-micro-length}
 \int_{y_n}^{y_n^+}\frac{dy}{\rho_n(y)}
 \le A_{r_n}\delta_{r_n}\le C r_n^22^{-r_n}.
\end{equation}
For large $n$ this arc is shorter than $1/2$, so the same quantity
bounds the circle distance between the two points.
Because $\sum_ne_n=\infty$, infinitely many indices $n$ satisfy
$e_n\ge(n+1)^{-2}$; otherwise the series would converge.  Along
these indices, for any fixed $\beta>0$,
\eqref{B:eq:limit-oscillation}--\eqref{B:eq:physical-micro-length} give
\[
 \frac{|g'(x_n)-g'(x_n^+)|}
      {\operatorname{dist}_{\T}(x_n,x_n^+)^\beta}
 \ge\frac{c e_n}{(Cr_n^22^{-r_n})^\beta}
 \ge \frac{c_\beta 2^{\beta r_n}}
 {(n+1)^2r_n^{2\beta}}
 \longrightarrow\infty.
\]
For the final step, \eqref{B:eq:checkpoints} gives
$r_n=k_{n+1}-1\ge n+1$ after increasing $n_0$, so the last expression is
bounded below by $c_\beta2^{\beta r_n}/r_n^{2+2\beta}$, which tends to
infinity.  Therefore $g'$ has no positive H\"older exponent.
\end{proof}

\subsection{Completion to the disk}\label{B:sec:disk-completion}

We now extend the annulus map to the whole disk. We first show that the conjugacies and their inverses stabilize on every compact subset away from the outer boundary. This gives a smooth conjugacy to the rigid rotation in the punctured interior. We then attach a rotating disk along the inner collar, following the construction in \cite[\S7]{LWY}.

\begin{proposition}\label{B:prop:punctured}
On $\T\times(0,1]$ there is a smooth area-preserving diffeomorphism $\mathcal H_\infty$ such that
\begin{equation}\label{B:eq:punctured-conj}
 F=\mathcal H_\infty^{-1}\mathcal R_\alpha\mathcal H_\infty.
\end{equation}
Moreover $F$ is exactly the rigid rotation on a full inner collar.
\end{proposition}

\begin{proof}
As in \cite[\S7, proof of Proposition~7.1]{LWY}, it suffices to track
where the conjugacies and their inverses change from one stage to the
next.  The support of $\mathcal M_n$ lies in the rotation
collar $v\le\delta_{a_n}$, and by \eqref{B:eq:H-cotangent} and
\eqref{B:eq:Bder},
\begin{align*}
 \{z:\mathcal H_{n+1}(z)\ne\mathcal H_n(z)\}
   &\subset\T\times[0,B_{\rm der}\delta_{a_n}],\\
 \{z:\mathcal H_{n+1}^{-1}(z)\ne\mathcal H_n^{-1}(z)\}
   &\subset\T\times[0,\delta_{a_n}].
\end{align*}
The second inclusion uses
$\mathcal H_{n+1}^{-1}=\mathcal H_n^{-1}\mathcal M_n^{-1}$.
Since $a_n\to\infty$, both
$\delta_{a_n}$ and
$B_{\rm der}\delta_{a_n}$
tend to zero.

For any compact
$K_0\subset\T\times(0,1]$,
the radial coordinate has
a positive minimum on $K_0$.
Hence, for all sufficiently
large $n$, both change regions
displayed above are disjoint
from $K_0$.
Thus $\mathcal H_n$ and
$\mathcal H_n^{-1}$ eventually
stabilize on every such compact
set.

Denote their limits by
$\mathcal H_\infty$ and
$\mathcal K_\infty$.
Locally each limit agrees with
a finite-stage smooth
area-preserving diffeomorphism.
In particular, both limits are
smooth, preserve area, and take
values in $\T\times(0,1]$.

Let $y=\mathcal H_\infty(x)$ and take $N$ beyond the stabilization
thresholds at $x$ and at $y$; then
\[
 \mathcal K_\infty(\mathcal H_\infty(x))
 =\mathcal H_N^{-1}(\mathcal H_N(x))=x.
\]
The identity in the other order is obtained the same way, so
$\mathcal K_\infty=\mathcal H_\infty^{-1}$.
Passing to the limit in the conjugacy identity, after stabilizing
also at $\mathcal R_\alpha\mathcal H_\infty(x)$, yields
\eqref{B:eq:punctured-conj}.

Finally, every micro-factor is the identity on the fixed inner collar
$s>\delta_{a_{n_0}}$, so $\mathcal H_\infty=\id$ and
$F=\mathcal R_\alpha$ there.
\end{proof}

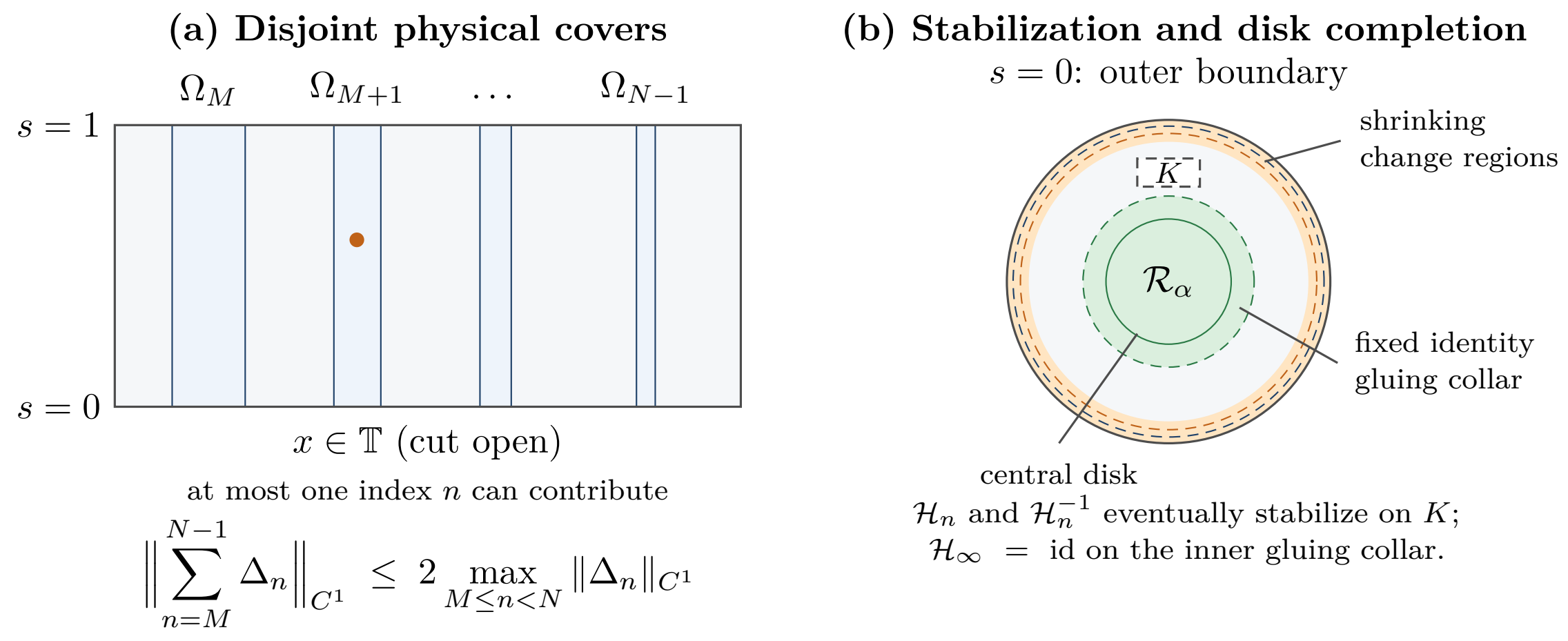
\begin{figure}[htbp]
\centering
\begin{tikzpicture}[x=1cm,y=1cm,font=\small]
 \node[font=\bfseries\small] at (2.9,3.6)
   {(a) Disjoint physical covers};
 \fill[softgray] (0,0) rectangle (6.0,2.7);
 \foreach \a/\b in {.55/1.25,2.1/2.55,3.5/3.8,5.0/5.18}{
   \filldraw[fill=softblue,draw=deepblue] (\a,0) rectangle (\b,2.7);
 }
 \draw[nipl] (0,0) rectangle (6.0,2.7);
 \node[left] at (0,0) {$s=0$};
 \node[left] at (0,2.7) {$s=1$};
 \node[above] at (.9,2.75) {$\Omega_M$};
 \node[above] at (2.325,2.75) {$\Omega_{M+1}$};
 \node[above] at (3.65,2.75) {$\cdots$};
 \node[above] at (5.09,2.75) {$\Omega_{N-1}$};
 \fill[frontline] (2.32,1.6) circle (2pt);
 \node[align=center,font=\scriptsize]
   at (3,-.79) {at most one index $n$ can contribute};
 \node at (3,-.35) {$x\in\T$ (cut open)};
 \node[align=center,text width=6.7cm] at (2.9,-1.60)
   {$\displaystyle
    \Bigl\|\sum_{n=M}^{N-1}\Delta_n\Bigr\|_{C^1}
       \le 2\max_{M\le n<N}\|\Delta_n\|_{C^1}$};
 \begin{scope}[xshift=7.6cm]
  \node[font=\bfseries\small] at (2.65,3.6)
    {(b) Stabilization and disk completion};
  \fill[softgray] (2.5,1.2) circle (1.55);
  \fill[frontorange] (2.5,1.2) circle (1.55);
  \fill[softgray] (2.5,1.2) circle (1.34);
  \draw[frontline,densely dashed] (2.5,1.2) circle (1.42);
  \draw[deepblue,densely dashed] (2.5,1.2) circle (1.49);
  \fill[coregreen] (2.5,1.2) circle (.82);
  \draw[coregreenline,densely dashed] (2.5,1.2) circle (.82);
  \filldraw[fill=coregreen,draw=coregreenline] (2.5,1.2) circle (.60);
  \draw[nipl] (2.5,1.2) circle (1.55);
  \node at (2.5,1.2) {$\mathcal R_\alpha$};
  \draw[nipl] (2.2,.70)--(1.45,-.35);
  \node[anchor=north,font=\scriptsize] at (1.45,-.41) {central disk};
  \draw[nipl,densely dashed,fill=white] (2.2,2.11) rectangle (2.8,2.38);
  \node[font=\scriptsize] at (2.5,2.245) {$K$};
  \node[font=\scriptsize,align=left,anchor=west] at (4.2,2.55)
    {shrinking\\change regions};
  \draw[nipl] (3.42,2.35)--(4.15,2.55);
  \node[above] at (2.5,2.9) {$s=0$: outer boundary};
  \draw[nipl] (3.18,.95)--(4.12,.42);
  \node[anchor=west,font=\scriptsize,align=left] at (4.15,.42)
    {fixed identity\\gluing collar};
  \node[align=center,text width=6.5cm,font=\scriptsize] at (2.65,-1.18)
    {$\mathcal H_n$ and $\mathcal H_n^{-1}$ eventually stabilize on $K$;\\
     $\mathcal H_\infty=\id$ on the inner gluing collar.};
 \end{scope}
\end{tikzpicture}
\caption{Two localization steps in \S\ref{B:sec:completion}.
Left: the columns are disjoint covers $\Omega_n\times[0,1]$, not
prescribed shapes for $\supp\Delta_n$; the largest increment controls
the $C^1$ tail.  Right: the change regions retreat to the outer
boundary, so the conjugacies and their inverses stabilize on each
fixed compact set $K$ away from that boundary.  The common inner
identity collar allows attachment of a rotating central disk.
No extension of $\mathcal H_\infty$ to the outer boundary is asserted.
All geometry is schematic.}
\label{B:fig:localization-completion}
\end{figure}

The two localization steps are summarized in
\cref{B:fig:localization-completion}.

The annulus is completed to a disk in the fixed action coordinates of
\cite[\S7]{LWY}.  Fix $I_*>0$ once and attach a central disk carrying
the same rotation; put $I=I_*+1-s$ on the annulus and
$I=\pi r^2$ on the central disk, with $x=\theta/(2\pi)$.  Along the
gluing collar, on each piece the area form and the map are
\[
 dI\wedge dx,\qquad (x,I)\longmapsto(x+\alpha,I).
\]
The pieces therefore glue smoothly and area is preserved.  A fixed
radial identification and homothety produce the standard unit disk;
their constant area normalization cancels under conjugation.  These
identifications are independent of $n_0$ and, within the bounded
$C^1$ neighborhood used here, bound the disk $C^1$ perturbation norm by a
fixed multiple of the annulus norm.

Because $\mathcal H_\infty=\id$ on the gluing collar, it extends by the
identity across the central disk.  In the final disk coordinates this
gives a single smooth area-preserving conjugacy $\Psi$ on
$\operatorname{int}\D\setminus\{0\}$ with
$F=\Psi^{-1}\mathcal R_\alpha\Psi$ there.

\begin{lemma}\label{B:lem:pseudo}
The disk map $F$ is a $C^1$ area-preserving irrational pseudo-rotation with
\[
 \Per(F)=\{0\},
 \qquad
 \rot(F)=\alpha.
\]
On $\operatorname{int}(\D)\setminus\{0\}$ it is smoothly conjugate to the rigid rotation.
\end{lemma}

\begin{proof}
Away from the center the punctured conjugacy rules out periodic
points, while the boundary map has irrational rotation number $\alpha$;
the center itself is fixed since the map is a rotation near it.  Hence
$\Per(F)=\{0\}$ and $\rot(F)=\alpha$, and the regularity and area
preservation asserted follow from the fixed completion just described.
\end{proof}

\Needspace{6\baselineskip}
The wandering boundary interval prevents uniform return on the whole disk.

\begin{lemma}\label{B:lem:nonrigid}
The map $F$ is not $C^0$-rigid.
\end{lemma}

\begin{proof}
Let $J=\operatorname{int}I_0$ be a wandering boundary interval and pick
$x\in J$.  Since $g^m(x)\notin J$ for every $m\ge1$,
\[
 d_\T(g^m(x),x)\ge d_\T(x,\T\setminus J)>0.
\]
Hence no sequence of positive iterates of $g$ converges uniformly
to the identity.  Such convergence for $F$ would restrict to the
boundary, and is therefore impossible; this is the boundary
obstruction of \cite[\S7]{LWY}.
\end{proof}

\begin{proof}[Proof of \cref{thm:nonrigidity,rem:near-rotation}]
The arithmetic of the prescribed $\alpha$ enters through
\cref{B:lem:arithmetic,B:lem:width}.  The recursive support selection
of \cref{B:prop:avoidance} and the stage increment bounds then give
the $C^1$ annulus limit of \cref{B:prop:annulus-limit}.  The disk
completion and the pseudo-rotation properties are
\cref{B:prop:punctured,B:lem:pseudo}, and non-rigidity is
\cref{B:lem:nonrigid}.  By \cref{B:prop:sharp-increment}, the derivative
of the boundary map has no positive H\"older exponent.  The outer annulus
boundary is identified with the boundary of the round disk by a fixed
smooth parametrization, so this persists after disk completion.  A
$C^{1,\beta}$ disk map would have a $C^{1,\beta}$ boundary
restriction, which is impossible.

For the additional proximity assertion in \cref{rem:near-rotation},
start the same construction at a sufficiently late index $n_0$.
All constants in the increment estimate and in the fixed disk
identification are independent of $n_0$, and the localized
Cauchy criterion \eqref{B:eq:tail-sup}, with
$F_{n_0}=\mathcal R_\alpha$, gives
\[
 \|F-\mathcal R_\alpha\|_{C^1(\mathbb A)}
 \le C\sup_{n\ge n_0}e_n\longrightarrow0
 \qquad(n_0\to\infty).
\]
The fixed disk identification described above transfers this estimate
to the disk, giving the prescribed $\varepsilon$.
All the stated properties persist:
the construction and stabilization arguments apply for every
sufficiently large $n_0$, and removing a finite initial segment
preserves $\sum_{n\ge n_0}e_n=\infty$, used to exclude $C^{1,\beta}$
regularity for every $\beta>0$.
\end{proof}

\FloatBarrier
\section{Boundary rigidity and arithmetic obstruction}\label{B:sec:scope}
In this section, we study the rigidity on the boundary for a pseudo rotation and discuss some results for the whole disk. We will first prove  \cref{thm:boundary-rigidity} by estimating the derivatives of iterates near the boundary. And then, we consider the class \(\mathcal I_\alpha\) and show that, under the additional assumptions defining this class, the boundary estimates can be used to obtain rigidity on the whole disk. Finally, we show that the Denjoy construction can also be realized on a smooth invariant circle in the interior of the disk.

Throughout,
\[
 \|DF^m\|_\infty
 :=\sup_{z\in\D}\|D(F^m)_z\|
 =\sup_{z\in\D}\sup_{|\xi|=1}|D(F^m)_z\xi|
 =\|D(F^m)\|_{0;\D}.
\]
For derivatives of disk maps, the infinity subscript denotes the
supremum over the whole disk of the Euclidean operator norm;
suprema restricted to the boundary are written explicitly.
\Needspace{5\baselineskip}

\begin{lemma}\label{B:lem:boundary-growth}
Let $g$ be an orientation-preserving $C^1$ circle diffeomorphism with
irrational rotation number.  Then
\[
 \frac1m\log(g^m)'\longrightarrow0
 \quad\text{uniformly on }\T.
\]
If $F$ is a $C^1$ area-preserving disk diffeomorphism whose boundary
restriction is $g$, then
\[
 \lim_{m\to\infty}
 \frac1m\log\sup_{z\in\partial\D}\|DF^m(z)\|=0.
\]
In fact, the same bound holds simultaneously in both time directions:
\begin{equation}\label{B:eq:boundary-subexp}
 \lim_{m\to\infty}
 \frac1m\log\left(
 \max_{\sigma\in\{-1,1\}}\sup_{z\in\partial\D}
 \|DF^{\sigma m}(z)\|\right)=0.
\end{equation}
\end{lemma}

\begin{proof}
Since an irrational circle homeomorphism is uniquely ergodic, the
Birkhoff averages of the continuous function $\log g'$ converge
uniformly to a constant, so that $m^{-1}\log(g^m)'\to c$ uniformly.
As $\int_{\T}(g^m)'(x)\,dx=1$, neither $c>0$ nor $c<0$ is
possible, so $c=0$; see also \cite[Corollary~25]{Chenciner}.  For
every $\varepsilon>0$ there is $C_\varepsilon\ge1$ such that
\begin{equation}\label{B:eq:boundary-1d}
 C_\varepsilon^{-1}e^{-\varepsilon m}
 \le(g^m)'(x)\le C_\varepsilon e^{\varepsilon m}
 \qquad(x\in\T,\ m\ge0).
\end{equation}

Now express the derivative at a boundary point in the frames formed
by the unit tangent and the inward normal, at the source point and at
its image.  Tangent-line invariance and $\det DF=1$ give
\[
 DF(x)=
 \begin{pmatrix}g'(x)&b(x)\\0&1/g'(x)\end{pmatrix},
\]
with $b$ continuous and bounded.
In the moving orthonormal frames,
write
\[
A_j=
\begin{pmatrix}
	a_j&b_j\\
	0&a_j^{-1}
\end{pmatrix},
\qquad
a_j=g'(g^jx),\quad b_j=b(g^jx).
\]
Then
\[
A_{m-1}\cdots A_0
=
\begin{pmatrix}
	(g^m)'(x)&B_m(x)\\
	0&((g^m)'(x))^{-1}
\end{pmatrix},
\]
where
\[
B_m(x)=
\sum_{i=0}^{m-1}
\left(\prod_{j=i+1}^{m-1}a_j\right)
b_i
\left(\prod_{j=0}^{i-1}a_j^{-1}\right).
\]
Empty products are understood as one.

By \eqref{B:eq:boundary-1d},
each summand has absolute value
at most
$C_\varepsilon^2\|b\|_\infty
e^{\varepsilon m}$.
Consequently,
\[
|B_m(x)|
\le mC_\varepsilon^2
\|b\|_\infty e^{\varepsilon m},
\]
and hence
\[
\sup_{z\in\partial\D}\|DF^m(z)\|
\le C'_\varepsilon(m+1)e^{\varepsilon m}.
\]
Since $\varepsilon>0$ is arbitrary,
this gives an upper exponential
growth rate of at most zero.
The determinant-one condition
also gives $\|DF^m(z)\|\ge1$,
so the growth rate is exactly zero.

For the backward direction,
a two-dimensional matrix of
determinant one and its inverse
have the same Euclidean operator
norm. Since $F^m$ maps
$\partial\D$ bijectively onto
itself,
\[
\sup_{z\in\partial\D}\|DF^{-m}(z)\|
=
\sup_{z\in\partial\D}\|DF^m(z)\|.
\]
This proves
\eqref{B:eq:boundary-subexp}
in both time directions.
\end{proof}

\Needspace{5\baselineskip}
A fixed-block argument extends the boundary estimate to exponentially thin orbit tubes for every prescribed finite return time.

\begin{lemma}\label{B:lem:boundary-tubes}
Let $F$ be a $C^1$ diffeomorphism of $\D$ satisfying
\eqref{B:eq:boundary-subexp}.  For every $s>0$ there are constants
$\delta_s,c_s>0$ and $C_s\ge1$ such that, for every integer $Q\ge1$,
$x\in\partial\D$, and $y\in\D$ with
\[
 |y-x|\le c_s e^{-sQ},
\]
one has, for $\sigma\in\{-1,1\}$ and $0\le j\le Q$,
\begin{equation}\label{B:eq:boundary-tube}
 \begin{aligned}
 |F^{\sigma j}(y)-F^{\sigma j}(x)|
   &\le C_s e^{sj}|y-x|,\\
 \dist(F^{\sigma j}(y),\partial\D)&\le\delta_s/2.
 \end{aligned}
\end{equation}
All distances and derivative norms here are Euclidean.
\end{lemma}

\begin{proof}
Fix $s>0$ and use \eqref{B:eq:boundary-subexp} to choose an integer
$k\ge1$ with
\[
 \max_{\sigma\in\{-1,1\}}\sup_{z\in\partial\D}
       \|DF^{\sigma k}(z)\|\le e^{sk/2}.
\]
Since $F^k$ and $F^{-k}$ are two fixed $C^1$ maps, the continuity of
their derivatives yields $0<\delta_s<1/4$ such that
\begin{equation}\label{B:eq:boundary-block}
 \|DF^{\sigma k}(z)\|\le e^{sk}
 \quad\text{if }\dist(z,\partial\D)<\delta_s,\quad \sigma\in\{-1,1\}.
\end{equation}
We then set
\[
 C_s=\max_{\substack{\sigma\in\{-1,1\}\\0\le r<k}}
                   \|DF^{\sigma r}\|_\infty\ge1,
 \qquad c_s=\frac{\delta_s}{2C_s}.
\]
These constants depend only on $F$ and $s$.

Now fix $Q,x,y,\sigma$ as in the statement and put $d_0=|y-x|$.
We prove by induction that, whenever $\ell k\le Q$,
\begin{equation}\label{B:eq:boundary-block-induction}
 |F^{\sigma\ell k}(y)-F^{\sigma\ell k}(x)|
       \le e^{s\ell k}d_0\le c_s.
\end{equation}
The case $\ell=0$ is immediate.  Suppose
\eqref{B:eq:boundary-block-induction} holds and $(\ell+1)k\le Q$.
By convexity, the straight segment
joining $F^{\sigma\ell k}(x)\in\partial\D$
to $F^{\sigma\ell k}(y)$ lies in $\D$.
Every point $z$ of this segment
satisfies
\[
\dist(z,\partial\D)
\le |z-F^{\sigma\ell k}(x)|
\le c_s<\delta_s.
\]
Thus \eqref{B:eq:boundary-block}
applies along the entire segment.
The mean-value inequality for $F^{\sigma k}$ along this segment,
combined with \eqref{B:eq:boundary-block}, gives the first inequality
in \eqref{B:eq:boundary-block-induction} at $\ell+1$; the second
follows from $d_0\le c_s e^{-sQ}$ and $(\ell+1)k\le Q$.

For a general $j=\ell k+r\le Q$ with $0\le r<k$, the global
mean-value bound for $F^{\sigma r}$ gives
\[
 \begin{split}
 |F^{\sigma j}(y)-F^{\sigma j}(x)|
 &\le C_s e^{s\ell k}d_0
 \le C_s e^{sj}d_0,\\
 |F^{\sigma j}(y)-F^{\sigma j}(x)|
 &\le C_s c_s=\delta_s/2.
 \end{split}
\]
Since $F^{\sigma j}(x)\in\partial\D$, these are the two assertions
of \eqref{B:eq:boundary-tube}.  The induction controls the orbit over
the prescribed time horizon; at no point is the collar itself assumed
to be invariant.
\end{proof}

\begin{proof}[Proof of \cref{thm:boundary-rigidity}]
Put $g=F|_{\partial\D}$ in angular coordinates.
We prove a quantitative version.  Suppose $Q_j\to\infty$ are positive
integers and, for some $\eta>0$,
\[
 \varepsilon_j:=\normT{Q_j\alpha}\le e^{-\eta Q_j}.
\]
For every $0<s<\eta/2$ there is $K_{F,s}>0$ such that, for all
sufficiently large $j$,
\begin{equation}\label{B:eq:boundary-rigidity-rate}
 \sup_{x\in\partial\D}|F^{Q_j}(x)-x|
 \le K_{F,s}e^{sQ_j}\sqrt{\varepsilon_j}
 \le K_{F,s}e^{-(\eta/2-s)Q_j}.
\end{equation}

Normalize area by $\lambda=\pi^{-1}dX\,dY$, so that $\lambda(\D)=1$.
We use the following form of \cite[Lemma~2.1]{AFLXZ}: if $H$ is a
continuous area-preserving disk pseudo-rotation with
$\rot(H)=\varepsilon+\Z$, $0<\varepsilon<1$, and
$U\subset\operatorname{int}\D\setminus\{0\}$ is an open topological
disk, then
\begin{equation}\label{B:eq:free-disk-capacity}
 H(U)\cap U=\varnothing\quad\Longrightarrow\quad\lambda(U)\le\varepsilon.
\end{equation}
This is the free-disk estimate normalized to total area one, and it
requires no differentiability.

Fix $0<s<\eta/2$.  By \cref{B:lem:boundary-growth},
\cref{B:lem:boundary-tubes} applies to $F$, and we write $c_s,C_s$
for the constants it provides.  As $\alpha$ is irrational,
$\varepsilon_j>0$; choose $\sigma_j\in\{-1,1\}$ so that
\[
 H_j:=F^{\sigma_jQ_j}
 \quad\text{satisfies}\quad
 \rot(H_j)=\varepsilon_j+\Z.
\]
Each $H_j$ preserves $\lambda$ and has the origin as its only
periodic point, so \eqref{B:eq:free-disk-capacity} applies to it.

Put $r_j=2\sqrt{\varepsilon_j}$.  For all sufficiently large $j$ we
have $r_j<1/4$ and
\begin{equation}\label{B:eq:boundary-radius-budget}
 3r_j e^{sQ_j}
 \le 6e^{-(\eta/2-s)Q_j}\le c_s.
\end{equation}
Given $x\in\partial\D$, consider the open Euclidean disk
\[
 U_j(x)=B((1-2r_j)x,r_j).
\]
Every $y\in U_j(x)$ satisfies
\[
 1-3r_j<|y|<1-r_j,
 \qquad |y-x|<3r_j.
\]
Hence $U_j(x)$ lies entirely in $\operatorname{int}\D\setminus\{0\}$,
while
\[
 \lambda(U_j(x))=r_j^2=4\varepsilon_j>\varepsilon_j.
\]
By \eqref{B:eq:free-disk-capacity} there is therefore some
$y\in U_j(x)$ with $H_j(y)\in U_j(x)$.  Thanks to the radius
condition \eqref{B:eq:boundary-radius-budget} we may apply
\eqref{B:eq:boundary-tube} with $Q=Q_j$ and $\sigma=\sigma_j$, which
gives
\[
 \begin{split}
 |H_j(x)-x|
 &\le |H_j(x)-H_j(y)|+|H_j(y)-x|\\
 &\le 3r_j(C_s e^{sQ_j}+1)\\
 &\le 6(C_s+1)e^{sQ_j}\sqrt{\varepsilon_j}.
 \end{split}
\]

The bound is uniform in $x$.  Since a bijection of $\partial\D$ and
its inverse have the same uniform displacement from the identity,
\[
 \sup_{x\in\partial\D}|F^{-Q_j}(x)-x|
 =\sup_{x\in\partial\D}|F^{Q_j}(x)-x|.
\]
The sign $\sigma_j$ therefore plays no role in the conclusion, and
\eqref{B:eq:boundary-rigidity-rate} follows with $K_{F,s}=6(C_s+1)$.

Suppose now that $\mathcal E(\alpha)>0$.
Choose a finite number
$0<\eta<\mathcal E(\alpha)$
and a subsequence $Q_j=q_{n_j}$
such that
$\log q_{n_j+1}\ge\eta Q_j$. The
convergent estimate
\[
 \normT{Q_j\alpha}\le q_{n_j+1}^{-1}\le e^{-\eta Q_j}
\]
puts us in the situation just treated, and boundary rigidity follows.

It remains to show that $g$ has
no wandering interval.
If $J\subset\T$ were such an
interval, choose $x\in J$.
Since the positive iterates of
$J$ are pairwise disjoint,
$g^m(x)\notin J$ for every $m\ge1$.
Consequently,
\[
d_{\T}(g^m(x),x)
\ge d_{\T}(x,\T\setminus J)>0
\qquad(m\ge1),
\]
contradicting the $C^0$-rigidity
just proved.

Thus $g$ has no wandering
interval. By the standard
classification of irrational
orientation-preserving circle
homeomorphisms, $g$ is
topologically conjugate to
$R_\alpha$.
Conversely, such a conjugacy gives $C^0$-rigidity along convergent
denominators, so the three circle conditions are equivalent.
\end{proof}

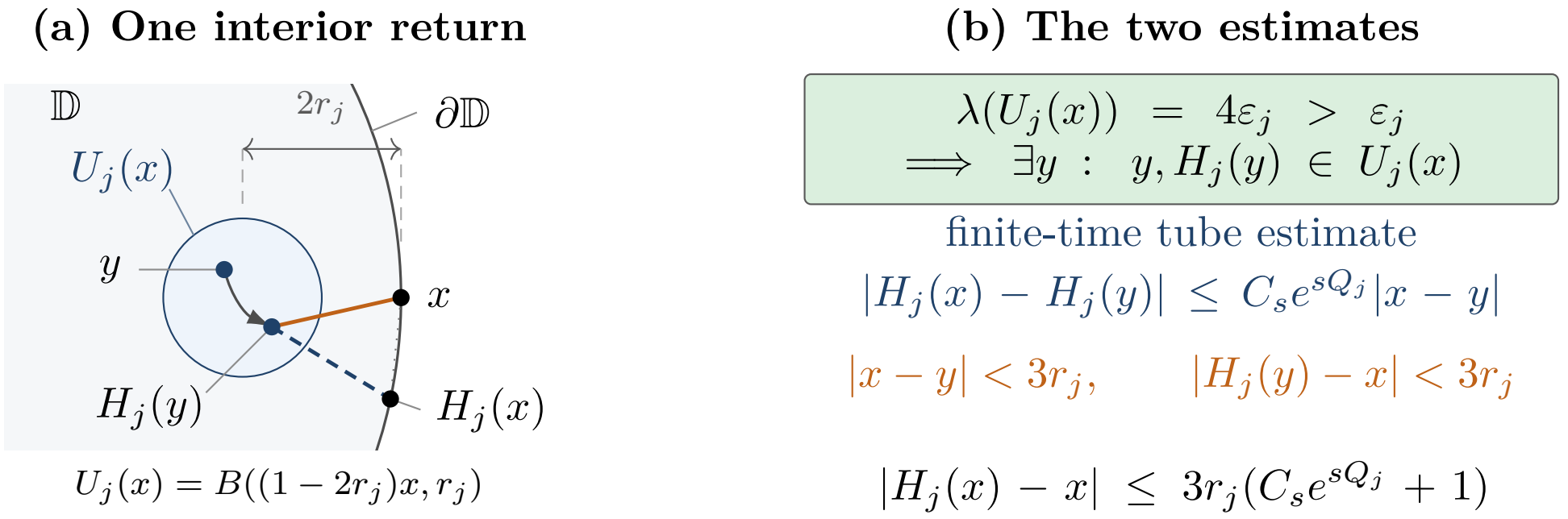
\begin{figure}[!htbp]
\centering
\begin{tikzpicture}[x=1cm,y=1cm,font=\small]
 \node[font=\bfseries\small] at (2.4,2.1)
   {(a) One interior return};
 \begin{scope}
  \clip (.15,-1.35) rectangle (4.6,1.65);
  \fill[softgray] (-.6,-.1) circle (4);
  \draw[nipl] (-.6,-.1) circle (4);
 \end{scope}
 \coordinate (bc) at (2.1,-.1);
 \coordinate (bx) at (3.4,-.1);
 \coordinate (by) at (1.95,.13);
 \coordinate (bhy) at (2.34,-.34);
 \coordinate (bhx) at ({-.6+4*cos(12)},{-.1-4*sin(12)});
 \filldraw[fill=softblue,draw=deepblue] (bc) circle (.65);
 \draw[frontline,line width=.9pt] (bx)--(bhy);
 \draw[deepblue,densely dashed,line width=.9pt] (bhy)--(bhx);
 \draw[black!55,dotted] (bx)--(bhx);
 \draw[nipa,bend right=25] (by) to (bhy);
 \draw[black!45,line width=.4pt] (1.25,.13)--(by);
 \draw[black!45,line width=.4pt] (1.8,-.88)--(bhy);
 \draw[black!45,line width=.4pt] (bhx)--(3.56,-1.02);
 \draw[deepblue!65,line width=.4pt] (1.5,.79)--(1.7,.412);
 \fill (bx) circle (2pt);
 \fill[deepblue] (by) circle (2pt);
 \fill[deepblue] (bhy) circle (2pt);
 \fill (bhx) circle (2pt);
 \node[anchor=west,inner sep=1pt] at (3.57,-.1) {$x$};
 \node[anchor=east,inner sep=1pt] at (1.15,.13) {$y$};
 \node[inner sep=1pt] at (1.34,-1.0) {$H_j(y)$};
 \node[anchor=west,inner sep=1pt] at (3.64,-1.02) {$H_j(x)$};
 \node[deepblue,inner sep=1pt] at (1.12,.94) {$U_j(x)$};
 \node at (.65,1.48) {$\D$};
 \node[anchor=west,inner sep=1pt] at (3.63,1.42) {$\partial\D$};
 \draw[black!45,line width=.4pt]
   ({-.6+4*cos(20)},{-.1+4*sin(20)})--(3.53,1.42);
 \draw[<->,black!65] (2.1,1.12)--(3.4,1.12)
   node[midway,above=2pt,font=\scriptsize] {$2r_j$};
 \draw[black!35,densely dashed] (2.1,.67)--(2.1,1.12);
 \draw[black!35,densely dashed] (3.4,.36)--(3.4,1.12);
 \node[align=center,font=\scriptsize] at (2.4,-1.65)
   {$U_j(x)=B((1-2r_j)x,r_j)$};
 \begin{scope}[xshift=6.9cm]
  \node[font=\bfseries\small] at (2.9,2.1)
    {(b) The two estimates};
  \node[nipbox,fill=coregreen,text width=5.9cm] (capacity) at (2.9,1.2)
    {$\lambda(U_j(x))=4\varepsilon_j>\varepsilon_j$\\
     $\Longrightarrow\ \exists y:\ y,H_j(y)\in U_j(x)$};
  \node[deepblue,align=center,text width=6.3cm] at (2.9,.15)
    {finite-time tube estimate\\[2pt]
     $|H_j(x)-H_j(y)|\le C_s e^{sQ_j}|x-y|$};
  \node[frontline,align=center] at (2.9,-.75)
    {$|x-y|<3r_j,\qquad |H_j(y)-x|<3r_j$};
  \node[align=center,text width=6.4cm] at (2.9,-1.65)
    {$\displaystyle |H_j(x)-x|
       \le 3r_j(C_s e^{sQ_j}+1)$};
 \end{scope}
\end{tikzpicture}
\caption{How an interior return controls a boundary point in the proof
of \cref{thm:boundary-rigidity}.  Here $r_j=2\sqrt{\varepsilon_j}$,
$H_j=F^{\sigma_jQ_j}$ has positive rotation representative
$\varepsilon_j$, and $U_j(x)$ lies entirely inside the disk and misses
$0$.  The free-disk estimate supplies one pair $y,H_j(y)\in U_j(x)$;
it does not make $U_j(x)$ invariant.  The dashed blue segment is
controlled by the finite-time tube estimate, and the orange segment
has length less than $3r_j$.  The drawing is schematic; no invariant
boundary collar is assumed.}
\label{B:fig:boundary-return}
\end{figure}

The geometric step in the boundary argument is illustrated in
\cref{B:fig:boundary-return}.

In \cref{thm:boundary-rigidity} the global condition
$\Per(F)=\{0\}$ enters through the free-disk estimate
\eqref{B:eq:free-disk-capacity}.  The boundary proof assumes neither
a rotation model near $0$ nor a conjugacy on the punctured interior.
Full-disk rigidity in the next proposition uses the additional
interior assumptions defining $\mathcal I_\alpha$.

\begin{proposition}
\label{B:prop:arithmetic-obstruction}
Let $F$ be a $C^1$ area-preserving pseudo-rotation of $\D$ with
irrational rotation number $\alpha$.  Suppose that $F=\mathcal R_\alpha$
near $0$ and that $F$ is smoothly conjugate to $\mathcal R_\alpha$ on
$\operatorname{int}(\D)\setminus\{0\}$.  Then
\begin{equation}\label{B:eq:uniform-subexp}
 \lim_{m\to\infty}\frac1m\log\|DF^m\|_\infty=0.
\end{equation}
Consequently, if $\mathcal E(\alpha)>0$,
then $F$ is $C^0$-rigid.  In particular, non-rigidity with this
interior structure requires zero exponential type.
\end{proposition}

\begin{proof}
\textit{1. Pointwise growth.}
Write the conjugacy on the punctured interior as
$F=\Psi^{-1}\mathcal R_\alpha\Psi$, so that for
$z\in\operatorname{int}\D\setminus\{0\}$
\[
 DF^m(z)=D\Psi^{-1}(\mathcal R_\alpha^m\Psi(z))\,
             D\mathcal R_\alpha^m\,D\Psi(z).
\]
Since the orbit of $\Psi(z)$ stays on a compact circle inside the
punctured disk, $\|DF^m(z)\|\le C(z)$ for every $m$.  At the origin
the same boundedness holds because $F$ is a rotation nearby, and on
the boundary \cref{B:lem:boundary-growth} applies.  Since the
determinant is one, every operator norm $\|DF^m(z)\|$ is at least
one, and putting the three regions together we get
\begin{equation}\label{B:eq:pointwise-subexp}
 \lim_{m\to\infty}\frac1m\log\|DF^m(z)\|=0
 \qquad(z\in\D).
\end{equation}

\textit{2. From pointwise to uniform growth.}
Set
\[
 a_m(z)=\log\|DF^m(z)\|,\qquad a_0=0,\qquad
 M=\log\max\{1,\|DF\|_\infty\}.
\]
The determinant-one condition and the chain rule give
\[
 0\le a_m(z)/m\le M\quad(m\ge1),\qquad
 a_{m+n}(z)\le a_m(F^nz)+a_n(z).
\]
Since the constant $M$ is integrable with respect to any
$F$-invariant probability measure $\mu$, dominated convergence lets
us exchange limit and integral:
\begin{equation}\label{B:eq:invariant-mean-subexp}
 \lim_{m\to\infty}\frac1m\int_\D a_m\,d\mu
 =\int_\D\lim_{m\to\infty}\frac{a_m(z)}m\,d\mu(z)=0.
\end{equation}
We now run the subadditive compactness argument that underlies
\cite[Proposition~1]{AvilaBochi}.
For fixed $k\ge1$ and $n\ge k$,
consider each offset
$r\in\{0,\ldots,k-1\}$ and write
\[
n=r+\ell_rk+t_r,
\qquad 0\le t_r<k.
\]
By subadditivity and
$0\le a_j\le jM$,
\[
a_n(z)\le\sum_{h=0}^{\ell_r-1}a_k(F^{r+hk}z)+(r+t_r)M\le\sum_{h=0}^{\ell_r-1}a_k(F^{r+hk}z)+2kM.
\]
Summing over all $k$ offsets,
every block starting index
$j=0,\ldots,n-k$ occurs
exactly once. Thus
\[
a_n(z)\le\frac1k\sum_{j=0}^{n-k}a_k(F^jz)+2kM\le\frac1k\sum_{j=0}^{n-1}a_k(F^jz)+2kM.
\]
where the last inequality
uses $a_k\ge0$.
Suppose uniform subexponential growth failed.  Then there would be
$n_i\to\infty$, $z_i\in\D$, and $\epsilon>0$ with
$a_{n_i}(z_i)\ge\epsilon n_i$, and along a subsequence the empirical
measures
\[
 \mu_i=\frac1{n_i}\sum_{j=0}^{n_i-1}\delta_{F^jz_i}
\]
would converge weakly to an $F$-invariant probability $\mu$.  By
continuity of $a_k$ and the block inequality, for every fixed $k$,
\[
 \epsilon\le\frac1k\int a_k\,d\mu.
\]
Letting $k\to\infty$ contradicts \eqref{B:eq:invariant-mean-subexp}.
This proves \eqref{B:eq:uniform-subexp}, and no minimality of the
dynamics was used.

\textit{3. Arithmetic returns.}
Suppose $\mathcal E(\alpha)>0$ and fix a finite number
$0<\eta<\mathcal E(\alpha)$, also when $\mathcal E(\alpha)=\infty$.
By the definition of the limit superior the index set
\[
 \mathcal N_\eta:=\{n:\log q_{n+1}\ge\eta q_n\}
\]
is infinite, and for its members $q_{n+1}^{-1/2}\le e^{-\eta q_n/2}$.
On the other hand, \eqref{B:eq:uniform-subexp} gives
\[
 \log\|DF^{q_n}\|_\infty=o(q_n),\qquad
 \|DF^{q_n}\|_\infty\le e^{\eta q_n/4}
 \quad\text{for all sufficiently large }n.
\]
The map $F^{q_n}$ is a pseudo-rotation with rotation number
$q_n\alpha\pmod\Z$, and $\|q_n\alpha\|_\T\le q_{n+1}^{-1}$, so
\cref{H:lem:displacement} yields
\begin{equation}\label{B:eq:displacement-subexp}
 \begin{aligned}
 \|F^{q_n}-\id\|_{C^0}
 &\le C_0(1+\|DF^{q_n}\|_\infty)\|q_n\alpha\|_\T^{1/2}\\
 &\le C_0(1+\|DF^{q_n}\|_\infty)q_{n+1}^{-1/2}.
 \end{aligned}
\end{equation}
Hence, for all sufficiently large $n\in\mathcal N_\eta$,
\[
 \|F^{q_n}-\id\|_{C^0}
 \le C_0(1+e^{\eta q_n/4})e^{-\eta q_n/2}
 \le2C_0e^{-\eta q_n/4}\longrightarrow0.
\]
As $q_n\to\infty$ along this infinite set, $F$ is $C^0$-rigid.
\end{proof}

\Needspace{5\baselineskip}
Together with the Brjuno examples of \cite{LWY},
\cref{thm:nonrigidity,B:prop:arithmetic-obstruction} give an exact
existence criterion within $\mathcal I_\alpha$.

\begin{corollary}
\label{cor:criterion}
For an irrational $\alpha$, the class $\mathcal I_\alpha$ contains
a non-$C^0$-rigid map if and only if $\mathcal E(\alpha)=0$.
\end{corollary}

\begin{proof}[Proof of \cref{cor:criterion}]
If $\mathcal E(\alpha)=0$ and $\alpha$ is non-Brjuno,
\cref{thm:nonrigidity} applies.
If $\alpha$ is Brjuno, we use
\cite[Theorem~1(i) and the
construction in \S7]{LWY},
which gives a non-$C^0$-rigid
disk pseudo-rotation belonging
to $\mathcal I_\alpha$. Together these two
cases cover every number of zero exponential type.  Conversely, when
$\mathcal E(\alpha)>0$, every member of $\mathcal I_\alpha$ is rigid
by \cref{B:prop:arithmetic-obstruction}.
\end{proof}

\subsection{An internal Denjoy circle}

\Needspace{5\baselineskip}
Reflecting the annulus construction across its Denjoy boundary moves the same mechanism to a smooth invariant circle in the interior of the disk.

\begin{proposition}
\label{B:prop:internal-denjoy-circle}
For every $\alpha$ satisfying the hypotheses of
\cref{thm:nonrigidity}, there is an area-preserving $C^1$ irrational
pseudo-rotation $\widetilde F$ of $\D$, with rotation number $\alpha$,
which equals $\mathcal R_\alpha$ near both $0$ and $\partial\D$, and
has a smooth embedded invariant circle
$C\subset\operatorname{int}\D\setminus\{0\}$ on which its
restriction is a Denjoy diffeomorphism.  The map $\widetilde F$ is
not $C^0$-rigid and is not $C^{1,\beta}$ for any $\beta>0$.
\end{proposition}

\begin{proof}
Let $A:\T\times[0,1]\to\T\times[0,1]$ be the annulus map from
\cref{B:prop:annulus-limit}, and write $g=A|_{s=0}$.  It has no
periodic points, and it equals the rigid rotation near $s=1$, by
\cref{B:prop:punctured}.
We first check the first jet needed
for the reflection.  For each fixed
$n$, the exact cotangent formula
\eqref{B:eq:cotangent-local}
holds on a sufficiently small
boundary collar, possibly depending
on $n$.  Thus
\[
F_n(x,s)=
\left(g_n(x),\frac{s}{g_n'(x)}\right)
\]
on that collar, and consequently
\[
DF_n(x,0)=
\begin{pmatrix}
	g_n'(x)&0\\
	0&1/g_n'(x)
\end{pmatrix}.
\]
The $C^1$ convergence in
\cref{B:prop:annulus-limit} gives
\begin{equation}\label{B:eq:reflection-boundary-jet}
	DA(x,0)=
	\begin{pmatrix}
		g'(x)&0\\
		0&1/g'(x)
	\end{pmatrix}.
\end{equation}

On $\T\times[-1,1]$ set
$J(x,s)=(x,-s)$ and define
\[
\widetilde A(x,s)=
\begin{cases}
	A(x,s),&s\ge0,\\
	JAJ(x,s),&s\le0.
\end{cases}
\]
The two traces agree on $s=0$.
Moreover, with
$P=\operatorname{diag}(1,-1)$,
\[
D(JAJ)(x,0^-)
=P\,DA(x,0^+)\,P
=DA(x,0^+).
\]
Thus the first jets match,
and $\widetilde A$ is a
$C^1$ diffeomorphism.
Since $J^*(dx\wedge ds)=-dx\wedge ds$, both pieces preserve
$dx\wedge ds$, as does their union.  The two half-annuli are invariant,
each is conjugate to $A$, and their common circle carries $g$;
therefore $\widetilde A$ has no periodic points.  Moreover,
$\widetilde A=\mathcal R_\alpha$ on full collars of $s=1$ and $s=-1$.

Use the action-coordinate embedding
\[
 \Phi(x,s)=\sqrt{\tfrac12-\tfrac14s}\,e^{2\pi i x}.
\]
Its image is the round annulus
$1/2\le |z|\le\sqrt3/2$, and
$\Phi^*(dX\wedge dY)=(\pi/4)dx\wedge ds$.
Define $\widetilde F=\Phi\widetilde A\Phi^{-1}$ there and extend
it by $\mathcal R_\alpha$ on the inner disk and the outer annulus.
The collar equalities make this a global $C^1$ standard-area-preserving
diffeomorphism of $\D$.

The irrational rotations on the complementary
regions and the absence of periodic points for $\widetilde A$ give
$\Per(\widetilde F)=\{0\}$.
The circle $C=\{|z|=1/\sqrt2\}$ is invariant and its restriction is
smoothly conjugate to $g$.  A wandering interval on $C$ prevents
uniform convergence of any sequence of positive iterates to the
identity, exactly as in \cref{B:lem:nonrigid}.  Finally, a
$C^{1,\beta}$ map restricted to this smooth circle would give a
$C^{1,\beta}$ circle map, contrary to
\cref{B:prop:sharp-increment}.
\end{proof}

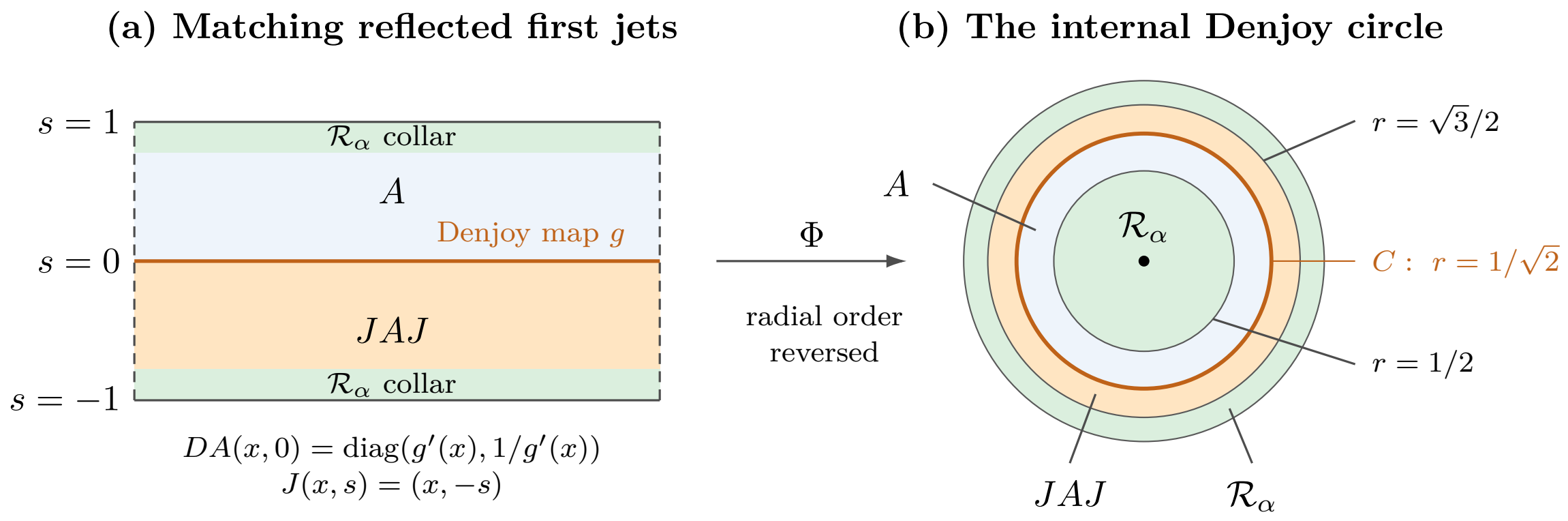
\begin{figure}[!htbp]
\centering
\begin{tikzpicture}[x=1cm,y=1cm,font=\small]
 \node[font=\bfseries\small] at (2.5,3.6)
   {(a) Matching reflected first jets};
 \fill[softblue] (0,1.35) rectangle (5.1,2.7);
 \fill[frontorange] (0,0) rectangle (5.1,1.35);
 \fill[coregreen] (0,2.4) rectangle (5.1,2.7);
 \fill[coregreen] (0,0) rectangle (5.1,.3);
 \draw[nipl] (0,0)--(5.1,0);
 \draw[nipl] (0,2.7)--(5.1,2.7);
 \draw[nipl,densely dashed] (0,0)--(0,2.7) (5.1,0)--(5.1,2.7);
 \draw[frontline,line width=1pt] (0,1.35)--(5.1,1.35);
 \node[left] at (0,2.7) {$s=1$};
 \node[left] at (0,1.35) {$s=0$};
 \node[left] at (0,0) {$s=-1$};
 \node at (2.5,2.03) {$A$};
 \node at (2.5,.68) {$JAJ$};
 \node[frontline,above,font=\scriptsize] at (3.85,1.35) {Denjoy map $g$};
 \node[font=\scriptsize] at (2.5,2.55) {$\mathcal R_\alpha$ collar};
 \node[font=\scriptsize] at (2.5,.15) {$\mathcal R_\alpha$ collar};
 \node[align=center,font=\scriptsize] at (2.5,-.65)
   {$DA(x,0)=\operatorname{diag}(g'(x),1/g'(x))$\\
    $J(x,s)=(x,-s)$};
 \draw[nipa] (5.65,1.35)--(7.5,1.35)
   node[midway,above] {$\Phi$};
 \node[align=center,font=\scriptsize] at (6.7,.65)
   {radial order\\reversed};
 \begin{scope}[xshift=9.8cm,yshift=1.35cm]
  \node[font=\bfseries\small] at (.25,2.25)
    {(b) The internal Denjoy circle};
  \filldraw[fill=coregreen,draw=black!65] (0,0) circle (1.75);
  \filldraw[fill=frontorange,draw=black!65] (0,0) circle (1.516);
  \fill[softblue] (0,0) circle (1.237);
  \draw[frontline,line width=1.1pt] (0,0) circle (1.237);
  \filldraw[fill=coregreen,draw=black!65] (0,0) circle (.875);
  \fill (0,0) circle (1.5pt);
  \node[above] at (0,.05) {$\mathcal R_\alpha$};
  \draw[nipl] (-1.05,.3)--(-2.05,.75);
  \node[anchor=east] at (-2.14,.75) {$A$};
  \draw[nipl] (-.47,-1.29)--(-.72,-1.96);
  \node[anchor=north] at (-.72,-2.0) {$JAJ$};
  \draw[nipl] (.83,-1.43)--(1.05,-1.96);
  \node[anchor=north] at (1.05,-2.0) {$\mathcal R_\alpha$};
  \draw[nipl] (1.162,.975)--(2.05,1.36);
  \node[anchor=west,font=\scriptsize] at (2.08,1.36) {$r=\sqrt3/2$};
  \draw[frontline] (1.237,0)--(2.05,0);
  \node[anchor=west,font=\scriptsize,text=frontline] at (2.08,0)
    {$C:\ r=1/\sqrt2$};
  \draw[nipl] (.67,-.562)--(2.05,-1.0);
  \node[anchor=west,font=\scriptsize] at (2.08,-1.0) {$r=1/2$};
 \end{scope}
\end{tikzpicture}
\caption{Moving the Denjoy mechanism to an internal circle
(\cref{B:prop:internal-denjoy-circle}).  The zero off-diagonal entries
in the displayed boundary jet make the reflected first jets agree;
the dashed vertical sides are identified.  Under
$\Phi(x,s)=\sqrt{1/2-s/4}\,e^{2\pi ix}$, the levels $s=1,0,-1$
go to radii $1/2,1/\sqrt2,\sqrt3/2$, respectively.
Thus $A$ occupies the inner half-annulus and $JAJ$ the outer one.
The rigid collars allow extension by the same irrational rotation
on the complementary disk and annulus, without adding periodic
points other than the center.  The circle $C$ carries $g$.
Radial thicknesses are schematic.}
\label{B:fig:internal-reflection}
\end{figure}

The reflection and the radial order of the embedding are shown in
\cref{B:fig:internal-reflection}.

\Needspace{5\baselineskip}
 We next explain why the Denjoy mechanism is also restricted when the map or the invariant circle has higher regularity.

\begin{remark}
\label{B:rem:internal-circle-obstructions}
The preceding relocation preserves zero exponential type.  At positive
exponential type, the proof of \cref{thm:boundary-rigidity} also rules
out Denjoy dynamics on every $C^1$ embedded invariant circle
$C\subset\operatorname{int}\D\setminus\{0\}$.  Here are the details
needed for this extension.  Such a circle surrounds $0$: otherwise
its invariant bounded complementary disk would contain an additional
fixed point by Brouwer's fixed-point theorem.  Its circle rotation
number is $\alpha$, by the common rotation number property of disk
pseudo-rotations \cite[\S1.1]{AFLXZ}.  The tangent--normal calculation
of \cref{B:lem:boundary-growth} gives uniform subexponential growth of
$DF^{\pm n}$ on $C$.  The fixed-block proof of
\cref{B:lem:boundary-tubes} then gives the same orbit-tube estimate
near $C$; take its width smaller than
$\dist(C,\partial\D)/4$ so that the joining segments remain in the
original round disk.
For $x\in C$ and the arithmetic sequence in the proof of the
theorem, the disk
$U_j(x)=B(x,2\sqrt{\varepsilon_j})$ lies in
$\operatorname{int}\D\setminus\{0\}$ for large $j$, uniformly in
$x$, and has normalized area $4\varepsilon_j$.
The same whole-disk free-disk estimate therefore gives
$y,H_j(y)\in U_j(x)$, where $H_j=F^{\pm Q_j}$ is chosen to have
rotation representative $\varepsilon_j$.  The orbit-tube bound yields
\[
 \sup_{x\in C}|H_j(x)-x|
 \le2(C_s+1)e^{sQ_j}\sqrt{\varepsilon_j}\longrightarrow0.
\]
The two signs have the same displacement on the invariant circle.
Thus $F|_C$ is rigid and cannot have a wandering interval.
For a $C^2$ ambient map, a Denjoy restriction to a $C^2$ embedded
invariant circle is already excluded by Denjoy's theorem.
These statements leave open geometrically rough invariant carriers
and other possible internal mechanisms of non-rigidity.
\end{remark}

\FloatBarrier
\section{Non-rigidity below the \texorpdfstring{$C^1$}{C1} endpoint}
\label{B:sec:bilip-positive-type}
Continuity of the derivative is essential in the boundary rigidity result. At finite positive exponential type, bi-Lipschitz regularity alone does not give the same conclusion.

\begin{proof}[Proof of \cref{thm:bilip-positive-type}]
We retain the local maps and protected squeezes of
\cref{B:sec:circle-construction,B:sec:hamiltonian-selection},
modifying the checkpoints and phase allocation as follows.

\Needspace{5\baselineskip}
\textit{1. Grouped local maps.}
Fix an integer $m\ge1$.  For each integer $r\ge K_{\min}$
divisible by $m$, put
\[
 \Psi_r(t)=G_{r,r+m}(mt),\qquad
 \widehat\Psi_r(t)=\mathcal G_{r,r+m}(mt),\qquad 0\le t\le1.
\]
Here $G$ is defined by \eqref{B:eq:block-composite}, and $\mathcal G$
is its Hamiltonian analogue preceding \eqref{B:eq:ham-moving-front}.
These paths go from the identity to the full block carrying $P_r$ to
$P_{r+m}$, and retain the protection implication
\eqref{B:eq:block-protect}.  Their quotients
$Q_{r;s,t}=\Psi_r(t)^{-1}\Psi_r(s)$ and
$\widehat Q_{r;s,t}=\widehat\Psi_r(t)^{-1}\widehat\Psi_r(s)$
are supported in $J_r$ and $Q_r$, respectively, preserve those sets,
and obey
\begin{equation}\label{B:eq:grouped-phase-bounds}
 \begin{aligned}
 \|\widehat Q_{r;s,t}-\id\|_0&\le C_m\delta_r|s-t|,\\
 \|D\widehat Q_{r;s,t}-I\|_0&\le C_m|s-t|,\\
 \|D^2\widehat Q_{r;s,t}\|_0&\le C_m\delta_r^{-1}|s-t|.
 \end{aligned}
\end{equation}
The analogous one-dimensional bounds hold, and the fixed-center
identity is
\begin{equation}\label{B:eq:grouped-center}
 Q_{r;s,t}(\theta_*)=\theta_*,\qquad
 Q_{r;s,t}'(\theta_*)=e^{m(t-s)}.
\end{equation}
Indeed all elementary factors are supported in nested rectangles
$Q_{r+j}\subset Q_r$, $0\le j<m$.  Rescaling both variables by
$\delta_r$, their spatial scales lie in $[2^{-(m-1)},1]$.
The smooth flow bounds of \cref{B:lem:ham-squeeze}, and their inverse
versions give uniform bounds, depending only on $m$, for the
displacements from the identity and for spatial derivatives through
order three.
On every clipped phase piece only one flow time varies.  Differentiating
in that time gives a uniform $C^2$ bound, and integration gives a
Lipschitz $C^2$ estimate across the finitely many phase endpoints as
well.  Composing with the inverse and undoing the rescaling proves
\eqref{B:eq:grouped-phase-bounds}.  Each elementary angular flow fixes
$\theta_*$ and has derivative $e^{-u}$ there; its clipped times sum
to $mt$.  This proves \eqref{B:eq:grouped-center}.

\Needspace{5\baselineskip}
\textit{2. Arithmetic and support selection.}
Set $E=\mathcal E(\alpha)$ and choose a fixed $m$ so large that
\begin{equation}\label{B:eq:grouped-arithmetic-choice}
 \frac{E}{m\log2}<\frac{\vartheta}{16},\qquad
 \vartheta=\frac1{10000}.
\end{equation}
Define $k_n$ by the formula \eqref{B:eq:checkpoints}, without imposing
zero exponential type, and replace the stage scales by
\[
 a_n=m\bigl(\lceil k_n/m\rceil+n\bigr),\qquad
 b_n=a_{n+1},\qquad
 w_n=(b_n-a_n)/m,\qquad e_n=w_n/q_n.
\]
Thus $a_n,b_n$ are multiples of $m$ and $w_n\ge1$.  Since $E<\infty$,
$\log q_{n+1}=O(q_n)$ eventually.  Consequently the logarithmic
corrections in $k_n$, divided by $q_n$, tend to zero, as do the
rounding errors and $n/q_n$.  Subtraction yields
\begin{equation}\label{B:eq:grouped-arithmetic-asymptotic}
 e_n-\frac{\log q_{n+1}}{m q_n\log2}\longrightarrow0,
 \qquad \limsup e_n=\frac{E}{m\log2}.
\end{equation}
The small-support estimate also survives:
\begin{equation}\label{B:eq:grouped-checkpoint-support}
	q_na_n^22^{-a_n}
	\le
	C_m^{\rm chk}
	(1+C_{\rm chk})^2
	2^{-C_{\rm chk}}
	\frac{2^{-n}}{(n+1)^2},
	\qquad
	a_n=O_m(\log q_n)\quad(n\to\infty).
\end{equation}
Here $C_m^{\rm chk}$ is
independent of $C_{\rm chk}$
and $n$.
Consequently, after fixing $m$,
the coefficient can be made
arbitrarily small by increasing
$C_{\rm chk}$.

For the first inequality use $a_n\ge k_n$, the boundedness of
$a_n/k_n$, and the direct checkpoint estimate
\eqref{B:eq:checkpoint-support}.  That estimate follows from the
formula for $k_n$ and does not require $e_n\to0$.

Keep the tower heights \eqref{B:eq:tower-height} and
$D_n=\lfloor\vartheta q_n\rfloor$.  By
\eqref{B:eq:grouped-arithmetic-choice}--\eqref{B:eq:grouped-arithmetic-asymptotic},
after increasing the starting index we have $D_n\ge2w_n$.
Distribute these $D_n$ edges as evenly as possible among the $w_n$
groups with initial scales
\[
 r=a_n,a_n+m,\ldots,b_n-m.
\]
If $d_{n,r}$ is the number assigned to one group, then
\begin{equation}\label{B:eq:grouped-allocation}
 \sum_r d_{n,r}=D_n,\qquad
 \frac{e_n}{2\vartheta}\le\frac1{d_{n,r}}
 \le\frac{4e_n}{\vartheta},\qquad \frac1{d_{n,r}}\le1.
\end{equation}
Use the sparse profile of \cref{B:def:sparse-front-profile} with total
height $w_n$, and jumps $1/d_{n,r}$ on the edges assigned to group
$r$.  On a tower level replace $G_{a_n,b_n}(\sigma_{n,j})$ by
$G_{a_n,b_n}(m\sigma_{n,j})$, and make the same replacement in the
Hamiltonian factors.  The old plateau is still the full original
block, so it carries $P_{a_n}$ onto $P_{b_n}$ exactly.

Let
\[
\mathscr R_n
=\{a_n+mh:h=0,\ldots,w_n-1\}
\]
be the set of group initial scales.
The selected edges on each front
are allocated to these groups in
outward order, as in
\cref{B:def:sparse-front-profile},
with $d_{n,r}$ edges assigned to
each group $r\in\mathscr R_n$.

The phase profile has total height
$w_n$, with the same plateau and
endpoint values and the same signed
jump convention as before.

For $r=a_n+mh$, the corresponding
group-phase interval is
\[
(w_n-h-1,w_n-h).
\]
On every edge assigned to group $r$,
both endpoint phases belong to its
closure.

Writing
\[
t_j=\sigma_{n,j}-(w_n-h-1),
\qquad
t_{j+1}
=\sigma_{n,j+1}-(w_n-h-1),
\]
we have $t_j,t_{j+1}\in[0,1]$ and
\[
|t_j-t_{j+1}|=d_{n,r}^{-1}.
\]

After replacing the original block
phase by $m\sigma_{n,j}$, the
completed finer factors cancel,
and the moving-front quotient is
exactly
\[
\begin{aligned}
	&\mathcal G_{a_n,b_n}
	(m\sigma_{n,j+1})^{-1}
	\mathcal G_{a_n,b_n}
	(m\sigma_{n,j})\\
	&\qquad=
	\widehat Q_{r;t_j,t_{j+1}}.
\end{aligned}
\]
Thus the grouped local estimates
apply to every selected edge.

The resulting quotient is supported in $Q_r$ and maps $Q_r$
onto itself. Consequently, after restoring the rotations in
\eqref{B:eq:ham-moving-front}, the map $\mathcal S_n$ sends
$\mathcal R_\alpha^jQ_r$ onto $\mathcal R_\alpha^{j+1}Q_r$.
Thus \eqref{B:eq:grouped-phase-bounds}
replaces the single-flow bounds
\eqref{B:eq:S2d0}--\eqref{B:eq:S2d2}.
Every original dyadic micro-scale is still used exactly once.  Hence
the upper derivative bound \eqref{B:eq:Bder}, the higher derivative
bounds \eqref{B:eq:h-high}, the protected-density induction of
\cref{B:prop:protected-density}, and the nested cotangent argument of
\cref{B:lem:nested-cotangent} remain valid.  In particular the
calculation proving \cref{B:prop:ledger}, now using the grouped bounds,
gives on an active group-$r$ physical component
\begin{equation}\label{B:eq:grouped-increments}
 \|\Delta_n\|_0\le C_m\delta_r/d_{n,r},\qquad
 \|D\Delta_n\|_0\le C_m/d_{n,r}.
\end{equation}
This is a finite-stage calculation and uses no convergence of
$1/d_{n,r}$ to zero.

The support-selection proof of \cref{B:prop:avoidance} uses only
$2D_n$ source arcs, protected density, the uniform upper bound on
$h_n'$, padded separation, and
$q_na_n^22^{-a_n}=O_m(2^{-n}(n+1)^{-2})$.  These hold by
\eqref{B:eq:grouped-checkpoint-support}; increasing the checkpoint
constant first and the starting index second supplies the same
packing margins.  It therefore selects pairwise disjoint closed
physical angular covers $\Omega_n$ with
\begin{equation}\label{B:eq:grouped-disjoint-supports}
 \supp\Delta_n\subset\Omega_n\times[0,1].
\end{equation}
\Needspace{5\baselineskip}
\textit{3. The bi-Lipschitz limit.}
In fixed physical coordinates at most one increment, and at most one
of its first derivatives, is nonzero at a given point.  Since
$F_{n_0}=\mathcal R_\alpha$, \eqref{B:eq:grouped-increments} and
\eqref{B:eq:grouped-disjoint-supports} imply
\[
 \sup_N\|DF_N\|_0\le1+C_m,\qquad
 \|F_N-F_M\|_0\le C_m\max_{M\le n<N}2^{-a_n}\longrightarrow0.
\]
The derivative bound is a maximum over disjoint physical supports,
not a product of successive Lipschitz constants.

Complete each stage to the same rotating central disk using the fixed
action coordinates described after \cref{B:prop:punctured}.  On the
annulus those coordinates and their inverses have bounded derivatives,
independent of the stage, and the maps are exactly rotations on the
gluing collar.  The completed smooth disk maps therefore have a common
Euclidean derivative bound.  They preserve standard area; in dimension
two the Euclidean operator norms of a determinant-one matrix and its
inverse agree.  On the convex disk both the maps and their inverses
consequently have a common Lipschitz constant.  Their uniform limit
$F$ is a bi-Lipschitz disk homeomorphism.  Moreover, for every continuous
$u:\D\to\R$,
\[
 \int_\D u\circ F\,dA
 =\lim_N\int_\D u\circ F_N\,dA=\int_\D u\,dA,
\]
so $F$ preserves standard area.

The punctured-interior proof in \cref{B:prop:punctured} requires only
that the conjugacies change in the physical collar
$s\le B_{\rm der}\delta_{a_n}$ and their inverses in the rotation
collar $v\le\delta_{a_n}$.  These inclusions still hold.  Both
sequences therefore stabilize smoothly on compact subsets of the
open annulus and yield mutually inverse smooth area conjugacies there.
After disk completion $F$ is smoothly conjugate by an area-preserving map to $\mathcal R_\alpha$
throughout $\operatorname{int}\D\setminus\{0\}$ and is the rotation
near $0$.

The interval proof in \cref{B:lem:stabilized} uses the unchanged full
plateau and hence also survives.  It supplies nondegenerate pairwise
disjoint boundary intervals $I_j$, $j\in\Z$, with
$g(I_j)=I_{j+1}$ for $g=F|_{\partial\D}$.
Since $g_N\to g$ uniformly and $\rot(g_N)=\alpha$, we have
$\rot(g)=\alpha$.  The irrational boundary rotation and the interior
conjugacy exclude every periodic point except $0$.
For $x\in\operatorname{int}I_0$ and every $n\ge1$,
\[
 |F^n(x)-x|\ge
 \dist\bigl(x,\partial\D\setminus\operatorname{int}I_0\bigr)>0,
\]
which gives the stated uniform lower bound and proves non-rigidity.

\Needspace{5\baselineskip}
\textit{4. Failure of $C^1$ regularity.}
Since $E>0$, the wandering boundary interval already excludes $C^1$
regularity by \cref{thm:boundary-rigidity}; the following calculation
gives a direct derivative-oscillation check.
At an active positive-front source center $x_n$, old support avoidance makes
$g_n=R_\alpha$ in a neighborhood.  If
$y_n=\theta_*+j\alpha$, then
$x_n=h_n^{-1}(y_n)$ and $z_n=h_n^{-1}(y_n+\alpha)$ satisfy
$h_n'(x_n)/h_n'(z_n)=g_n'(x_n)=1$.  Since the grouped quotient
fixes its center, these two conjugacy derivatives cancel in the
chain rule.  Formula \eqref{B:eq:grouped-center} therefore gives
\[
 g_{n+1}(x_n)=R_\alpha(x_n),\qquad
 g_{n+1}'(x_n)=e^{-m/d_{n,r}}.
\]
Choose $x_n^+$ in the interior
of the same physical source arc,
sufficiently close to its lateral
endpoint that the grouped quotient
is the identity on a neighborhood
of its rotation-coordinate image.
Then
\[
g_{n+1}'(x_n^+)=1.
\]
Since this entire source arc is
contained in the current physical
cover $\Omega_n$, all later
increments vanish on the arc.
Hence the limit $g$ agrees there
with the smooth map $g_{n+1}$,
and both derivative values are
retained. The physical source-arc
length is at most $C_ma_n^22^{-a_n}\to0$, by protected density and
\eqref{B:eq:grouped-checkpoint-support}.
On infinitely many stages \eqref{B:eq:grouped-arithmetic-asymptotic}
gives $e_n\ge E/(2m\log2)$, and
\eqref{B:eq:grouped-allocation} consequently gives
\[
 |g'(x_n)-g'(x_n^+)|
 \ge1-\exp\!\left(-\frac{E}{4\vartheta\log2}\right)>0.
\]
By compactness the paired points have a common accumulation point;
the displayed derivative oscillation contradicts continuity of $g'$.
Thus $g$, and hence $F$, is not $C^1$.
\end{proof}

\begin{figure}[htbp]
\centering
\begin{tikzpicture}[x=1cm,y=1cm,font=\small]
  \node[anchor=west,font=\bfseries\small] at (0,3.15)
    {(a) Disjoint stabilized source arcs};
  \draw[nipl] (0,2.0)--(6.1,2.0);
  \foreach \a/\b in {0.3/1.55,2.0/2.82,3.2/3.72,4.08/4.39}{
    \draw[line width=3pt,deepblue] (\a,2.0)--(\b,2.0);
  }
  \foreach \a/\b in {0.67/1.42,2.24/2.75,3.34/3.68,4.16/4.36}{
    \fill[frontline] (\a,2.0) circle (1.7pt);
    \fill[white] (\b,2.0) circle (2.0pt);
    \draw[deepblue] (\b,2.0) circle (2.0pt);
  }
  \node[below] at (.67,1.94) {$x_n$};
  \node[below] at (1.42,1.94) {$x_n^+$};
  \node[above] at (.93,2.05) {$J_n^{\rm phys}$};
  \node at (4.77,2.0) {$\cdots$};
  \draw[densely dashed] (5.35,1.6)--(5.35,2.5);
  \node[above] at (5.35,2.5) {$x_*$};
  \node[align=center,text width=6.1cm,font=\footnotesize] at (3.0,.68)
    {$|J_n^{\rm phys}|\to0$ along a subsequence;\\
     the limit is smooth in the interior of each displayed arc.};

  \begin{scope}[xshift=7.2cm]
    \node[anchor=west,font=\bfseries\small] at (0,3.15)
      {(b) Two retained derivative values};
    \draw[nipa] (.2,.48)--(.2,2.7);
    \draw[nipa] (.2,.48)--(5.75,.48);
    \node[anchor=east] at (.12,2.25) {$1$};
    \draw[gray,densely dashed] (.25,2.25)--(5.1,2.25);
    \fill[white] (4.2,2.25) circle (2.2pt);
    \draw[deepblue] (4.2,2.25) circle (2.2pt);
    \fill[frontline] (1.3,1.0) circle (2pt);
    \draw[densely dashed,gray] (1.3,.5)--(1.3,1.0);
    \draw[densely dashed,gray] (4.2,.5)--(4.2,2.2);
    \node[below] at (1.3,.48) {$x_n$};
    \node[below] at (4.2,.48) {$x_n^+$};
    \node[anchor=west,font=\scriptsize] at (1.45,.98)
      {$g'(x_n)=e^{-m/d_{n,r}}$};
    \node[above,font=\scriptsize] at (4.1,2.3) {$g'(x_n^+)=1$};
    \draw[{Latex[length=1.5mm]}-{Latex[length=1.5mm]}]
      (5.1,1.0)--(5.1,2.25);
    \node[anchor=west,font=\scriptsize]
      at (5.23,1.6) {$\ge c_E$};
    \node[align=center,text width=5.5cm,font=\scriptsize]
      at (2.95,-.3) {No derivative at $x_*$ is assumed.};
  \end{scope}
\end{tikzpicture}
\caption{The obstruction to $C^1$ regularity in
\cref{thm:bilip-positive-type}.  Here $J_n^{\rm phys}$ denotes one
selected physical source arc at stage $n$.  These arcs are pairwise
disjoint, their lengths tend to zero, and all subsequent increments
vanish there.  Thus the limit has the two indicated local
derivatives, whose gap is at least
$c_E=1-\exp(-E/(4\vartheta\log2))>0$ on an infinite subsequence,
where $E=\mathcal E(\alpha)$ and the group size $m$ is fixed.
The schematic does not assert differentiability at the common
accumulation point $x_*$.  Meanwhile
$\|\Delta_n\|_0\le C_m2^{-a_n}\to0$ and the disjoint supports give a
uniform derivative bound for the finite-stage maps.}
\label{B:fig:bilip-oscillation}
\end{figure}
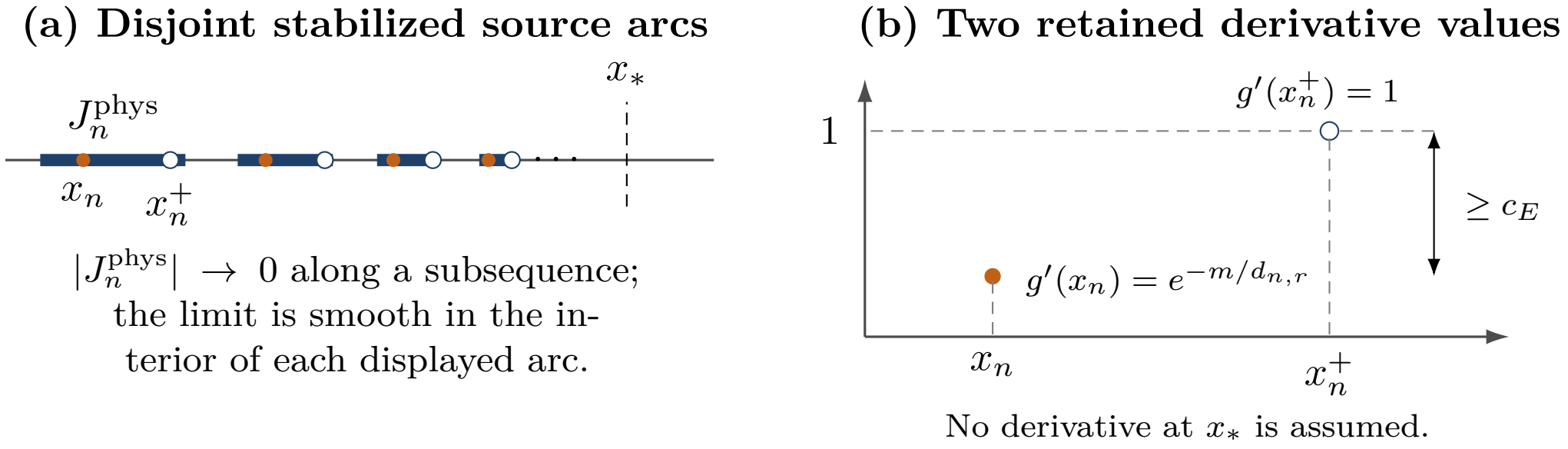

The separation between uniform convergence and derivative continuity
is illustrated in \cref{B:fig:bilip-oscillation}.

\Needspace{8\baselineskip}
\begin{remark}\label{B:rem:bilip-scope}
The conclusion of \cref{thm:bilip-positive-type} cannot extend to
$\mathcal E(\alpha)=\infty$: rigidity there is the Lipschitz case
$a=1$ of \cite[Theorem~2]{AFLXZ}.  For completeness, the free-disk argument
underlying \cref{H:lem:displacement} uses only the Lipschitz constant,
and gives
\[
 \|H-\id\|_0\le C(1+L)\sqrt{\|\rot(H)\|_{\T}}
\]
for a bi-Lipschitz conservative pseudo-rotation, where
$L=\max\{\Lip(H),\Lip(H^{-1})\}$.  To see this directly, set
$\varepsilon=\|\rot(H)\|_{\T}$ and assume
$r=2\sqrt\varepsilon<1/4$.  For any $z\in\D$ use the disk
$B((1-2r)z,r)$, which lies in $\operatorname{int}\D$ for small $r$;
if it contains $0$, move the center radially outward so that its
center has radius $2r$, taking any direction if $z=0$.
In either case its center is within $2r$
of $z$, its area is $\pi r^2>\pi\varepsilon$, and it misses $0$.
The free-disk capacity estimate forces a point $y$ in it with $H(y)$
in it, giving $|H(z)-z|\le3(1+L)r$; apply $H^{-1}$ for the negative
rotation representative.  When $2\sqrt\varepsilon\ge1/4$, the same
estimate follows, after enlarging $C$, from $\|H-\id\|_0\le2$.
If $F$ is bi-Lipschitz with constant $L_F$
in both time directions, use $H=F^{q_n}$, whose corresponding
constant is at most $L_F^{q_n}$.  When $\mathcal E(\alpha)=\infty$
a subsequence has $\log q_{n+1}/q_n\to\infty$, so
$\|F^{q_n}-\id\|_0\to0$ along that subsequence.
\end{remark}

\section*{Acknowledgments}
In preparing this work, we used ChatGPT (GPT-6) as an AI-assisted
research tool. We developed the overall proof strategy and communicated
it to the AI system, which assisted us in working out the technical
details. We subsequently checked the resulting arguments manually and
take full responsibility for the correctness of the proofs and for the
contents of this paper.

\Needspace{12\baselineskip}

\end{document}